\documentclass[11pt,a4paper,usenames,dvipsnames]{article}

\usepackage{url,amsmath,amssymb,latexsym,pstricks,mathrsfs,comment,amsthm,graphicx,tikz,tikz-cd,enumerate,accents,pgffor,cite,wrapfig,multicol,float,cases,calc,geometry}
\usepackage[colorlinks]{hyperref}

\usepackage{xcolor}

\usepackage{enumitem}

\usepackage[T1]{fontenc}
\usepackage{ wasysym }
\usepackage{stmaryrd}

\renewcommand{\arraystretch}{1.2}
\newcommand{\nc}{\newcommand}
\nc{\rnc}{\renewcommand}

\allowdisplaybreaks

\DeclareMathSymbol{\widehatsym}{\mathord}{largesymbols}{"62}
\newcommand\lowerwidehatsym{%
  \text{\smash{\raisebox{-1.3ex}{%
    $\widehatsym$}}}}
\newcommand\fixwidehat[1]{%
  \mathchoice
    {\accentset{\displaystyle\lowerwidehatsym}{#1}}
    {\accentset{\textstyle\lowerwidehatsym}{#1}}
    {\accentset{\scriptstyle\lowerwidehatsym}{#1}}
    {\accentset{\scriptscriptstyle\lowerwidehatsym}{#1}}
}

\nc{\gLa}{\mathrel{\mathscr L^a}}
\nc{\gRa}{\mathrel{\mathscr R^a}}
\nc{\gHa}{\mathrel{\mathscr H^a}}
\nc{\gJa}{\mathrel{\mathscr J^a}}
\nc{\gDa}{\mathrel{\mathscr D^a}}
\nc{\gKa}{\mathrel{\mathscr K^a}}
\nc{\gLh}{\mathrel{\fixwidehat{\mathscr L}}}
\nc{\gRh}{\mathrel{\fixwidehat{\mathscr R}}}
\nc{\gHh}{\mathrel{\fixwidehat{\mathscr H}}}
\nc{\gJh}{\mathrel{\fixwidehat{\mathscr J}}}
\nc{\gDh}{\mathrel{\fixwidehat{\mathscr D}}}
\nc{\gKh}{\mathrel{\fixwidehat{\mathscr K}}}
\nc{\Rh}{\fixwidehat{R}}
\nc{\Lh}{\fixwidehat{L}}
\nc{\Hh}{\fixwidehat{H}}
\nc{\Dh}{\fixwidehat{D}}
\nc{\Jh}{\fixwidehat{J}}
\nc{\Kh}{\fixwidehat{K}}

\newcommand{\darcxx}[4]{\draw[#4](#1,0)arc(180:90:#3) (#1+#3,#3)--(#2-#3,#3) (#2-#3,#3) arc(90:0:#3);}
\newcommand{\uarcxx}[4]{\draw[#4](#1,2)arc(180:270:#3) (#1+#3,2-#3)--(#2-#3,2-#3) (#2-#3,2-#3) arc(270:360:#3);}
\newcommand{\stlinex}[3]{\draw[#3](#1,2)--(#2,0);}
\newcommand{\alt}{\widetilde{\al}}
\newcommand{\stlines}[1]{{\foreach \x/\y in {#1} { \stline{\x}{\y} }}}

\nc{\lar}[1]{ \xrightarrow {\ #1\ }}
\nc{\relrank}[2]{\rank(#1\hspace{0.5truemm}{:}\hspace{0.5truemm}#2)}

\nc\E{\mathbb E}

\nc\Pmn{\P_{mn}}
\nc\Pnk{\P_{nk}}
\nc\Pmk{\P_{mk}}
\nc\Pkl{\P_{kl}}
\nc\Dem{\De_{[m]}}
\nc\Den{\De_{[n]}}
\nc\RR{\mathbb R}
\nc\TL{\mathcal{T}\!\mathcal{L}}
\nc\KK{\mathcal K}
\nc\Mod[1]{\ (\textup{mod}\ #1)}
\nc\Kmn{\KK_{mn}}
\nc\Bmn{\B_{mn}}
\nc\Kmns{\KK_{mn}^\si}
\nc\Pmns{\P_{mn}^\si}
\nc\PBmns{\PB_{mn}^\si}
\nc\Bmns{\B_{mn}^\si}
\nc\Bklt{\B_{kl}^\tau}
\nc\Mmns{\M_{mn}^\si}
\nc\PPmns{\PP_{mn}^\si}
\nc\TLmns{\TL_{mn}^\si}
\nc\Sone{S^{(1)}}
\nc\Tone{T^{(1)}}
\nc\ve\varepsilon
\nc\leqJa{\leq_{\J^a}}
\nc\leqJsi{\leq_{\J^\si}}

\nc\Pre{\operatorname{Pre}}
\nc\Post{\operatorname{Post}}
\nc\Reg{\operatorname{Reg}}
\nc\RI{\operatorname{RI}}
\nc\MI{\operatorname{MI}}
\nc\RP{\operatorname{RP}}
\nc\even{{\operatorname{even}}}
\nc{\Set}{\operatorname{{\bf Set}}}

\nc\bd{{\bf d}}
\nc\br{{\bf r}}

\nc{\set}[2]{\{#1:#2\}}
\nc{\bigset}[2]{\big\{#1:#2\big\}}
\nc\sub\subseteq
\nc\mt\mapsto
\nc\al\alpha
\nc\be\beta
\nc\ga\gamma
\nc\de\delta
\nc\si\sigma
\nc\ep\epsilon
\nc\lam\lambda
\nc\ze\zeta
\nc\vk\varkappa
\nc\vt{\widetilde{\nu}}
\nc\vs\varsigma
\nc\ka\kappa
\rnc\th\theta
\nc\om\omega
\nc\ups\upsilon
\nc\Ga\Gamma
\nc\De\Delta
\nc\Si\Sigma
\nc\Om\Omega
\rnc\O{\mathbb O}
\nc\bit{\begin{itemize}}
\nc\eit{\end{itemize}}
\nc\ben{\begin{enumerate}[label=\textup{(\roman*)},leftmargin=7mm]}
\nc\bena{\begin{enumerate}[label=\textup{(\alph*)},leftmargin=7mm]}
\nc\een{\end{enumerate}}
\nc\bmc{\begin{multicols}}
\nc\emc{\end{multicols}}
\nc{\leqH}{\leq_{\H}}
\nc{\leqR}{\leq_{\R}}
\nc{\leqL}{\leq_{\L}}
\nc{\leqJ}{\leq_{\J}}
\nc{\leqK}{\leq_{\K}}
\nc{\geqK}{\geq_{\K}}
\nc{\geqR}{\geq_{\R}}
\nc{\geqL}{\geq_{\L}}
\nc{\geqJ}{\geq_{\J}}

\nc\bp{{\bf p}}
\nc\bq{{\bf q}}
\rnc\iff{\ \Leftrightarrow\ }
\rnc\implies{\ \Rightarrow\ }

\nc\pf{\begin{proof}}
\nc\epf{\end{proof}}
\nc\epfres{\hfill\qed}
\nc\epfreseq{\tag*{\qed}}

\let\oldproofname=\proofname
\renewcommand{\proofname}{\rm\bf{\oldproofname}}

\nc\AND{\qquad\text{and}\qquad}
\nc\OR{\qquad\text{or}\qquad}
\nc\WHERE{\qquad\text{where}\qquad}
\nc\ANd{\quad\text{and}\quad}
\nc\anD{\ \ \ \text{and}\ \ \ }
\nc\ANDSIM{\qquad\text{and similarly}\qquad}
\nc{\COMMA}{,\qquad}
\nc{\COMMa}{,\quad}
\nc\ul\underline
\nc\ol\overline
\nc\permdec[1]{#1^{\natural}}
\nc\ext[1]{#1^\textup{E}}
\nc\wh\widehat
\nc\normal\unlhd
\rnc\emptyset\varnothing
\nc{\firstpfitem}[1]{#1.}
\nc{\pfitem}[1]{\medskip \noindent #1.}
\nc{\pfcase}[1]{\medskip\noindent {\bf Case #1.}}
\nc{\pfsubcase}[1]{\medskip\noindent {\bf Subcase #1.}}
\nc\im{\operatorname{im}}
\nc\LSUB{\operatorname{LSUB}}

\nc\pre\preceq
\nc\Y{\mathcal Y}
\nc\B{\mathcal B}
\nc\Z{\mathbb Z}
\nc\F{\mathfrak F}
\nc\T{\mathcal T}
\nc\TT{\mathscr T}
\nc\PP{\mathscr P\mathcal P}
\nc\C{\mathscr C}
\nc\I{\mathcal I}
\nc\Eq{\mathbb{E}}
\nc\Part{\mathbb{P}}
\nc\cg[2]{(#1,#2)^\sharp}
\nc\Rev{\operatorname{Rev}}
\nc\cR{\mathcal R}
\nc\Ptop{P^\top}
\nc\Qtop{Q^\top}
\nc\la\langle
\nc\ra\rangle
\nc\tb[1]{\operatorname{Seq}(#1)}
\nc\RevX{{\Rev}\big([0,|X|],[0,|X|^+]\big)}

\makeatletter
\DeclareRobustCommand\widecheck[1]{{\mathpalette\@widecheck{#1}}}
\def\@widecheck#1#2{%
    \setbox\z@\hbox{\m@th$#1#2$}%
    \setbox\tw@\hbox{\m@th$#1%
       \widehat{%
          \vrule\@width\z@\@height\ht\z@
          \vrule\@height\z@\@width\wd\z@}$}%
    \dp\tw@-\ht\z@
    \@tempdima\ht\z@ \advance\@tempdima2\ht\tw@ \divide\@tempdima\thr@@
    \setbox\tw@\hbox{%
       \raise\@tempdima\hbox{\scalebox{1}[-1]{\lower\@tempdima\box
\tw@}}}%
    {\ooalign{\box\tw@ \cr \box\z@}}}
\makeatother

\nc\ch\widecheck

\newcommand{\uv}[1]{\fill (#1,2)circle(.17);}
\newcommand{\lv}[1]{\fill (#1,0)circle(.17);}
\newcommand{\uvs}[1]{{\foreach \x in {#1} { \uv{\x}}}}
\newcommand{\lvs}[1]{{\foreach \x in {#1} { \lv{\x}}}}
\newcommand{\darcx}[3]{\draw(#1,0)arc(180:90:#3) (#1+#3,#3)--(#2-#3,#3) (#2-#3,#3) arc(90:0:#3);}
\newcommand{\darc}[2]{\darcx{#1}{#2}{.4}}
\newcommand{\uarcx}[3]{\draw(#1,2)arc(180:270:#3) (#1+#3,2-#3)--(#2-#3,2-#3) (#2-#3,2-#3) arc(270:360:#3);}
\newcommand{\uarc}[2]{\uarcx{#1}{#2}{.4}}
\newcommand{\stline}[2]{\draw(#1,2)--(#2,0);}

\newcommand{\uvx}[2]{\fill (#1,2)circle(#2);}
\newcommand{\lvx}[2]{\fill (#1,0)circle(#2);}
\newcommand{\uvxs}[2]{{\foreach \x in {#1} { \uvx{\x}{#2}}}}
\newcommand{\lvxs}[2]{{\foreach \x in {#1} { \lvx{\x}{#2}}}}

\nc{\buv}[1]{\fill (#1,2)circle(.18);}
\nc{\buvs}[1]{{
\foreach \x in {#1}
{ \buv{\x}}
}}
\nc{\blv}[1]{\fill (#1,0)circle(.18);}
\nc{\blvs}[1]{{
\foreach \x in {#1}
{ \blv{\x}}
}}

\nc{\uarcs}[1]{
{\foreach \x/\y in {#1}
{ \uarc{\x}{\y} }
}
}

\nc{\darcs}[1]{
{\foreach \x/\y in {#1}
{ \darc{\x}{\y} }
}
}
\nc{\darcxhalf}[3]{\draw(#1,0)arc(180:90:#3) (#1+#3,#3)--(#2,#3) ;}
\nc{\darchalf}[2]{\darcxhalf{#1}{#2}{.4}}
\nc{\uarcxhalf}[3]{\draw(#1,2)arc(180:270:#3) (#1+#3,1.5-#3)--(#2,1.5-#3) ;}
\nc{\uarchalf}[2]{\uarcxhalf{#1}{#2}{.4}}
\nc{\colv}[3]{\fill[#3] (#1,#2)circle(.17);}

\nc{\uvert}[1]{\fill (#1,2)circle(.2);}
\rnc{\lvert}[1]{\fill (#1,0)circle(.2);}

\nc{\custpartn}[3]{{\lower1.4 ex\hbox{
\begin{tikzpicture}[scale=.3]
\foreach \x in {#1}
{ \uvert{\x}  }
\foreach \x in {#2}
{ \lvert{\x}  }
#3 \end{tikzpicture}
}}}

\renewcommand{\P}{\mathcal P} 
\newcommand{\PB}{\mathcal{PB}} 
 
\renewcommand{\S}{\mathcal{S}}
\newcommand{\M}{\mathcal M}

\nc\MYZ{\mathcal M_{Y\cup Z}}

\renewcommand{\H}{\mathrel{\mathscr H}}
\renewcommand{\L}{\mathrel{\mathscr L}}
\newcommand{\R}{\mathrel{\mathscr R}}
\newcommand{\D}{\mathrel{\mathscr D}}
\newcommand{\J}{\mathrel{\mathscr J}}
\newcommand{\K}{\mathrel{\mathscr K}}

\newcommand{\N}{\mathbb{N}}

\nc\HH{\mathcal H}

\newcommand{\coker}{\operatorname{coker}}
\newcommand{\dom}{\operatorname{dom}} 
\newcommand{\codom}{\operatorname{codom}}
\newcommand{\rank}{\operatorname{rank}}
\newcommand{\idrank}{\operatorname{idrank}}

\newcommand{\id}{\operatorname{id}}

\renewcommand{\c}{@{}c@{}}
\newcommand{\cend}{@{}c@{\hspace{1.5truemm}}}

\newcommand{\partI}[8]{
\Big( 
{ \scriptsize \renewcommand*{\arraystretch}{1}
\begin{array} {\c|\c|\c|\c|\c|\cend}
 #1 \:&\: \cdots \:&\: #2 \:&\: #3 \:&\: \cdots \:&\: #4 \\ \cline{4-6}
 #5 \:&\: \cdots \:&\: #6 \:&\: #7 \:&\: \cdots \:&\: #8 
\rule[0mm]{0mm}{2.7mm}
\end{array} 
}
\hspace{-1.5 truemm} \Big)
}

\newcommand{\partpermIII}[6]{
\Big(
{ \scriptsize \renewcommand*{\arraystretch}{1}
\begin{array} {\c|\c|\c}
 #1 \:&\: #2 \:&\: #3 \\ 
 #4 \:&\: #5 \:&\: #6 
\rule[0mm]{0mm}{2.7mm}
\end{array} 
}
\Big)
}

\newcommand{\partpermV}[8]{
\Big(
{ \scriptsize \renewcommand*{\arraystretch}{1}
\begin{array} {\c|\c|\c|\c|\c}
 #1 \:&\: \cdots \:&\: #2 \:&\: #3 \:&\: #4 \\ 
 #5 \:&\: \cdots \:&\: #6 \:&\: #7 \:&\: #8 
\rule[0mm]{0mm}{2.7mm}
\end{array} 
}
\Big)
}
\newcommand{\partpermVI}[8]{
\Big(
{ \scriptsize \renewcommand*{\arraystretch}{1}
\begin{array} {\c|\c|\c|\c|\c|\c}
 #1 \:&\: \cdots \:&\: #2 \:&\: #3 \:&\: \cdots \:&\: #4 \\ 
 #5 \:&\: \cdots \:&\: #6 \:&\: #7 \:&\: \cdots \:&\: #8 
\rule[0mm]{0mm}{2.7mm}
\end{array} 
}
\Big)
}

\numberwithin{equation}{section}

\newtheorem{thm}[equation]{Theorem}
\newtheorem{lemma}[equation]{Lemma}
\newtheorem{cor}[equation]{Corollary}
\newtheorem{prop}[equation]{Proposition}

\theoremstyle{definition}

\newtheorem{rem}[equation]{Remark}

\newtheorem{eg}[equation]{Example}

\newcommand{\sm}{\setminus}

\begin{document}

\title{Sandwich semigroups in diagram categories}
\author{
Ivana \DJ ur\dj ev\footnote{Department of Mathematics and Informatics, University of Novi Sad, Novi Sad, Serbia. Mathematical Institute of the Serbian Academy of Sciences and Arts, Beograd, Serbia. {\it Emails:} {\tt ivana.djurdjev\,@\,dmi.uns.ac.rs}, {\tt ivana.djurdjev\,@\,mi.sanu.ac.rs}}, 
Igor Dolinka\footnote{Department of Mathematics and Informatics, University of Novi Sad, Novi Sad, Serbia. {\it Email:} {\tt dockie\,@\,dmi.uns.ac.rs}},
James East\footnote{Centre for Research in Mathematics, School of Computing, Engineering and Mathematics, Western Sydney University, Sydney, Australia. {\it Email:} {\tt j.east\,@\,westernsydney.edu.au}}
}
\date{}

\maketitle

\begin{abstract}
This paper concerns a number of diagram categories, namely the partition, planar partition, Brauer, partial Brauer, Motzkin and Temperley-Lieb categories.  If $\mathcal K$ denotes any of these categories, and if $\sigma\in\mathcal K_{nm}$ is a fixed morphism, then an associative operation $\star_\sigma$ may be defined on~$\mathcal K_{mn}$ by $\alpha\star_\sigma\beta=\alpha\sigma\beta$.  The resulting semigroup $\mathcal K_{mn}^\sigma=(\mathcal K_{mn},\star_\sigma)$ is called a sandwich semigroup.  We conduct a thorough investigation of these sandwich semigroups, with an emphasis on structural and combinatorial properties such as Green's relations and preorders, regularity, stability, mid-identities, ideal structure, (products of) idempotents, and minimal generation.  It turns out that the Brauer category has many remarkable properties not shared by any of the other diagram categories we study.  Because of these unique properties, we may completely classify isomorphism classes of sandwich semigroups in the Brauer category, calculate the rank (smallest size of a generating set) of an arbitrary sandwich semigroup, enumerate Green's classes and idempotents, and calculate ranks (and idempotent ranks, where appropriate) of the regular subsemigroup and its ideals, as well as the idempotent-generated subsemigroup.  Several illustrative examples are considered throughout, partly to demonstrate the sometimes-subtle differences between the various diagram~categories.


\end{abstract}

\tableofcontents

\section{Introduction}\label{s:intro}

Categories and algebras of diagrams (which are visual representations of set partitions) arise in many branches of mathematics, including representation theory \cite{HR2005,Martin2008}, statistical mechanics \cite{Martin1994,Jones1994_2,TL1971,Kauffman1987}, knot theory \cite{Stroppel2005,Jones1987,Jones1983_2,Kauffman1990,Kauffman1997}, classical groups \cite{Brauer1937}, invariant theory \cite{LZ2015,LZ2012}, and more.  An excellent overview may be found in the survey \cite{Martin2008}, including a detailed discussion of the role played by diagram categories in theoretical physics.  

Sandwich semigroups were used in \cite{DE2018,Sandwich1,Sandwich2} to give algebraic structures to arbitrary hom-sets in (locally small) categories, and indeed in more general structures called partial semigroups.  If $\C$ is a category, and if $X$ and $Y$ are fixed objects of $\C$, then elements of the hom-set $\C_{XY}$ can only be directly composed if $X=Y$, in which case $\C_{XY}=\C_X$ is an endomorphism monoid of $\C$.  However, if we fix a morphism $a\in\C_{YX}$, then an associative operation $\star_a$ can be defined on $\C_{XY}$ by $f\star_ag=fag$, and we obtain a \emph{sandwich semigroup} $\C_{XY}^a=(\C_{XY},\star_a)$.
This construction generalises many families of examples previously only studied as isolated special cases \cite{Brown1955,Munn1955,Lyapin,MS1975,Thornton1982,Chase1979}.  For more background and references, see the introductions to \cite{DE2018,Sandwich1,DE2015}.

A general theory of sandwich semigroups in (locally small) categories has been developed in \cite{DE2018,Sandwich1}, and this has been applied to several concrete categories of (linear) transformations in \cite{DE2018,Sandwich2}.  The current article continues to develop the general theory, and moves towards applications to a number of diagram categories, namely the partition category $\P$ \cite{Martin1994,Jones1994_2}, the planar partition category~$\PP$~\cite{HR2005,Jones1994_2}, the Brauer category~$\B$ \cite{LZ2015,Brauer1937}, the partial Brauer category $\PB$ \cite{MarMaz2014,Maz1998}, the Temperley-Lieb category~$\TL$~\cite{Stroppel2005,TL1971} and the Motzkin category~$\M$ \cite{BH2014}.  While sandwich semigroups in these diagram categories share some common properties with those in (linear) transformation categories, there are some striking differences that will be emphasised during the article.  There are also some intriguing dichotomies between different diagram categories.  For example, the behaviour of maximal elements in the division order on sandwich semigroups in the planar categories~$\PP$,~$\TL$ and~$\M$ is quite different to that in their nonplanar counterparts $\P$, $\B$ and~$\PB$.

Perhaps most surprising of all is the long list of neat structural and combinatorial features possessed by the Brauer category, but by none of the other diagram categories studied here.  An entire section of the paper is devoted to this exceptional category.  The kinds of problems we solve mostly involve the computation of combinatorial invariants, related to Green's classes, minimal sizes of generating sets, isomorphism classes of sandwich semigroups, and so on.  These take their inspiration from influential works of Howie and others on combinatorial (transformation) semigroup theory; see for example \cite{Howie1966,Howie1978,GH1992,HM1990,Gray2007,Gray2008,Gray2014}, and in particular \cite{DEG2017,EG2017} for related studies of diagram monoids.  The introduction to \cite{EG2017} contains many more references.

The article is organised as follows.  Section \ref{s:S} contains preliminary material on (partial) semigroups and categories.  Section \ref{s:sandwich} provides some general framework for working with sandwich semigroups; as well as revising some of the key results from \cite{DE2018,Sandwich1}, we also develop a theory of Green's preorders in arbitrary sandwich semigroups, with a focus on maximal classes, and connections to mid-identities and one-sided invertibility.  Section~\ref{s:DC} introduces the diagram categories that will be the focus of our study, and proves a number of structural and combinatorial results concerning them.  Section \ref{s:sandwichK} proves many results on sandwich semigroups in diagram categories, including characterisations of the regular elements, determination of Green's relations and preorders, classification of maximal classes, description of the idempotent-generated subsemigroup, and criteria for idempotent-generation of certain ideals.  Section \ref{s:B} exclusively concerns the Brauer category $\B$.  As noted above, $\B$ has many special properties not shared by any of the other diagram categories.  These allow us to solve several additional problems for~$\B$; in particular, we completely classify isomorphism classes of sandwich semigroups, calculate the rank (minimum size of a generating set) of an arbitrary sandwich semigroup, enumerate Green's classes and idempotents, and calculate ranks (and idempotent ranks, where appropriate) of the regular subsemigroup and its ideals, as well as the idempotent-generated subsemigroup.  We also explain why many of these results do not hold in the other diagram categories.

\section{\boldmath Preliminaries on (partial) semigroups and regular $*$-categories}\label{s:S}

\subsection{Basic definitions}\label{ss:basic}

We begin by recalling some ideas from \cite{DE2018,Sandwich1}), slightly adapting notation to suit our present purposes.  A \emph{partial semigroup} is a 5-tuple $(S,I,\bd,\br,\cdot)$, where $S$ and $I$ are sets, $\bd,\br:S\to I$ are mappings, and $(x,y)\mt x\cdot y$ is a partial binary operation (defined on a subset of $S\times S$), such that for all $x,y,z\in S$:
\ben
\item \label{C1} 
$x\cdot y$ is defined if and only if $\br(x)=\bd(y)$, in which case $\bd(x\cdot y)=\bd(x)$ and $\br(x\cdot y)=\br(y)$,
\item \label{C2} 
if $x\cdot y$ and $y\cdot z$ are defined, then $(x\cdot y)\cdot z=x\cdot (y\cdot z)$.
\een
If the context is clear, we usually write $xy$ for a well-defined product $x\cdot y$.  If the set~$I$, the mappings~$\bd,\br$, and the product $\cdot$ are all understood, we will often refer to ``the partial semigroup~$S$'' instead of ``the partial semigroup $(S,I,\bd,\br,\cdot)$''.  By \ref{C2}, we may omit parentheses on products of length greater than two.  
If $|I|=1$, then of course $S$ is just a semigroup.
For $i,j\in I$, we define the set
\[
S_{ij}=\set{x\in S}{\bd(x)=i,\ \br(x)=j}.
\]
The sets $S_i=S_{ii}$ are semigroups.

As in \cite[page 21]{MacLane1998}, an element $x\in S$ is \emph{(von Neumann) regular} if there exists $a\in S$ such that $x=xax$; note then that $\bd(a)=\br(x)$ and $\br(a)=\bd(x)$.  
Following standard semigroup terminology, we then call $a$ a \emph{pre-inverse} of $x$, and $x$ a \emph{post-inverse} of $a$.  We write $\Pre(x)$ and $\Post(x)$ for the sets of all pre- or post-inverses of $x$, respectively.  Note that $x\in\Post(a)\iff a\in\Pre(x)$.  As usual, the elements of $\Pre(x)\cap\Post(x)$ are called \emph{inverses} of $x$; we denote this set by $V(x)$.  If $a\in\Pre(x)$, then $axa\in V(x)$, so $\Pre(x)\not=\emptyset\iff V(x)\not=\emptyset$.
It is possible to have $\Pre(x)=\emptyset\not=\Post(x)$.
If every element of $S$ is regular, we say $S$ is (von Neumann) regular.  (Note that the term ``regular'' has another meaning within category theory; here we always mean \emph{von Neumann} regular.)

We call the partial semigroup $S$ a \emph{(small) category} if it additionally satisfies the following:
\ben \setcounter{enumi}{2}
\item \label{C3} 
for all $i\in I$ there exists $e_i\in S_i$ such that for all $x\in S$, $xe_{\br(x)}=e_{\bd(x)}x=x$.
\een
In the usual meaning of the word ``category'', 
\bit
\item the objects of $S\equiv(S,I,\bd,\br,\cdot)$ are the elements of $I$,
\item the hom-sets (a.k.a.~morphism sets) are the $S_{ij}$ ($i,j\in I$),
\item the semigroups $S_i=S_{ii}$ ($i\in I$) are monoids; the elements of $S_i$ are endomorphisms, and invertible elements of $S_i$ are automorphisms,
\item $\bd(x)$ and $\br(x)$ are the domain and range (or source and target) of the morphism $x\in S$,
\item $\cdot$ is the composition operation (and morphisms are composed left to right).
\eit
Even if a partial semigroup $S$ is not a category, we will still usually refer to the $S_{ij}$ as hom-sets, the~$S_i$ as endomorphism semigroups, and so on.

By a \emph{partial $*$-semigroup}, we mean a 6-tuple $(S,I,\bd,\br,\cdot,*)$ such that $(S,I,\bd,\br,\cdot)$ is a partial semigroup, and $*:S\to S:x\mt x^*$ a mapping such that for all $x,y\in S$,
\ben \setcounter{enumi}{3}
\item \label{SC1} $\bd(x^*)=\br(x)$, $\br(x^*)=\bd(x)$ and $(x^*)^*=x$,
\item \label{SC2} if $xy$ is defined, then $(xy)^*=y^*x^*$.
\een
Analogously to \cite{NS1978}, by a \emph{regular partial $*$-semigroup}, we mean a partial $*$-semigroup $S$ such that 
\ben \setcounter{enumi}{5}
\item \label{SC3}  $x=xx^*x$ for all $x\in S$.
\een
In particular, any regular partial $*$-semigroup is regular.  If a regular partial $*$-semigroup happens to be a category, we call it a \emph{regular $*$-category}.  All of our motivating examples in Sections \ref{s:DC}--\ref{s:B} are regular $*$-categories.  It follows from \ref{SC3} and \ref{SC1} that $x^* = x^*(x^*)^*x^* = x^*xx^*$ for all $x\in S$.

\subsection{Green's relations and preorders}\label{ss:Green}

As in \cite{DE2018}, for a partial semigroup $S\equiv(S,I,\bd,\br,\cdot)$, we denote by $\Sone$ the category obtained by adjoining an identity morphism at each object that did not already have one.  Green's preorders on $S$ are defined, for $x,y\in S$, by
\bit
\item $x\leqR y \iff x=ya$ for some $a\in \Sone$,
\item $x\leqL y \iff x=ay$ for some $a\in \Sone$,
\item $x\leqJ y \iff x=ayb$ for some $a,b\in \Sone$,
\item $x\leqH y \iff x\leqR y$ and $x\leqL y$.
\eit
Note that $x\leqR y$ implies that $\bd(x)=\bd(y)$; other such implications hold, but we will not state them all.  Note that $\Sone$ can be replaced by $S$ in all of the above if $S$ is regular and/or a category.

It is worth noting that $x\leqJ a$ for any $a\in\Pre(x)$, that $x\geqJ a$ for any $a\in\Post(x)$, and that $x\J a$ for any $a\in V(x)$.

If $\K$ is any of $\R$, $\L$, $\J$ or $\H$, then Green's $\K$ relation is defined by ${\K}={\leqK}\cap{\geqK}$.  Green's~$\D$ relation is defined to be ${\D}={\R}\vee{\L}$, the join of $\R$ and $\L$ in the lattice of equivalences on $S$: i.e., the least equivalence containing ${\R}\cup{\L}$.  As with semigroups, we have ${\D}={\R}\circ{\L}={\L}\circ{\R}\sub{\J}$; cf.~\cite[Lemma 2.6]{DE2018}.

The next result shows that the $\R$, $\L$ and $\H$ relations and preorders in a regular partial $*$-semigroup may be characterised equationally (instead of asserting the existence of further elements).
Elements of a regular partial $*$-semigroup of the form $xx^*$ are called \emph{projections}; they may also be characterised as the elements $z$ for which $z^2=z=z^*$; projections are always endomorphisms.

\newpage

\begin{lemma}\label{l:green_S}
If $S$ is a regular partial $*$-semigroup, and if $x,y\in S$, then
\ben\bmc2
\item \label{GS1} $x\leqR y\iff xx^*=yy^*xx^*$,
\item \label{GS2} $x\R y\iff xx^*=yy^*$,
\item \label{GS3} $x\leqL y\iff x^*x=x^*xy^*y$,
\item \label{GS4} $x\L y\iff x^*x=y^*y$.
\emc\een
\end{lemma}

\pf
\firstpfitem{\ref{GS1}}  If $x\leqR y$, then $x=ya$ for some $a\in S$, in which case
\[
xx^* = yaa^*y^* = yy^*yaa^*y^* = yy^*xx^*.
\]
Conversely, if $xx^*=yy^*xx^*$, then $x=xx^*x = yy^*xx^*x \leqR y$.

\pfitem{\ref{GS2}} If $x\R y$, then $x\leqR y$ and $y\leqR x$, and so $xx^*=yy^*xx^*$ and $yy^*=xx^*yy^*$, by \ref{GS1}, which gives
\[
xx^* = (xx^*)^* = (yy^*xx^*)^* = xx^*yy^* = yy^*.
\]
The converse is quickly checked.

\pfitem{\ref{GS3} and \ref{GS4}}  These are dual to \ref{GS1} and \ref{GS2}.
\epf

If $\K$ is any of Green's relations $\R$, $\L$, $\H$, $\J$ or $\D$, then for any $x\in S_{ij}$ we write
\[
K_x = \set{y\in S_{ij}}{x\K y},
\]
and we call these the \emph{$\K$-classes} of $S_{ij}$.
For ${\K}\not={\D}$, Green's preorder $\leqK$ on $S$ induces a partial order on $\K$-classes; we generally denote these partial orders by the same symbol, so for $x,y\in S_{ij}$, we write $K_x\leqK K_y \iff x\leqK y$.  The $\leqJ$ order on $\J$-classes is usually denoted simply by $\leq$.

The structure of a semigroup $T$ can be visualised by means of \emph{eggbox diagrams}.  We draw the elements of a $\D$-class so that all $\R$-related elements are in the same row, $\L$-related elements in the same column, and $\H$-related elements in the same cell.  If an $\H$-class contains an idempotent ($x=x^2$), then that $\H$-class is a group, and it is shaded grey; the group is usually labelled by a standard representative of its isomorphism class (or sometimes we simply list the elements of $T$ in the appropriate cells).  In the case that $T$ is finite, we may draw all the ${\J}={\D}$-classes like this, and illustrate the ${\leq}={\leqJ}$ order by including an edge from $J_x$ up to $J_y$ for each cover $J_x< J_y$.  Many examples of eggbox diagrams are given in this paper; see Figures~\ref{f:emax4}--\ref{f:emax8},~\ref{f:TL4} and~\ref{f:B}.

\subsection{Rectangular bands and groups}\label{ss:RBG}

Recall that a \emph{left-zero semigroup} is a semigroup $U$ with multiplication $u_1u_2=u_1$.  A \emph{left-group} is (isomorphic to) a direct product of a left-zero semigroup and a group; the \emph{degree} of the left-group is the size of the associated left-zero semigroup.    \emph{Right-zero semigroups} and \emph{right-groups} are defined analogously

A \emph{$\rho\times\lam$ rectangular band} is (isomorphic to) a semigroup of the form $U\times V$, where $U$ is a left-zero semigroup of size $\rho$, and $V$ a right-zero semigroup of size $\lam$.  A \emph{$\rho\times\lam$ rectangular group over a group~$G$} is (isomorphic to) a direct product of a $\rho\times\lam$ rectangular band with $G$.  Note that a left-group of degree~$d$ is a $d\times1$ rectangular group.
The $\R$-classes of the rectangular band $U\times V$ are the sets $\{u\}\times V$ ($u\in U$), and a similar statement holds for $\L$-classes.  Thus, a $\rho\times\lam$ rectangular band has $\rho$ $\R$-classes and $\lam$ $\L$-classes.

An \emph{idempotent} of a semigroup $S$ is an element $e$ of $S$ such that $e=e^2$.  For any subset $A$ of $S$, we write $E(A)$ for the set of all idempotents contained in $A$.  An idempotent $e$ of $S$ is \emph{primitive} if $ef=fe=f\implies e=f$ for all idempotents $f$ of $S$.  A semigroup is \emph{simple} if it has a single $\J$-class, and \emph{completely simple} if it is simple and contains a primitive idempotent.  

It is known that a semigroup $S$ is a rectangular band if and only if $x=xyx$ for all $x,y\in S$; see \cite[page 7]{Howie}.
It is known that a semigroup is a rectangular group if and only if it is completely simple, regular and its idempotents form a subsemigroup; see \cite[Exercise~10, page~139]{Howie}.
The next result is stated without proof in \cite[Proposition 1.6]{Ault1973}; we provide one for convenience.

\begin{lemma}\label{l:RG}
A semigroup $S$ is a rectangular group if and only if it is regular and $E(S)$ is a rectangular band.
\end{lemma}

\pf
The forwards implication being clear, suppose $S$ is regular and that $E(S)$ is a rectangular band.  As noted above, it suffices to show that $S$ is completely simple.  If $e,f\in E(S)$ are such that $ef=fe=f$, then since $e=efe$ (also noted above), we have $e=efe=ef=f$, so that in fact every idempotent is primitive.  So it remains to show that $S$ is simple; we do this by showing that~$S$ is a single $\D$-class (meaning that $S$ is in fact \emph{bisimple}).

To do so, let $x,y\in S$.  Since $S$ is regular, we have $x\R e$ and $y\L f$ for some idempotents $e,f\in E(S)$.  Since $E(S)$ is a rectangular band, $e$ and $f$ are $\D$-related in $E(S)$, and hence also in $S$.  But then $x\R e\D f\L y$, which gives $x\D y$.
\epf

Recall from \cite{DE2018} that a (partial) semigroup~$S$ is \emph{stable} if for all $x,u\in S$,
\begin{equation}\label{e:stab}
x\J xu\implies x\R xu \AND x\J ux\implies x\L ux.
\end{equation}
Stability is an extremely useful property; for example, we have ${\J}={\D}$ in any stable (partial) semigroup; see \cite[Lemma 2.6]{DE2018}.
It is possible to prove the next two results using the Rees Theorem (cf.~\cite[Theorem~3.2.3]{Howie}), but we prefer the more direct approach here.

\begin{lemma}\label{l:ED}
Let $D$ be a regular $\D$-class of a stable semigroup $S$.  If $E(D)$ is a subsemigroup of $S$, then $E(D)$ is a rectangular band, and $D$ is a rectangular group.
\end{lemma}

\pf
We show that $E(D)$ is a rectangular band by showing that $xyx=x$ for all $x,y\in E(D)$.  So fix some $x,y\in E(D)$.  Since $xy\in E(D)$ we have $xy\D x$, so $xy\R x$ by stability.  Thus, $x=xyu$ for some $u\in\Sone$.  Since $xy\in E(D)$, it follows that $x=(xy)^2u=xy(xyu)=xyx$, as required.

By Lemma \ref{l:RG}, and since $D$ is regular, it remains to show that $D$ is a semigroup.  To do so, let $x,y\in D$.  Fix inverses $a\in V(x)$ and ${b\in V(y)}$.  Then $x y=xax yby \leqJ ax yb \leqJ x y$, so that $x y\J ax yb$, whence ${xy\in J_{axyb}= D_{axyb}}$ (the latter since ${\J}={\D}$, as $S$ is stable).  But $ax,yb\in E(D)$, so $ax yb\in E(D)\sub D$.  It follows that $D_{axyb}=D$, and so $xy\in D$.
\epf

\begin{rem}
Stability is a crucial assumption in Lemma \ref{l:ED}.  For example, the bicyclic semigroup is a single regular $\D$-class, and its idempotents form a subsemigroup; but it is certainly not a rectangular group.
\end{rem}

\begin{lemma}\label{l:EL}
If a regular $\D$-class of a semigroup is an $\L$-class, then it is a left-group.
\end{lemma}

\pf
Let $D$ be the $\D$-class in question.  Since a semigroup is a left-group if and only if it is regular and $\L$-universal~\cite[Section 1.11]{CPbook1}, it is enough to show that $D$ is a semigroup.  So let $x,y\in D$.  Since~$D$ is regular, we have $y\R e$ for some idempotent $e\in E(D)$.  Since $D$ is an $\L$-class, we have $x\L e$.  Thus, $L_x\cap R_y$ contains an idempotent, so \cite[Proposition 2.3.7]{Howie} tells us that $xy\in R_x\cap L_y\sub D$.
\epf

\section{Sandwich semigroups}\label{s:sandwich}

For the duration of Section \ref{s:sandwich}, we fix a partial semigroup $S\equiv(S,I,\bd,\br,\cdot)$.  We also fix some element $a\in S_{ji}$, where $i,j\in I$.  As in \cite{DE2018}, an associative operation $\star_a$ may be defined on $S_{ij}$ by $x\star_ay=xay$ for $x,y\in S_{ij}$.  The resulting semigroup $(S_{ij},\star_a)$ is called a \emph{sandwich semigroup}, and is denoted $S_{ij}^a$.  When $|I|=1$, $S$ is a semigroup, and any sandwich semigroup in $S$ is a so-called \emph{variant} $S^a=(S,\star_a)$ of $S$; cf.~\cite{DE2015,Hickey1986,Hickey1983,KL2001}.  In the case that~$S$ is a partial $*$-semigroup, it is easy to check that the map
\[
S_{ij}\to S_{ji}: x\mt x^*
\]
determines an anti-isomorphism $S_{ij}^a\to S_{ji}^{a^*}$.

This section outlines some general machinery for working with arbitrary sandwich semigroups that will be used when studying diagram categories in Sections \ref{s:DC}--\ref{s:B}.  Sections \ref{ss:GreenSija} and \ref{ss:nonsandwich} mostly revise results we need from \cite{DE2018,Sandwich1}, while Sections \ref{ss:<Sija}--\ref{ss:RI} develop a substantial new theory of Green's preorders in sandwich semigroups, with a focus on maximal elements and classes in these orderings.

\subsection{Green's relations and regularity}\label{ss:GreenSija}

Let $\K$ be any of $\R$, $\L$, $\J$, $\H$ or $\D$.  We denote Green's $\K$ relation on the sandwich semigroup~$S_{ij}^a$ by $\gKa$.  Note for example that if $x,y\in S_{ij}$, then
\[
x\R^a y \iff x=y \text{ or } [x=yau \text{ and } y=xav \text{ for some } u,v\in S_{ij}.]
\]
For $x\in S_{ij}$, we write
\[
K_x = \set{y\in S_{ij}}{x\K y} \AND K_x^a = \set{y\in S_{ij}}{x\gKa y}
\]
for the $\K$- and $\gKa$-classes of $x$ in $S_{ij}$, respectively; note that $K_x^a\sub K_x$.

To understand the $\gKa$ relations, a crucial role is played by the sets $P_1^a$, $P_2^a$, $P_3^a$ and ${P^a=P_1^a\cap P_2^a}$, where
\[
P_1^a = \set{x\in S_{ij}}{xa\R x} \COMMA
P_2^a = \set{x\in S_{ij}}{ax\L x} \COMMA
P_3^a = \set{x\in S_{ij}}{axa\J x} .
\]
Note that $x\in P_1^a$ if and only if $x=xav$ for some $v\in\Sone$, in which case $x=(xav)av=xa(vav)$ with $vav\in S$ (not just $\Sone$); note then that in fact $vav\in S_{ij}$.  With similar reasoning for the other sets, it follows that
\[
P_1^a = \set{x\in S_{ij}}{x\in xaS_{ij}} \COMMA
P_2^a = \set{x\in S_{ij}}{x\in S_{ij}ax} \COMMA
P_3^a = \set{x\in S_{ij}}{x\in S_{ij}axaS_{ij}} .
\]
We will use this fact frequently in what follows, usually without explicit reference.

The following is \cite[Theorem 2.13]{DE2018}.

\begin{thm}\label{t:Green_Sij}
Let $S\equiv(S,I,\bd,\br,\cdot)$ be a partial semigroup, and let $a\in S_{ji}$ where $i,j\in I$.  If $x\in S_{ij}$, then   
\ben
\begin{multicols}{2}
\item \label{GSij1} $R_x^a = \begin{cases}
R_x\cap P_1^a &\text{if $x\in P_1^a$}\\
\{x\} &\text{if $x\in S_{ij}\sm P_1^a$,}
\end{cases}$
\item \label{GSij2} $L_x^a = \begin{cases}
L_x\cap P_2^a &\hspace{0.7mm}\text{if $x\in P_2^a$}\\
\{x\} &\hspace{0.7mm}\text{if $x\in S_{ij}\sm P_2^a$,}
\end{cases}
$
\item \label{GSij3} $H_x^a = \begin{cases}
H_x &\hspace{7.4mm}\text{if $x\in P^a$}\\
\{x\} &\hspace{7.4mm}\text{if $x\in S_{ij}\sm P^a$,}
\end{cases}$
\item \label{GSij4} $D_x^a = \begin{cases}
D_x\cap P^a &\text{if $x\in P^a$}\\
L_x^a &\text{if $x\in P_2^a\sm P_1^a$}\\
R_x^a &\text{if $x\in P_1^a\sm P_2^a$}\\
\{x\} &\text{if $x\in S_{ij}\sm (P_1^a\cup P_2^a)$,}
\end{cases}$
\item \label{GSij5} $J_x^a = \begin{cases}
J_x\cap P_3^a &\hspace{2.2mm}\text{if $x\in P_3^a$}\\
D_x^a &\hspace{2.2mm}\text{if $x\in S_{ij}\sm P_3^a$.}
\end{cases}$
\end{multicols}
\een
Further, if $x\in S_{ij}\sm P^a$, then $H_x^a=\{x\}$ is a non-group $\gHa$-class of $S_{ij}^a$.  \epfres
\end{thm}

The set $P^a=P_1^a\cap P_2^a$ may also be used to describe the regular elements of the sandwich semigroup~$S_{ij}^a$.  For any semigroup $T$, we denote by $\Reg(T)$ the set of regular elements; this need not be a subsemigroup of $T$.  The next result is a special case of \cite[Proposition 2.7]{Sandwich1}:

\begin{prop}\label{p:P}
Let $a$ be an element of a regular partial semigroup $S$, say with $a\in S_{ji}$.  Then $\Reg(S_{ij}^a) = P^a$ is a (regular) subsemigroup of $S_{ij}^a$.  \epfres
\end{prop}

In the case of regular partial $*$-semigroups, the sets $P_1^a$ and $P_2^a$ may be described equationally as follows.  

\begin{lemma}\label{l:P}
Let $S$ be a regular partial $*$-semigroup, let $i,j\in I$, and let $a\in S_{ji}$.  Then
\ben
\item \label{P1} $P_1^a 
= \set{x\in S_{ij}}{x^*x\in \Post(aa^*)}
= \set{x\in S_{ij}}{aa^*\in \Pre(x^*x)}$,
\item \label{P2} $P_2^a 
= \set{x\in S_{ij}}{xx^*\in \Post(a^*a)}
= \set{x\in S_{ij}}{a^*a\in \Pre(xx^*)}$.
\een
\end{lemma}

\pf
We just prove \ref{P1}, as \ref{P2} is dual.  If $x\in P_1^a$, then $x\R xa$ so by Lemma \ref{l:green_S}\ref{GS2} we have ${xx^*=xa(xa)^*=xaa^*x^*}$; it follows that $x^*x = x^*(xx^*)x = x^*x aa^* x^*x$, and so $aa^*\in\Pre(x^*x)$.  The reverse containment is similar.
\epf

Although it is not necessary for our purposes, the next result seems to be of independent interest.  Recall that an element $a$ of a regular partial $*$-semigroup is a projection if $a^2=a=a^*$; such an element must belong to an endomorphism semigroup $S_i$.  

\begin{cor}\label{c:rss}
If $S$ is a regular partial $*$-semigroup, and if $a\in S_i$ is a projection, then $\Reg(S_i^a)$ is a regular $*$-semigroup (with involution inherited from $S$).
\end{cor}

\pf
Since $S$ is regular, Proposition \ref{p:P} tells us that $\Reg(S_i^a)=P^a$.  It remains to check that 
\[
x^{**}=x \COMMA (x\star_ay)^*=y^*\star_ax^* \COMMA x=x\star_ax^*\star_ax \qquad\text{for all $x,y\in P^a$.}
\]
The first is clear, and for the second we have $(x\star_ay)^*=(xay)^*=y^*a^*x^*=y^*ax^*=y^*\star_ax^*$.  For the third, let $x\in P^a$.  Since $a$ is a projection, we have $a=aa=aa^*$.  Since $x\in P_1^a$, it follows from Lemma~\ref{l:P}\ref{P1} that $x^*x=x^*xax^*x$, and so $xx^*=x(x^*x)x^*=x(x^*xax^*x)x^*=xax^*$.  A similar calculation gives $x^*x=x^*ax$.  Together these give $x\star_ax^*\star_ax=xax^*ax=xx^*ax=xx^*x=x$.
\epf

\begin{rem}
A special case of \cite[Proposition 5.1]{Sandwich1} says that if $S$ is an inverse category (a more restrictive class than regular partial $*$-semigroups \cite{Kastl1979}), then $\Reg(S_{ij}^a)$ is an inverse semigroup for \emph{any} $a\in S_{ji}$.  We cannot similarly strengthen Corollary \ref{c:rss} to work for arbitrary $a$ in a regular partial $*$-semigroup.  Indeed, consider the regular semigroup $\Reg(\B_{64}^{\si_2})$ in Remark \ref{r:B} below.  As can be seen from Figure \ref{f:B}, any $\D$-class in this semigroup contains unequal numbers of $\R$- and $\L$-classes; it follows from this that $\Reg(\B_{64}^{\si_2})$ does not even have an involution.  (See also Remark \ref{r:TL4} and Figure \ref{f:TL4}.)
\end{rem}

\subsection{Green's preorders}\label{ss:<Sija}

We now prove two results concerning Green's preorders in $S_{ij}^a$.  Both concern the case in which the sandwich element $a$ has a left- and right-identity in $S$: i.e., $a=ea=af$ for some $e,f\in S$.  This condition holds, for example, if $a$ is regular, or if $S$ is a category.
The first result compares the preorders $\leq_{\K^a}$ on $S_{ij}^a$ to the preorders $\leqK$ on $S$.  

\begin{lemma}\label{l:<}
If $a\in S_{ji}$ has a left- and right-identity in $S$, and if $x,y\in S_{ij}$, then 
\ben
\item \label{<11} $x\leq_{\R^a} y\iff x=y $ or $ x \leqR ya$,
\item \label{<12} $x\leq_{\L^a} y\iff x=y $ or $ x\leqL ay$.
\item \label{<13} $x\leq_{\J^a} y\iff x=y $ or $ x \leqR ya $ or $ x\leqL ay $ or $ x\leqJ aya$.  
\een
\end{lemma}

\pf
We just prove \ref{<13} as the others are similar.  First note that $x\leqJa y$ if and only if one of the following holds:
\bena
\bmc2
\item \label{<a} $x=y$,
\item \label{<b} $x=uay$ for some $u\in S_{ij}$,
\item \label{<c} $x=yav$ for some $v\in S_{ij}$,
\item \label{<d} $x=uayav$ for some $u,v\in S_{ij}$.
\emc
\een
Clearly \ref{<b}$\implies$$x\leqL ay$, while \ref{<c}$\implies$$x\leqR ya$ and \ref{<d}$\implies$$x\leqJ aya$.  The converses of these are easily established.  For example, if $x\leqJ aya$, then $x=sayat$ for some $s,t\in\Sone$; if $e,f\in S$ are such that $a=ea=af$, then $x=s(ea)y(af)t$, so that \ref{<d} holds (with $u=se$ and $v=ft$, both from $S_{ij}$).
\epf

The next result is a generalisation of \cite[Proposition 3.21]{Sandwich2}, which concerns the special case of the category of partial maps.  It shows how the $\leqJa$ preorder on $S_{ij}^a$ simplifies when elements of~$P_1^a$,~$P_2^a$ or~$P_3^a$ are involved, under the same hypothesis on $a$ having left- and right-identities.  (Note that $P^a=P_1^a\cap P_2^a\sub P_3^a$.)

\newpage

\begin{prop}\label{p:<}
Suppose $a\in S_{ji}$ has a left- and right-identity in $S$, and let $x,y\in S_{ij}$.
\ben
\item \label{<1} If $x\in P_1^a$, then $x\leqJa y \iff x\leqJ aya$ or $x\leqR ya$.
\item \label{<2} If $x\in P_2^a$, then $x\leqJa y \iff x\leqJ aya$ or $x\leqL ay$.
\item \label{<3} If $x\in P_3^a$, then $x\leqJa y \iff x\leqJ aya$.
\item \label{<4} If $y\in P_1^a$, then $x\leqJa y \iff x\leqJ ay$ or $x\leqR y$.
\item \label{<5} If $y\in P_2^a$, then $x\leqJa y \iff x\leqJ ya$ or $x\leqL y$.
\item \label{<6} If $y\in P_3^a$, then $x\leqJa y \iff x\leqJ y$.
\een
\end{prop}

\pf
Again note that $x\leqJa y$ if and only if one of \ref{<a}--\ref{<d} holds, as in the proof of Lemma \ref{l:<}.
We just prove \ref{<1} and \ref{<4}, as the others are very similar.  

\pfitem{\ref{<1}}  Suppose $x\in P_1^a$, so that $x=xaz$ for some $z\in S_{ij}$.  

Suppose first that $x\leqJa y$, so that one of \ref{<a}--\ref{<d} holds.  If \ref{<c} holds, then $x\leqR ya$.  If \ref{<d} holds, then $x\leqJ aya$.  If \ref{<a} holds, then $x=xaz=yaz\leqR ya$.  Similarly, \ref{<b} implies $x\leqJ aya$.

The converse follows from Lemma \ref{l:<}\ref{<13}.

\pfitem{\ref{<4}}  Suppose $y\in P_1^a$, so that $y=yaz$ for some $z\in S_{ij}$.  

If $x\leqJa y$, then one of \ref{<a}--\ref{<d} holds.  Clearly \ref{<a} and \ref{<c} each imply $x\leqR y$, while \ref{<b} and \ref{<d} each imply $x\leqJ ay$.

Conversely, if $x\leqJ ay$, then $x=sayt$ for some $s,t\in\Sone$; if $e\in S$ is such that $a=ea$, then $x=s(ea)(yaz)t=(se)aya(zt)$, with $se,zt\in S_{ij}$, so that \ref{<d} holds.  On the other hand, if $x\leqR y$, then $x=yt$ for some $t\in\Sone$, in which case $x=(yaz)t=ya(zt)$, and \ref{<d} holds.
\epf

\begin{rem}
In the proofs of Lemma \ref{l:<} and Proposition \ref{p:<}, the forwards implications did not require the assumption on $a$ having identities.

One could also prove similar simplifying statements for the ${\R^a}/{\R}$ and ${\L^a}/{\L}$ relations: e.g.,
\bit
\item If $x\in P_1^a$, then $x\leq_{\R^a} y\iff x\leqR ya$.
\item If $y\in P_1^a$, then $x\leq_{\R^a} y \iff x\leqR y$.
\eit
\end{rem}

\subsection{Maximal $\J^a$-classes}\label{ss:max}

As in Section \ref{ss:Green}, if $\K$ is any of $\R$, $\L$, $\H$ or $\J$, then the $\leq_{\K^a}$ preorder on $S_{ij}^a$ induces a partial order also denoted $\leq_{\K^a}$ on the $\K^a$-classes of $S_{ij}^a$: for $x,y\in S_{ij}$, we have $K_x^a\leq_{\K^a}K_y^a \iff x\leq_{\K^a}y$.  
We will typically denote the $\leqJa$ order on $\J^a$-classes simply by $\leq$.  By a \emph{maximal $\J^a$-class} in~$S_{ij}^a$, we mean a $\J^a$-class that is maximal with respect to this order.  These maximal $\J^a$-classes will play an important role in our investigation of diagram categories, and here we prove some general results concerning them.  

Our first result identifies a natural family of maximal $\J^a$-classes.

\begin{lemma}\label{l:max1}
If $x\in S_{ij}$ is such that $x\not\leqJ a$ in $S$, then $\{x\}$ is a maximal $\J^a$-class in~$S_{ij}^a$; additionally,~$\{x\}$ is a nonregular $\D^a$-class.
\end{lemma}

\pf
Fix some such element $x$.  To show that $\{x\}$ is a maximal $\J^a$-class, it suffices to show that for any $y\in S_{ij}$, $x\leqJa y$ implies $x=y$.  To do so, suppose $y\in S_{ij}$ and $x\leqJa y$.  Then one of~\ref{<a}--\ref{<d} holds, as in the proof of Lemma \ref{l:<}.  Any of \ref{<b}--\ref{<d} would imply $x\leqJ a$ (in $S$), contrary to assumption, so in fact $x=y$.

Now $\{x\}\sub D_x^a\sub J_x^a=\{x\}$, so that $\{x\}=D_x^a$ is indeed a $\D^a$-class.  If it was a regular $\D^a$-class, then it would contain an idempotent; but then $x=x\star_ax=xax\leqJ a$, a contradiction.
\epf

We will call a maximal $\J^a$-class of $S_{ij}^a$ \emph{trivial} if it is of the form described in Lemma \ref{l:max1}.  Any other $\J^a$-class will be called \emph{nontrivial}.
Nontrivial maximal $\J^a$-classes do not always exist (cf.~Example~\ref{e:max}).
The next result concerns nontrivial maximal $\J^a$-classes in the case that the sandwich element $a\in S_{ji}$ is regular.

\begin{lemma}\label{l:max2}
Suppose $a\in S_{ji}$ is regular.
\ben
\item \label{max21} There is at most one nontrivial maximal $\J^a$-class in $S_{ij}^a$.  
\item \label{max22} If a nontrivial maximal $\J^a$-class exists, then it contains $\Pre(a)$.  
\item \label{max23} If a nontrivial maximal $\J^a$-class exists, and if it is a $\D^a$-class, then it is regular.
\een
\end{lemma}

\pf
Since \ref{max21} clearly follows from \ref{max22}, it suffices to prove \ref{max22} and \ref{max23}.

\pfitem{\ref{max22}}  
Suppose $J$ is a nontrivial maximal $\J^a$-class.  So $J=J_x^a$ for some $x\in S_{ij}$ with $x\leqJ a$.  Let $b\in\Pre(a)$ be arbitrary: i.e., $a=aba$.  Since $x\leqJ a$, we have $x=uav$ for some $u,v\in\Sone$.  But then $x=uabababav=(uab)\star_ab\star_a(bav)$, with $uab,bav\in S_{ij}$, so that $J=J_x^a\leq J_b^a$.  Maximality then gives $J=J_b^a$: i.e., $b\in J$.  

\pfitem{\ref{max23}}  Suppose $J$ is a nontrivial maximal $\J^a$-class; it suffices to show that $J$ contains an idempotent.  Let $b\in V(a)$.  Then we have $b\in\Pre(a)$, so $J=J_b^a$ by \ref{max22}: i.e., $b\in J$.  But $b=bab=b\star_ab$.
\epf

We now give necessary and sufficient conditions for a nontrivial maximal $\J^a$-class to exist, again in the case that the sandwich element $a\in S_{ji}$ is regular.

\begin{prop}\label{p:Jba}
Suppose $a\in S_{ji}$ is regular.  Then the following are equivalent:
\ben
\item \label{Jba1} $S_{ij}^a$ has a nontrivial maximal $\J^a$-class,
\item \label{Jba2} for all $x\in S_{ij}$, $a \leqJ axa \implies x\leqJ a$,
\item \label{Jba3} for all $x\in S_{ij}$, $a \J axa \implies x\J a$.
\een
\end{prop}

\pf
Since $a$ is regular, we may fix some $b\in V(a)$.  Write $J=J_b^a$, and note that if $S_{ij}^a$ has a nontrivial maximal $\J^a$-class then it must be $J$ (cf.~Lemma \ref{l:max2}\ref{max22}).  Since $b\in V(a)$, we have $a\J b$.  Also, from $b=b(aba)b$ we deduce that $b\in P_3^a$ (indeed, $b\in P_1^a\cap P_2^a=P^a$).

\pfitem{\ref{Jba1}$\implies$\ref{Jba2}}  Aiming to prove the contrapositive, suppose \ref{Jba2} does not hold.  So there exists $x\in S_{ij}$ such that $a\leqJ axa$ and $x\not\leqJ a$.  
From $x\not\leqJ a$, Lemma \ref{l:max1} tells us that $J_x^a=\{x\}$ is a maximal $\J^a$-class, and that $x$ is nonregular.  
From $a\leqJ axa$ and $b\J a$, we have $b\leqJ axa$, so $b\leqJa x$ by Lemma~\ref{l:<}\ref{<13}.  This means that $J=J_b^a\leq J_x^a$.  But $J\not=J_x^a$, since~$J$ contains the regular element $b$, and since $x$ is nonregular.  It follows that $J<J_x^a$, and so $J$ is not maximal.

\pfitem{\ref{Jba2}$\implies$\ref{Jba3}}  Suppose condition \ref{Jba2} holds, and suppose $x\in S_{ij}$ satisfies $a\J axa$.  Then in particular, $a\leqJ axa$, so \ref{Jba2} gives $x\leqJ a$.  But also $a\leqJ axa\leqJ x$, so it follows that $x\J a$.

\pfitem{\ref{Jba3}$\implies$\ref{Jba1}}  Suppose condition \ref{Jba3} holds.  We will show that $J$ is maximal (and it is nontrivial since it contains the regular element $b$).  To do so, suppose $x\in S_{ij}$ is such that $J\leq J_x^a$.  We must show that $J=J_x^a$, and it suffices to show that $J_x^a\leq J$.  Since $J_b^a=J\leq J_x^a$, we have $b\leqJa x$.  Since $b\in P_3^a$, Proposition \ref{p:<}\ref{<3} gives $b\leqJ axa$.  Together with $a\J b$, it follows that $a\leqJ axa$, so that $a\J axa$.  It follows from \ref{Jba3} that $x\J a$.  In particular, $x\leqJ a$, and so $x\leqJ b$ (as $a\J b$).  Since $b\in P_3^a$, Proposition \ref{p:<}\ref{<6} then gives $x\leqJa b$, so that $J_x^a\leq J_b^a=J$, as required.
\epf

A simple consequence worth noting is as follows:

\begin{cor}\label{c:max}
If $a\in S_{ji}$ has a pre-inverse that is not $\J$-related to $a$ (in $S$), then $S_{ij}^a$ has only trivial maximal $\J^a$-classes.
\end{cor}

\pf
Suppose $x$ is such a pre-inverse.  Since $a=axa$, certainly $a\J axa$; since $(a,x)\not\in{\J}$, by assumption, the implication \ref{Jba3} in Proposition \ref{p:Jba} does not hold.  It follows by that proposition that~$S_{ij}^a$ has no nontrivial maximal $\J^a$-class.  
\epf

Although the converse of Corollary \ref{c:max} does not hold in general (see Example \ref{e:max}\ref{emax8} below), it does in a certain special case.  To state this (see Proposition \ref{p:max2}), we require the concept of stability, as defined in \eqref{e:stab}.  We have already noted that ${\J}={\D}$ in any stable partial semigroup \cite[Lemma~2.6]{DE2018}.  A special case of \cite[Lemma 2.6]{Sandwich1} says that if the partial semigroup $S$ is stable, then every sandwich semigroup $S_{ij}^a$ is stable.

\begin{lemma}\label{l:axa}
If $S$ is stable, and if $a\in S_{ji}$ and $x\in S_{ij}$, then
\ben
\item \label{axa1} $a\J axa \iff a\H axa$,
\item \label{axa2} if $a\J x$, then $x=xax\iff a=axa$.
\een
\end{lemma}

\pf
\firstpfitem{\ref{axa1}}
Since ${\H}\sub{\J}$, it suffices to prove the forwards implication, so suppose $a\J axa$.  From $a\J a(xa)$, stability gives $a\R a(xa)$; similarly, $a\L(ax)a$, and so $a\H axa$.

\pfitem{\ref{axa2}}
By symmetry it suffices to prove only the forwards implication, so suppose $x=xax$.  Now, ${ax\leqJ a \J x=xax \leqJ ax}$; it follows that $ax\J a$, so stability gives $ax\R a$.  Thus, ${a=axv}$ for some $v\in\Sone$.  But then ${axa = ax(axv) = a(xax)v = axv = a}$.
\epf

\begin{rem}\label{r:axa}
Stability is necessary in both parts of Lemma \ref{l:axa}.  For example, suppose $S$ is a monoid with identity~$1$, and that there is a nonidentity idempotent $e$ with $e\J1$ (the bicyclic monoid has this property).  Then $e=e1e$ yet $(1,1e1)=(1,e)\in{\J}\sm{\H}$.

In the case that $S$ is a semigroup, Lemma \ref{l:axa}\ref{axa1} is \cite[Exercise A.2.2.1]{RSbook}.
\end{rem}

\begin{prop}\label{p:max2}
If $S$ is stable and $\H$-trivial (i.e., $\H$ is the trivial relation), and if $a\in S_{ji}$ is regular, then the following are equivalent:
\ben
\item \label{pmax21} $S_{ij}^a$ has a nontrivial maximal $\J^a$-class,
\item \label{pmax22} every pre-inverse of $a$ is $\J$-related to $a$ (in $S$),
\item \label{pmax23} $\Pre(a)=V(a)$.
\een
\end{prop}

\pf
\firstpfitem{\ref{pmax21}$\iff$\ref{pmax22}}  
By Lemma \ref{l:axa}\ref{axa1} and $\H$-triviality, we have $a\J axa\iff a=axa$, for $x\in S_{ij}$.  Thus, condition \ref{Jba3} of Proposition \ref{p:Jba} is equivalent to:
\bit
\item for all $x\in S_{ij}$, $a=axa \implies x\J a$: i.e., every pre-inverse of $a$ is $\J$-related to $a$.
\eit
\pfitem{\ref{pmax22}$\implies$\ref{pmax23}} 
Clearly it suffices to show that $\Pre(a)\sub V(a)$.  With this in mind, fix some $x\in\Pre(a)$, so that $a=axa$.  By \ref{pmax22} we have $x\J a$; it follows from Lemma \ref{l:axa}\ref{axa2} that $x=xax$, and so $x\in V(a)$.

\pfitem{\ref{pmax23}$\implies$\ref{pmax22}} 
This is clear since every element of $V(a)$ is $\J$-related to $a$.
\epf

We now consider a number of examples, illustrating the above results on maximal $\J^a$-classes.  In each case, the partial semigroup $S$ is in fact a semigroup, so the sandwich semigroups are all \emph{variants}~$S^a$.  In particular, examples \ref{emax5} and \ref{emax6} show that Lemma \ref{l:max2}\ref{max21} is not true in general when~$a$ is nonregular; example \ref{emax8} shows that Proposition \ref{p:max2} need not hold if $S$ is not $\H$-trivial.

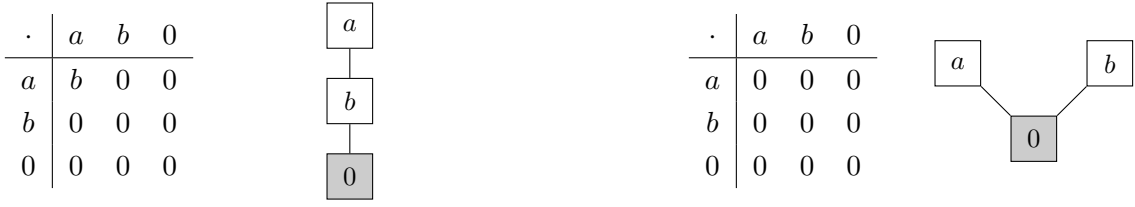
\begin{figure}[ht]
\begin{center}
\begin{tikzpicture}
\node () at (0,0) {
\begin{tabular}{c|ccc}
$\cdot$ & $a$ & $b$ & $0$ \\
\hline
$a$ & $b$ & $0$ & $0$ \\
$b$ & $0$ & $0$ & $0$ \\
$0$ & $0$ & $0$ & $0$ \\
\end{tabular}
};
\node () at (3.3,0) {
\begin{tikzpicture}[scale=1,
block/.style={
draw,
fill=white,
rectangle, 
minimum width=.6cm,
minimum height=.6cm,
font=\small}]
\node[block] (a) at (0,2) {$a$};
\node[block] (b) at (0,1) {$b$};
\node[block,fill=black!20] (0) at (0,0) {$0$};
\draw (a)--(b)--(0);
\end{tikzpicture}
};
\node () at (9,0) {
\begin{tabular}{c|ccc}
$\cdot$ & $a$ & $b$ & $0$ \\
\hline
$a$ & $0$ & $0$ & $0$ \\
$b$ & $0$ & $0$ & $0$ \\
$0$ & $0$ & $0$ & $0$ \\
\end{tabular}
};
\node () at (12.3,0) {
\begin{tikzpicture}[scale=1,
block/.style={
draw,
fill=white,
rectangle, 
minimum width=.6cm,
minimum height=.6cm,
font=\small}]
\node[block] (a) at (-1,1) {$a$};
\node[block] (b) at (1,1) {$b$};
\node[block,fill=black!20] (0) at (0,0) {$0$};
\draw (a)--(0)--(b);
\end{tikzpicture}
};
\end{tikzpicture}
\caption{Multiplication table and eggbox diagram for $S$ (left) and $S^a$ (right), as in Example \ref{e:max}\ref{emax5}.}
\label{f:emax5}
\end{center}
\end{figure}

\begin{figure}[ht]
\begin{center}
\includegraphics[height=4cm]{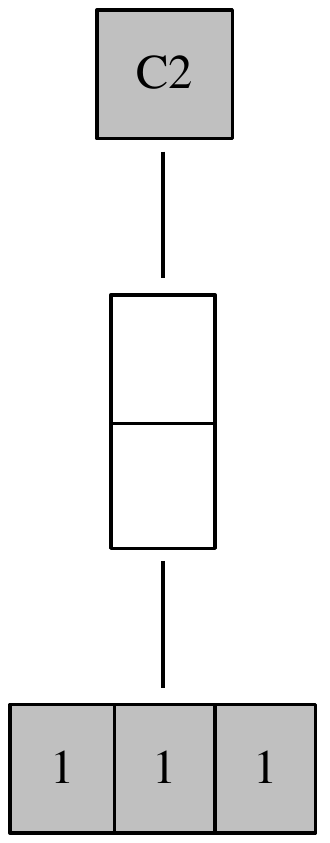} 
\qquad\qquad\qquad
\includegraphics[height=2cm]{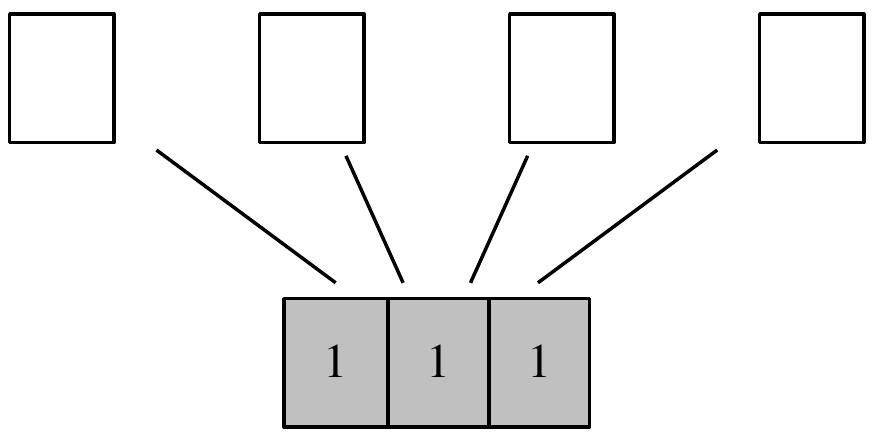} 
\caption[blah]{Eggbox diagrams of the semigroup $S$ (left) and variant $S^a$ (right) from Example \ref{e:max}\ref{emax6}.}
\label{f:emax6}
\end{center}
\end{figure}

\begin{figure}[ht]
\begin{center}
\includegraphics[height=5.1cm]{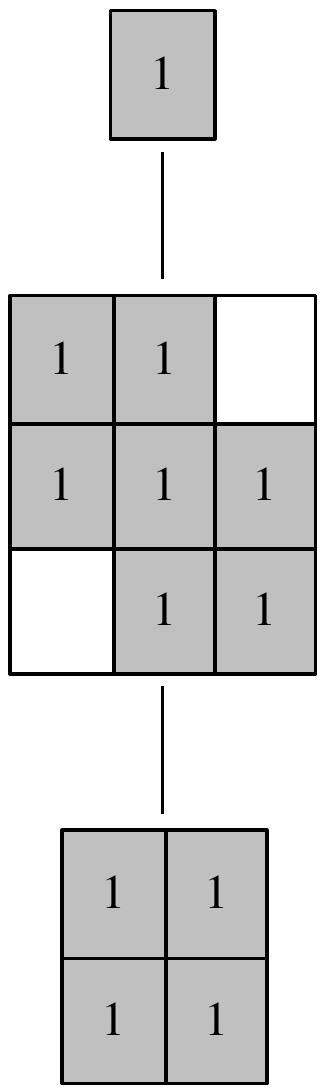} 
\qquad\qquad\qquad
\includegraphics[height=7.5cm]{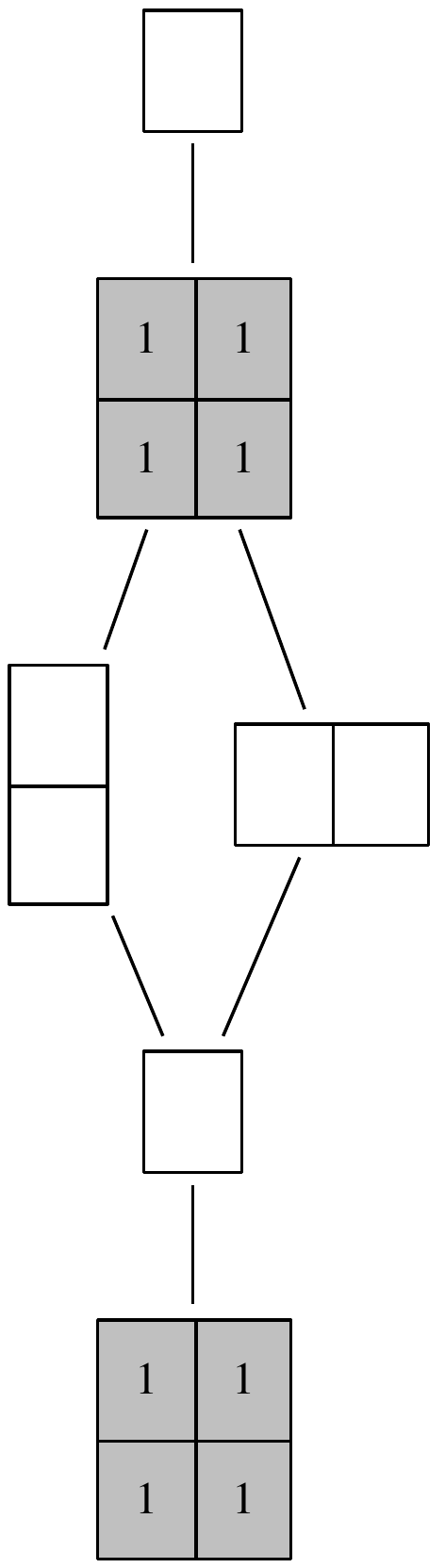} 
\qquad\qquad\qquad
\includegraphics[height=6.3cm]{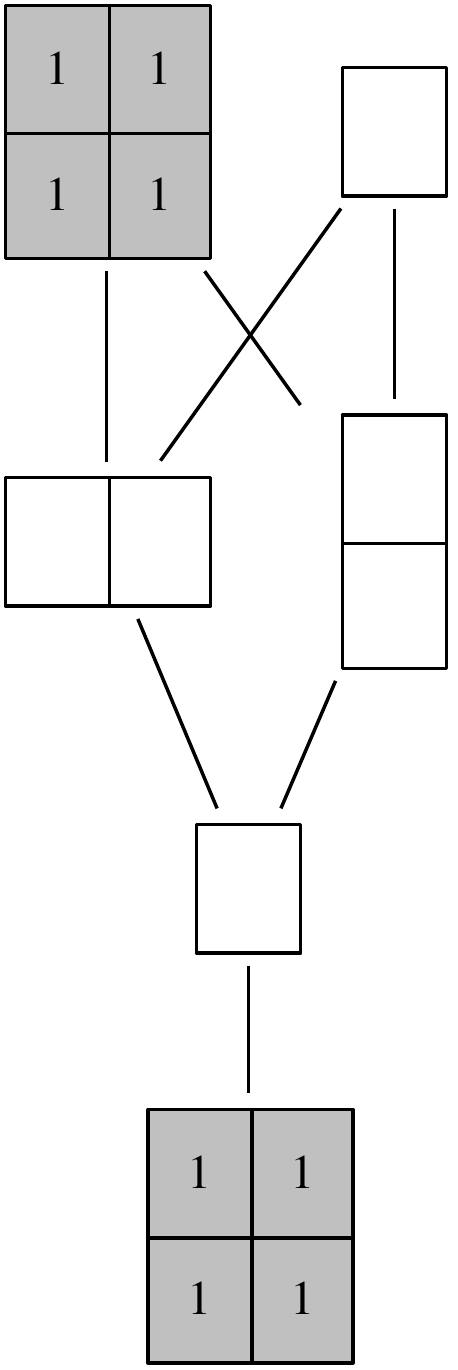} 
\qquad\qquad\qquad
\includegraphics[height=2cm]{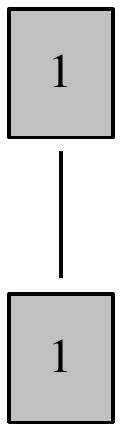} 
\caption[blah]{Eggbox diagrams of the Temperley-Lieb monoids $\TL_4$ (left) and $\TL_2$ (right), and the variants~$\TL_4^\si$ (middle-left) and~$\TL_4^\tau$ (middle-right) from Example \ref{e:max}\ref{emax7}.}
\label{f:emax7}
\end{center}
\end{figure}

\begin{figure}[ht]
\begin{center}
\includegraphics[height=7.3cm]{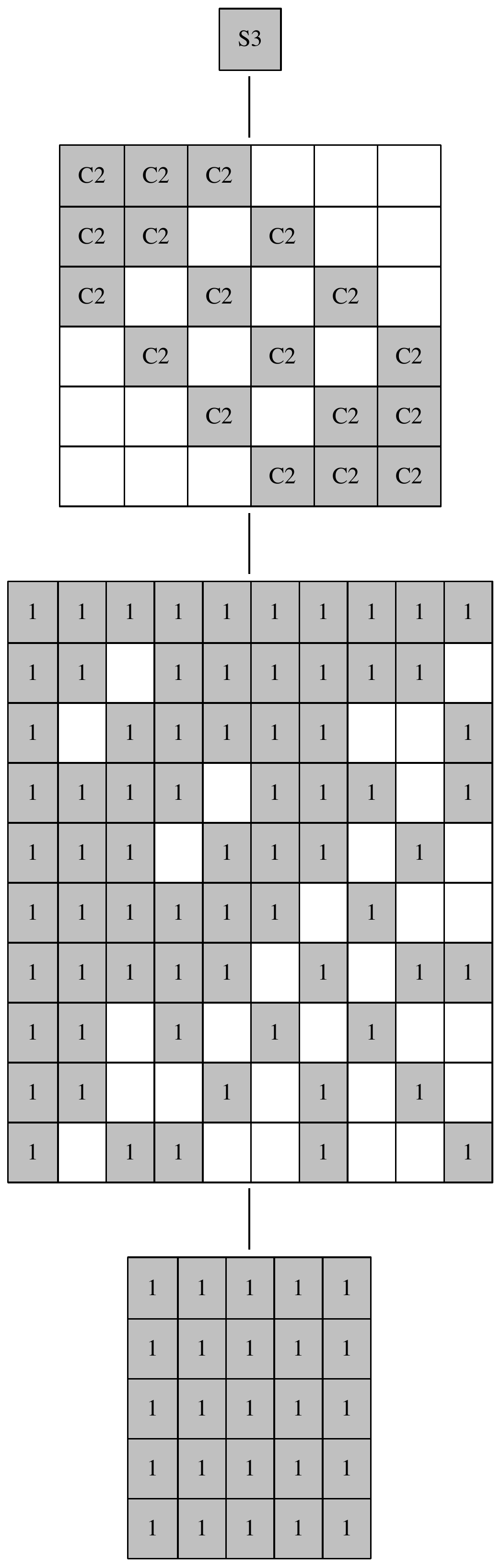} 
\qquad
\includegraphics[height=9cm]{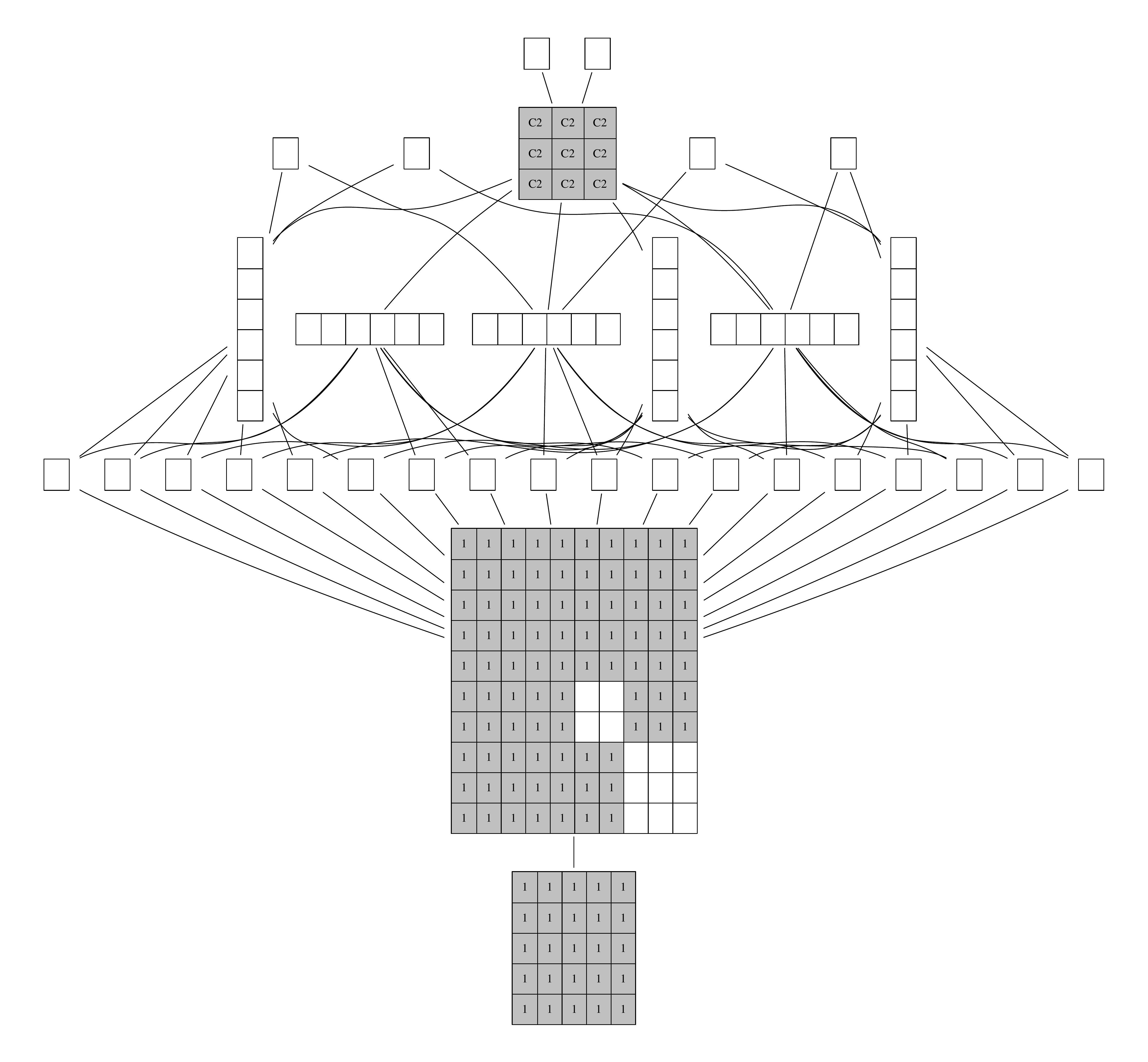} 
\qquad
\includegraphics[height=4cm]{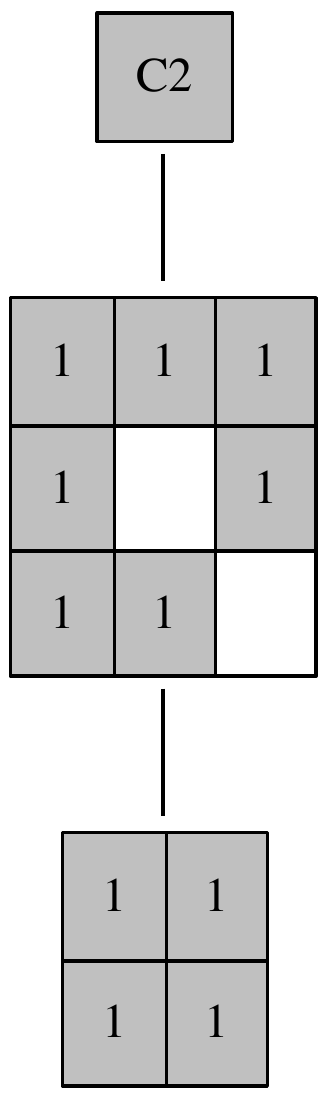} 
\caption[blah]{Eggbox diagrams of the partition monoids $\P_3$ (left) and $\P_2$ (right), and the variant $\P_3^\si$ (middle), from Example \ref{e:max}\ref{emax8}.}
\label{f:emax8}
\end{center}
\end{figure}

\begin{eg}\label{e:max}
\ben

\item \label{emax1} Consider the full transformation semigroup $\T_n$, which consists of all functions from $\{1,\ldots,n\}$ to itself, and let $a\in\T_n$ be a non-bijection.  By \cite[Proposition 4.4]{DE2015}, all maximal $\J^a$-classes in $\T_n^a$ are trivial.  See also Proposition \ref{p:maxPPBB} below, which shows that this is true in some diagram categories as well.

\item \label{emax2} If $M$ is a finite monoid with group of units $G$, then $G$ is the unique maximal $\J$-class of $M$.  If $a\in G$, then $M^a\cong M$, so in fact $G$ is the unique maximal $\J^a$-class of $M^a$, and it is regular; thus, $M^a$ has no trivial maximal $\J^a$-classes.  (Note that every element of $M$ is $\leqJ$-below $a$.)

\item \label{emax4} Consider the semigroup $S=\{a,b,0\}$ with multiplication table given in Figure \ref{f:emax4}.  The multiplication table for the variant $S^a$ is also given in Figure \ref{f:emax4}.  Here $\{a\}$ and $\{b\}$ are both maximal $\J$-classes of $S$, and they are both also maximal $\J^a$-classes of $S^a$.  Here $\{b\}$ is a trivial $\J^a$-class, and $\{a\}$ is nontrivial; note that $a$ is regular (indeed, an idempotent), and is its own unique inverse.  Figure \ref{f:emax4} also gives eggbox diagrams of $S$ and $S^a$.

\item \label{emax5} Consider the semigroup $S=\{a,b,0\}$ with multiplication table given in Figure \ref{f:emax5}; note that $a$ is nonregular.  The multiplication table for the variant $S^a$ is also given in Figure \ref{f:emax5}, as well as eggbox diagrams for both $S$ and $S^a$.  Here $\{a\}$ is the unique maximal $\J$-class of~$S$, and it is not regular.  The maximal $\J^a$-classes of $S^a$ are $\{a\}$ and $\{b\}$, and these are both nontrivial (as every element of~$S$ is $\leqJ$-below $a$), and they are both nonregular.

\item \label{emax6} 
Let $S$ be the subsemigroup of $\T_3$ generated by the transformations $f=[2,1,2]$ and $a=[3,1,3]$.  (Here $[x_1,x_2,x_3]$ is the transformation mapping $i\mt x_i$ for $i=1,2,3$.)  GAP \cite{GAP} tells us that $S$ has size $7$, and it displays the eggbox diagrams of the semigroup $S$ and its variant $S^a$ as shown in Figure \ref{f:emax6}.  Here~$a$ is a nonregular element of $S$ (it belongs to the middle $\D$-class).  The variant $S^a$ has two trivial maximal $\J^a$-classes (corresponding to the two elements of the top $\D$-class of $S$, which are not $\leqJ$-below $a$) and also two non-trivial maximal $\J^a$-classes (corresponding to the two elements of the nonregular $\D$-class of $S$).

\item \label{emax7} Denote by $\TL_4$ the Temperley-Lieb monoid of degree $4$ (which is stable, regular and $\H$-trivial), as defined in Section \ref{ss:K}, and consider the (regular) partitions $\si = \custpartn{1,2,3,4}{1,2,3,4}{\stline11\stline22\uarc34\darc34}$ and $\tau = \custpartn{1,2,3,4}{1,2,3,4}{\stline13\stline24\uarc34\darc12}$ from~$\TL_4$.  Figure \ref{f:emax7} shows eggbox diagrams for~$\TL_4$ (and $\TL_2$, though this can be ignored for now), as well as the variants $\TL_4^\si$ and $\TL_4^\tau$, again all generated by GAP.  From this figure, it can be seen that $\TL_4^\tau$ has a nontrivial maximal $\J^\tau$-class, but $\TL_4^\si$ does not have a nontrivial maximal $\J^\si$-class; both have a single trivial maximal ${\J^\si}/{\J^\tau}$-class.  Note that $\si\J \tau$ (in $\TL_4$), and that the only element of $\TL_4$ strictly $\leqJ$-above $\si$ and $\tau$ is the identity element; this is a pre-inverse of~$\si$ but not of $\tau$ (since $\si$ is an idempotent and $\tau$ is not); cf.~Corollary~\ref{c:max} and Proposition \ref{p:max2}.

\item \label{emax8} Denote by $\P_3$ the partition monoid of degree $3$ (which is stable and regular), as defined in Section~\ref{ss:P}, and consider the (regular) partition $\si = \custpartn{1,2,3}{1,2,3}{\stline11\stline21\stline32\stline33\uarc12\darc23}\in\P_3$.  Figure \ref{f:emax8} shows eggbox diagrams for~$\P_3$ (and $\P_2$, though this can be ignored for now), as well as the variant $\P_3^\si$, again all generated by GAP.  From this figure, it is immediate that $\P_3^\si$ has only trivial maximal $\J^\si$-classes.  However, one may easily check (by hand, or using GAP) that each pre-inverse of $\si$ is $\J$-related to $\si$.  This shows that the converse of Corollary \ref{c:max} is not true in general.  Comparing this to Proposition~\ref{p:max2}, note that while $\P_3$ is stable (as it is finite, but see also Proposition \ref{p:stab}), it is not $\H$-trivial.  See also Remark \ref{r:max2}, which puts this example into a more general context.

\een
\end{eg}

\begin{figure}[ht]
\begin{center}
\begin{tikzpicture}
\node () at (0,0) {
\begin{tabular}{c|ccc}
$\cdot$ & $a$ & $b$ & $0$ \\
\hline
$a$ & $a$ & $0$ & $0$ \\
$b$ & $0$ & $b$ & $0$ \\
$0$ & $0$ & $0$ & $0$ \\
\end{tabular}
};
\node () at (3.3,0) {
\begin{tikzpicture}[scale=1,
block/.style={
draw,
fill=white,
rectangle, 
minimum width=.6cm,
minimum height=.6cm,
font=\small}]
\node[block,fill=black!20] (a) at (-1,1) {$a$};
\node[block,fill=black!20] (b) at (1,1) {$b$};
\node[block,fill=black!20] (0) at (0,0) {$0$};
\draw (a)--(0)--(b);
\end{tikzpicture}
};
\node () at (9,0) {
\begin{tabular}{c|ccc}
$\cdot$ & $a$ & $b$ & $0$ \\
\hline
$a$ & $a$ & $0$ & $0$ \\
$b$ & $0$ & $0$ & $0$ \\
$0$ & $0$ & $0$ & $0$ \\
\end{tabular}
};
\node () at (12.3,0) {
\begin{tikzpicture}[scale=1,
block/.style={
draw,
fill=white,
rectangle, 
minimum width=.6cm,
minimum height=.6cm,
font=\small}]
\node[block,fill=black!20] (a) at (-1,1) {$a$};
\node[block] (b) at (1,1) {$b$};
\node[block,fill=black!20] (0) at (0,0) {$0$};
\draw (a)--(0)--(b);
\end{tikzpicture}
};
\end{tikzpicture}
\caption{Multiplication table and eggbox diagram for $S$ (left) and $S^a$ (right), as in Example \ref{e:max}\ref{emax4}.}
\label{f:emax4}
\end{center}
\end{figure}

\subsection{Mid-identities and regularity-preserving elements}\label{ss:MIRP}

In Section \ref{ss:max} we identified a natural family of maximal $\J^a$-classes of $S_{ij}^a$, called the \emph{trivial} ones, namely those of the form $J_x^a=\{x\}$, where $x\in S_{ij}$ and $x\not\leqJ a$.  In the case that $a$ is regular, there is at most one nontrivial maximal $\J^a$-class, and this contains $\Pre(a)$, so is of the form $J_b^a$ for any $b\in V(a)$.  
Even in the case that $S_{ij}^a$ only has trivial maximal $\J^a$-classes, it is easy to see that~$V(a)$ is contained in a single $\J^a$-class of $S_{ij}^a$; we will see in this section that this $\J^a$-class is still very important, even if it is not maximal.  (But note that if~$S_{ij}^a$ only has trivial maximal $\J^a$-classes, then $\Pre(a)$ may or may not be contained in a single $\J^a$-class of~$S_{ij}^a$.)  The main result of this section (Proposition \ref{p:Jba2}) shows that under certain stability and regularity assumptions, the $\J^a$-class containing $V(a)$ has a very nice structure.  Before we state this result, we must first gather some notions from \cite{Yamada1955,Hickey1983}.

Let $T$ be a regular semigroup.  Recall \cite{Yamada1955} that an element $u$ of $T$ is a \emph{mid-identity} if $xy=xuy$ for all $x,y\in T$; for such a mid-identity $u$, the variant semigroup $T^u=(T,\star_u)$ is precisely the original semigroup~$T$.  Recall \cite{Hickey1983} that an element $u$ of $T$ is \emph{regularity-preserving} if the variant $T^u$ is regular.  We denote by $\MI(T)$ and $\RP(T)$ the (possibly empty) sets of all mid-identities and regularity-preserving elements of~$T$.  Clearly $\MI(T)\sub\RP(T)$.  
If $T$ is a monoid, then $\RP(T)$ is the group of units \cite[Proposition 1]{KL2001}, and $\MI(T)$ contains only the identity element.  It was argued in \cite{KL2001} that $\RP(T)$ is a good ``substitute'' for the group of units in a non-monoid.

The next result involves rectangular bands and groups; cf.~Section \ref{ss:RBG}.  Recall that $E(T)$ denotes the set of idempotents of a semigroup $T$.

\begin{prop}\label{p:Jba2}
Suppose $S$ is stable, and $a\in S_{ji}$ is regular.  Fix some $b\in V(a)$.  Then
\ben
\item \label{Jba21} $J_b^a=D_b^a$,
\item \label{Jba22} $E(J_b^a)=V(a)$ is a rectangular band (under $\star_a$),
\item \label{Jba23} $J_b^a$ is a rectangular group (under $\star_a$),
\item \label{Jba24} if $S$ is regular, then $J_b^a=\RP(P^a)$ and $E(J_b^a)=\MI(P^a)$.
\een
\end{prop}

\pf
\firstpfitem{\ref{Jba21}}  Since $S$ is stable, \cite[Lemma 2.6]{DE2018} gives ${\J}={\D}$ (in $S$).  It then follows from \cite[Corollary~2.5]{Sandwich1} that ${\D^a}={\J^a}$ (in $S_{ij}^a$).  The result is now clear.

\pfitem{\ref{Jba22} and \ref{Jba23}}  Suppose $x\in E(J_b^a)$.  Then $x=x\star_ax=xax$, and also $x\in J_b^a\sub J_b$; the latter gives $x\J b\J a$.  
It then follows from Lemma \ref{l:axa} that $a=axa$, and so $x\in V(a)$.

Conversely, suppose $x\in V(a)$.  Then $x=xax$ and $a=axa$.  Then $x=xax=xabax$ and $b=bab=baxab$, which gives $x\J^ab$: i.e., $x\in J_b^a$.  Since $x=xax=x\star_ax$, it follows that $x\in E(J_b^a)$.

Since $J_b^a$ is a regular $\D^a$-class, the proof will be complete if we can show that $E(J_b^a)$ is a subsemigroup of $S_{ij}^a$ (cf.~Lemma \ref{l:ED}).  But for any $x,y\in E(J_b^a)=V(a)$, it is easy to verify that $x\star_ay\in V(a)$.

\pfitem{\ref{Jba24}}  It was proved in \cite[Proposition 4.5]{Sandwich1} that $\MI(P^a)=V(a)$ holds under the weaker assumption that every element of $\{a\}\cup aS_{ij}a$ is regular (and without assuming stability of $S$); this condition on $a$ was called \emph{sandwich regularity} in \cite{Sandwich1}.  Combined with part \ref{Jba22}, this completes the proof of the second assertion.

It was shown in \cite[Theorem 1.2]{BH1984} that if $T$ is a regular semigroup with a mid-identity, then $\RP(T)$ consists of all elements of $T$ that are $\H$-related to a mid-identity.  Thus, since $P^a=\Reg(S_{ij}^a)$ is regular, and since $\MI(P^a)=V(a)\not=\emptyset$, it follows that
\[
\RP(P^a) = \bigcup_{x\in \MI(P^a)}H_x^a = \bigcup_{x\in E(D_b^a)}H_x^a = D_b^a = J_b^a. \qedhere
\]
\epf

\begin{rem}
One may readily locate the $\J^a$-class containing $V(a)$ in Figures \ref{f:emax7} and \ref{f:emax8}; see also Figures \ref{f:TL4} and \ref{f:B}.  The rectangular group structure of this $\J^a$-class is evident from the figures.
\end{rem}

\subsection{Right-invertibility}\label{ss:RI}

As in \cite[Section 2.2]{Sandwich1}, we say $a\in S_{ji}$ is \emph{right-invertible} if there exists $b\in S_{ij}$ such that $ab$ is a right-identity for $S_{ij}$: i.e., $x=xab$ for all $x\in S_{ij}$.  Such an element $b$ is called a \emph{right-inverse} of $a$, and it is then apparent that $b$ is right-identity for $S_{ij}^a$.  Note that~$ab$ need not be a right-identity for every~$x$ with $\br(x)=j$; see Remark \ref{r:SXY} for an example.  We write $\RI(a)$ for the set of all right-inverses of $a$.

Left-invertibility and left-inverses are defined analogously.  Every result in this section has a dual statement, but we will omit these.

\begin{lemma}\label{l:RI}
\ben
\item \label{RI1} If $a\in S_{ji}$ is right-invertible, then $V(a)=\Pre(a)\sub\RI(a)\sub\Post(a)$.
\item \label{RI2} If $a\in S_{ji}$ is right-invertible and regular, then $V(a)=\Pre(a)=\RI(a)\sub\Post(a)$.
\een
\end{lemma}

\pf
\firstpfitem{\ref{RI1}}
Let $x\in\Pre(a)$, so that $x\in S_{ij}$ and $a=axa$.  Since $a$ is right-invertible, we may fix some $b\in\RI(a)$.  Since $ab$ is a right-identity for $S_{ij}$, we have $x=x(ab)$.  But then $ax=ax(ab)=(axa)b=ab$ is a right-identity for $S_{ij}$, so that $x\in\RI(a)$.

Let $x\in\RI(a)$.  Since $ax$ is a right-identity for $S_{ij}$, and since $x\in S_{ij}$, we have $x=x(ax)$: i.e., $x\in\Post(a)$.

We have proved that $\Pre(a)\sub\RI(a)\sub\Post(a)$, and so also $V(a)=\Pre(a)\cap\Post(a)=\Pre(a)$.

\pfitem{\ref{RI2}}
In light of \ref{RI1}, it is enough to show that $\RI(a)\sub\Pre(a)$, so let $x\in\RI(a)$.  Since $a$ is regular, we have $a=aya$ for some $y\in S_{ij}$.  Since $ax$ is a right-identity for $S_{ij}$, it follows that $a=aya=a(yax)a=(aya)xa=axa$, so that $x\in\Pre(a)$.
\epf

\begin{rem}\label{r:SXY}
It is possible for an element of a partial semigroup to be right-invertible but not regular.  For example, let $X$ and $Y$ be distinct nonempty sets.  Let $I=\{X,Y\}$, and consider the partial semigroup 
\[
S = S_{X,X}\cup S_{X,Y}\cup S_{Y,Y}\cup S_{Y,X},
\]
where $S_{Y,X}$ consists of all partial maps $Y\to X$, $S_{Y,Y}$ consists of all partial maps $Y\to Y$, and where each other $S_{U,V}$ contains only the empty map $\emptyset_{U,V}:U\to V$.  If $a\in S_{Y,X}$ is nonempty, then $a\emptyset_{X,Y}=\emptyset_{Y,Y}$ is a right-identity for $S_{X,Y}=\{\emptyset_{X,Y}\}$, so that $a$ is right-invertible; however, $a$ is nonregular.  Note also that $\emptyset_{Y,Y}$ is not a right-identity for $S_{Y,Y}$.

On the other hand, if we wished to study some sandwich semigroup $S_{ij}^a$ with $a$ right-invertible, then we could assume without loss of generality that $a$ is regular.  Indeed, fix some right-inverse $b\in \RI(a)$, and define $c=aba\in S_{ji}$.  Then since $\RI(a)\sub\Post(a)$, we have $bab=b$, so ${cbc=abababa=aba=c}$, so that $c$ is regular; also, $c$ is still right-invertible ($b$ is a right-inverse).  Since $xay=xcy$ for all $x,y\in S_{ij}$, it follows that the sandwich semigroups $S_{ij}^a$ and $S_{ij}^c$ are precisely the same semigroup.
\end{rem}

Next we wish to prove a result on maximal $\J^a$-classes in the case that $a$ is right-invertible (Proposition \ref{p:max3} below).  To do this, we need the following simple lemma.

\begin{lemma}\label{l:stabuv}
If $T$ is a stable partial semigroup, and if $u,v\in T$ satisfy $u\leqL v \leqJ u$, then $u\L v$.
\end{lemma}

\pf
From $u\leqL v \leqJ u$ we have $u\J v$, and also $u=xv$ for some $x\in \Tone$.  But then $v\J u=xv$, so stability gives $v\L xv=u$.
\epf

The next result concerns left-groups, which were defined in Section \ref{ss:RBG}.

\begin{prop}\label{p:max3}
Suppose $a\in S_{ji}$ is right invertible.   
\ben
\item \label{max31} The sandwich semigroup $S_{ij}^a$ has a maximum $\J^a$-class, and this contains $\RI(a)$.
\item \label{max32} If $S_{ij}^a$ is stable, then the maximum $\J^a$-class of $S_{ij}^a$ is in fact an $\L^a$-class, and is a left-group with set of idempotents $\RI(a)$.
\een
\end{prop}

\pf
\firstpfitem{\ref{max31}}  Let $b\in\RI(a)$ be arbitrary.  It is enough to show that every element $x\in S_{ij}$ is $\leqJa$-below~$b$; but this is clear because $x=x(ab)=x\star_ab\leq_{\L^a} b$, whence $x\leqJa b$.

\pfitem{\ref{max32}}  Again, let $b\in\RI(a)$, so that the maximum $\J^a$-class of $S_{ij}^a$ is $J_b^a$.  By definition we have $L_b^a\sub J_b^a$.  Conversely, let $x\in J_b^a$.  As above, $x\leq_{\L^a} b$; since also $b\J^ax$, we have $x\leq_{\L^a}b\leqJa x$.  Lemma~\ref{l:stabuv} (applied in the semigroup $T=S_{ij}^a$) then gives $x\L^ab$: i.e., $x\in L_b^a$.

Since $J_b^a$ is an $\L^a$-class, it is certainly a $\D^a$-class.  Since it is regular (as $b\in J_b^a$ is an idempotent), it follows from Lemma \ref{l:EL} that it is a left-group.

Finally, we prove that $E(J_b^a)=\RI(a)$.  First, if $x\in\RI(a)$, then part \ref{max31} gives $x\in J_b^a$; since $ax$ is a right identity for $x\in S_{ij}$, $x=x(ax)=x\star_ax$, so that in fact $x\in E(J_b^a)$.  Conversely, suppose $x\in E(J_b^a)$.  So $x=xax$, and since $x\L^ab$ (as $x\in J_b^a=L_b^a$), we have $b=uax$ for some $u\in S_{ij}$.  But then for any $z\in S_{ij}$ we have $z = zab = zauax = zauaxax = zabax = zax$, so that $ax$ is a right-identity for~$S_{ij}$: i.e., $x\in\RI(a)$.
\epf

\begin{rem}
Stability is essential in Proposition \ref{p:max3}\ref{max32}.  For example, let $S=\P\T$ be the category of partial transformations, as studied in \cite[Section 3]{Sandwich2}.  Let $X$ be an infinite set, let ${a\in\P\T_X=\P\T_{XX}}$ be an injective mapping with domain $X$ that is not surjective.  Then $a$ is right-invertible, so we may fix some right-inverse ${b\in\RI(a)}$.  Now, $J_b^a$ is the maximum $\J^a$-class of $\P\T_X^a$ (by Proposition \ref{p:max3}\ref{max31}), but $J_b^a$ is not even a $\D^a$-class, let alone an $\L^a$-class.  (This all follows from \cite[Lemma 3.7]{Sandwich2} and the proof of \cite[Proposition~3.17]{Sandwich2}.)
\end{rem}

The next result is analogous to Proposition \ref{p:max3}, but concerns $\K$-classes of the hom-set $S_{ij}$ rather than $\K^a$-classes of the semigroup $S_{ij}^a$.  The proof is virtually identical.

\begin{lemma}\label{l:max4}
Suppose $a\in S_{ji}$ is right-invertible.
\ben
\item \label{max41} The hom-set $S_{ij}$ has a maximum $\J$-class, and this contains $\RI(a)$.
\item \label{max42} If $S$ is stable, then the maximum $\J$-class of $S_{ij}$ is in fact an $\L$-class.  \epfres
\een
\end{lemma}

The final result of this section is somewhat technical, but it will be used in Section \ref{ss:RankB} when studying sandwich semigroups in the Brauer category.  To give the statement, we first introduce some notation.  

Suppose $S$ is stable, and that $a\in S_{ji}$ is right-invertible; fix some right-inverse $b\in\RI(a)$.  By \cite[Proposition 2.14]{DE2018}, $S_{ij}^a$ is stable.
By Lemma \ref{l:max4}, $J_b=L_b$ is the maximum $\J$-class of the hom-set $S_{ij}$.  By Proposition \ref{p:max3}, $J_b^a=L_b^a$ is the maximum $\J^a$-class of the sandwich semigroup $S_{ij}^a$; moreover, $J_b^a$ is a left-group over $H_b^a$.  Let $X$ be a cross-section of the $\H$-classes contained in $J_b$, by which we mean that $J_b=\bigcup_{x\in X}H_x$, with $H_x\cap H_y=\emptyset$ for $x\not=y$.

Now, $J_b^a=L_b^a$ is a regular $\D^a$-class (since it is a left-group), so it follows from Theorem~\ref{t:Green_Sij} that $H_x^a=H_x$ for all $x\in J_b^a$.  Thus, we have $J_b^a=\bigcup_{x\in X_1}H_x=\bigcup_{x\in X_1}H_x^a$ for some $X_1\sub X$.  Put $X_2=X\sm X_1$.  Note that for $x\in X_2$, $H_x$ is an $\H$-class of $S$, but is not an $\H^a$-class of $S_{ij}^a$ (unless~$|H_x|=1$); cf.~Theorem \ref{t:Green_Sij}\ref{GS3}.

Recall that the \emph{rank} of a semigroup $T$, denoted $\rank(T)$, is the minimum size of a generating set for $T$.  In the next statement and proof, if $\Om\sub S_{ij}$, we write $\la\Om\ra$ for the subsemigroup of $S_{ij}^a$ generated by $\Om$: i.e., the set of all (nonempty) products of the form $y_1\star_ay_2\star_a\cdots\star_ay_k$, with $y_1,y_2,\ldots,y_k\in\Om$.  It follows from results of Ru\v{s}kuc \cite{Ruskuc1994} (see also \cite[Proposition 4.11]{Sandwich1}) that the rank of a left-group of degree $\lam$ over a group $G$ is equal to $\max(\lam,\rank(G))$.  In particular, since $J_b^a$ (as above) is a left-group of degree $|X_1|=|J_b^a/{\H^a}|$ over $H_b^a$, we have
\begin{equation}\label{e:rankJba}
\rank(J_b^a) = {\max} \big( |J_b^a/{\H^a}|, \rank(H_b^a)\big).
\end{equation}
(Here, if $\ve$ is an equivalence relation on a set $A$, and if $B\sub A$ is a union of $\ve$-classes, we write $B/\ve$ for the set of all such $\ve$-classes.)

\begin{prop}\label{p:RI}
Suppose $S$ is stable, and that $a\in S_{ji}$ is right-invertible.  Keep the above notation (the right-inverse $b$, the set $X=X_1\cup X_2$, etc.).  Let $T=\la J_b\ra$ be the subsemigroup of $S_{ij}^a$ generated by~$J_b$.
\ben
\item \label{pRI1}  We have $T=\la J_b^a\cup X_2\ra$.
\item \label{pRI2}  If $\rank(H_b^a)\leq|J_b^a/{\H^a}|$, then $\rank(T)=|J_b/{\H}|$.
\een
\end{prop}

\pf
\firstpfitem{\ref{pRI1}}  Since $J_b^a\cup X_2\sub J_b$, it suffices to show that $J_b\sub\la J_b^a\cup X_2\ra$, so let $y\in J_b$ be arbitrary.  Let $x\in X$ be such that $y\in H_x$.  In particular, $y\R x$, and so $y=xu$ for some $u\in\Sone$.  Since $x\in S_{ij}$, and since $ab$ is a right-identity for $S_{ij}$, we have $y=xu=x(ab)u = x\star_a(bu)$.  Since $x\in X\sub J_b^a\cup X_2$, the proof of \ref{pRI1} will be complete if we can show that $bu\in J_b^a$.  First note that $b\J y = xabu \leqJ bu \leqJ b$, so that all these elements are $\J$-related; in particular, $bu\J b$, so that $bu\in J_b=L_b$.  Since also $bu=b(ab)u=b\cdot a(bu)$ we have $bu\L a(bu)$: i.e., $bu\in P_2^a$.  But then it follows (using Theorem \ref{t:Green_Sij}\ref{GSij2}) that $bu\in L_b\cap P_2^a=L_b^a=J_b^a$, as required.

\pfitem{\ref{pRI2}}  By \eqref{e:rankJba}, and by assumption, we have $\rank(J_b^a)=|J_b^a/{\H^a}|=|X_1|$.  Fix a generating set $\Om$ of $J_b^a$ with $|\Om|=|X_1|$.  By \ref{pRI1}, we have $T=\la J_b^a\cup X_2\ra = \la\Om\cup X_2\ra$, so it follows that
\[
\rank(T)\leq|\Om\cup X_2|=|\Om|+|X_2|=|X|=|J_b/{\H}|.
\]
On the other hand, since $J_b$ is a maximal (indeed, maximum) $\J$-class in $S_{ij}$, \cite[Lemma 6.1(ii)]{Sandwich1} tells us that any generating set for $S_{ij}^a$ has size at least $|J_b/{\R}|$, so that $\rank(S_{ij}^a)\geq|J_b/{\R}|$.  But since $J_b$ is an $\L$-class, every $\R$-class contained in $J_b$ is in fact an $\H$-class, so $|J_b/{\R}|=|J_b/{\H}|$.
\epf

\subsection{Structure of the regular subsemigroup, and connections to nonsandwich semigroups}\label{ss:nonsandwich}

Throughout Section \ref{ss:nonsandwich} we assume the partial semigroup $S\equiv(S,I,\bd,\br,\cdot)$ is regular.  We also fix~${a\in S_{ji}}$, and an inverse $b\in V(a)$.  (If $S$ was a regular partial $*$-semigroup, we could take $b=a^*$ throughout.)  As in \cite[Section~3]{Sandwich1}, 
\bit
\item $(S_{ij}a,\cdot)$ is a subsemigroup of the (ordinary) semigroup $S_i$,
\item $(aS_{ij},\cdot)$ is a subsemigroup of the (ordinary) semigroup $S_j$,
\item $(aS_{ij}a,\star_{b})$ is a subsemigroup of the sandwich semigroup $S_{ji}^{b}$.
\eit
In fact, \cite[Proposition 3.5]{Sandwich1} tells us that $(aS_{ij}a,\star_{b})$ is a regular monoid with identity $a$, and that we have commutative diagrams of semigroup surmorphisms (i.e., surjective homomorphisms):
\begin{equation}\label{e:CD1}
\begin{tikzcd} 
~ & (S_{ij},\star_a) \arrow[swap]{dl}{x\mt xa} \arrow{dr}{x\mt ax} & \\
(S_{ij}a,\cdot) \arrow[swap]{dr}{y\mt ay} & & (aS_{ij},\cdot) \arrow{dl}{y\mt ya}\\
& (aS_{ij}a,\star_{b}) 
\end{tikzcd}
\qquad\quad
\begin{tikzcd} [column sep=small]
~ & \Reg(S_{ij},\star_a) \arrow[swap]{dl}{x\mt xa} \arrow{dr}{x\mt ax} & \\
\Reg(S_{ij}a,\cdot) \arrow[swap]{dr}{y\mt ay} & & \Reg(aS_{ij},\cdot) \arrow{dl}{y\mt ya}\\
& (aS_{ij}a,\star_{b}) 
\end{tikzcd}
\end{equation}
From $S_{ij}a = S_{ij}aba \sub S_i ba \sub S_{ij}a$, it follows that $S_{ij}a=S_iba$ is in fact the principal left ideal of $S_i$ generated by the idempotent $ba$.  Similarly, $aS_{ij} = abS_j$ is the principal right ideal of $S_j$ generated by the idempotent $ab$.  

Of particular importance is the surmorphism
\begin{equation}\label{e:phi}
\phi:\Reg(S_{ij}^a)\to (aS_{ij}a,\star_{b}):x\mt axa
\end{equation}
induced by the second diagram in \eqref{e:CD1}.  This map was a key ingredient in many results of \cite{Sandwich1} that demonstrate close structural relationships between the regular semigroup $P^a=\Reg(S_{ij}^a)$ and the regular monoid $(aS_{ij}a,\star_{b})$.

If $\K$ denotes any of Green's relations, then we define a relation $\gKh^a$ on $P^a$ by
\[
x\gKh^a y \iff x\phi\K^{b} y\phi \text{ in $(aS_{ij}a,\star_{b})$.}
\]
Then ${\gDh^a}={\D^a}$ and ${\gJh^a}={\J^a}$, while ${\K^a}\sub{\gKh^a}\sub{\D^a}$ if $\K$ is any of $\R$, $\L$ or $\H$; cf.~\cite[Lemma~3.11]{Sandwich1}.  For $x\in P^a$, we write $\Kh_x^a$ for the $\gKh^a$-class of $x$.  The following is \cite[Theorem~3.14]{Sandwich1}.  (Rectangular bands and groups were defined in Section \ref{ss:RBG}.)

\begin{thm}\label{t:RG}
Let $x\in P^a=\Reg(S_{ij}^a)$, and put $\rho=|\Hh_x^a/{\R^a}|$ and $\lam=|\Hh_x^a/{\L^a}|$.  Then
\ben
\item \label{RG1} the restriction to $H_x^a$ of the map $\phi:P^a\to (aS_{ij}a,\star_{b})$ is a bijection $\phi|_{H_x^a}:H_x^a\to H_{x\phi}^b$,
\item \label{RG2} $H_x^a$ is a group if and only if $H_{x\phi}^b$ is a group, in which case these groups are isomorphic, 
\item \label{RG3} if $H_x^a$ is a group, then $\Hh_x^a$ is a $\rho\times\lam$ rectangular group over $H_x^a$, 
\item \label{RG4} if $H_x^a$ is a group, then $E_a(\Hh_x^a)$ is a $\rho\times\lam$ rectangular band.  \epfres
\een
\end{thm}

\begin{rem}
One of the most important consequences of Theorem \ref{t:RG} is that a $\D^a$-class $D_x^a$ of $P^a=\Reg(S_{ij}^a)$ can be thought of as a kind of ``inflation'' of the $\D^b$-class $D_{x\phi}^b$ of $(aS_{ij}a,\star_b)$.  This is explained at length in \cite[Remark 3.15]{Sandwich1}, so we will not repeat the full details here.  But one may get an idea of this kind of ``inflation'' by examining Figures \ref{f:emax7}, \ref{f:emax8} and \ref{f:B}.
\end{rem}

As in \cite[Remark 3.6]{Sandwich1}, the maps
\begin{equation}\label{e:aSaa*}
aS_{ij}a \to aS_{ij}ab=abS_jab:x\mt xb \AND aS_{ij}a \to baS_{ij}a=baS_iba:x\mt bx
\end{equation}
induce isomorphisms $(aS_{ij}a,\star_{b})\to(abS_jab,\cdot)$ and $(aS_{ij}a,\star_{b})\to(baS_iba,\cdot)$.  In particular, $(aS_{ij}a,\star_{b})$ is isomorphic to both
\bit
\item $abS_jab$, the local submonoid of $S_j$ with identity $ab$, and 
\item $baS_iba$, the local submonoid of $S_i$ with identity $ba$.
\eit
In the diagram categories studied in Sections \ref{s:DC}--\ref{s:B}, the monoid $abS_jab$ will always be naturally isomorphic to an ordinary diagram monoid.

Combining \eqref{e:CD1} with the first isomorphism in \eqref{e:aSaa*}, we obtain the following diagrams, with all maps surmorphisms, and with all semigroups other than $S_{ij}^a$ and $\Reg(S_{ij}^a)$ having $\cdot$ as their operation:
\begin{equation}\label{e:CD2}
\begin{tikzcd} 
~ & S_{ij}^a \arrow[swap]{dl}{x\mt xa} \arrow{dr}{x\mt ax} & \\
S_iba \arrow[swap]{dr}{y\mt ayb} & & abS_j \arrow{dl}{y\mt yab}\\
& abS_jab
\end{tikzcd}
\qquad\quad
\begin{tikzcd} [column sep=small]
~ & \Reg(S_{ij}^a) \arrow[swap]{dl}{x\mt xa} \arrow{dr}{x\mt ax} & \\
\Reg(S_iba) \arrow[swap]{dr}{y\mt ayb} & & \Reg(abS_j) \arrow{dl}{y\mt yab}\\
& abS_jab
\end{tikzcd}
\end{equation}
The map
\begin{equation}\label{e:Phi}
\Phi:\Reg(S_{ij}^a)\to abS_jab:x\mt axab
\end{equation}
induced by the second diagram in \eqref{e:CD2} is the composition of $\phi$ (as in \eqref{e:phi}) with the first map in~\eqref{e:aSaa*}.  As such, the results from \cite[Section 3]{Sandwich1} may be rephrased in terms of the $\Phi$ map, including Theorem \ref{t:RG} above.

A final result from \cite{Sandwich1} that will be important to note concerns idempotents, and the idempotent-generated subsemigroup.  For a semigroup $T$, recall that $E(T)=\set{x\in T}{x=x^2}$ is the set of all idempotents from $T$; we will also write ${\E(T)=\la E(T)\ra}$ for the subsemigroup of $T$ generated by the idempotents.  Note that $E(S_{ij}^a)=\set{x\in S_{ij}}{x=xax}$ is precisely the set $\Post(a)$ of all post-inverses of $a$.  It was shown in \cite[Lemma~3.13 and Theorem~3.17]{Sandwich1} that 
\[
E(S_{ij}^a) = E(P^a) = E(aS_{ij}a,\star_{b})\phi^{-1} \AND \E(S_{ij}^a) = \E(P^a) = \E(aS_{ij}a,\star_{b})\phi^{-1}.
\]
(If $\psi:U\to V$ is a semigroup surmorphism, then $E(U)\sub E(V)\psi^{-1}$ and $\E(U)\sub \E(V)\psi^{-1}$ always hold, but the reverse inclusions do not in general.)
Combining these with the first isomorphism from~\eqref{e:aSaa*}, it of course quickly follows that 
\begin{equation}\label{e:ESija}
E(S_{ij}^a) = E(P^a) = E(abS_jab)\Phi^{-1} \AND \E(S_{ij}^a) = \E(P^a) = \E(abS_jab)\Phi^{-1}.
\end{equation}

\section{Diagram categories}\label{s:DC}

We now wish to apply the results of Sections \ref{s:S} and \ref{s:sandwich} to a number of diagram categories.  In this section we recall the definitions of these categories, and then prove a number of structural results concerning them.  Section \ref{s:sandwichK} treats sandwich semigroups in these categories; Section \ref{s:B} concerns the particular case of the Brauer category, which turns out to be especially amenable to analysis.

\subsection{The partition category}\label{ss:P}

The diagram categories we wish to study are all subcategories of the \emph{partition categories}, so we define these first.  We write $\N=\{0,1,2,\ldots\}$ for the set of all natural numbers.  For $n\in\N$, we write $[n]=\{1,\ldots,n\}$, interpreting $[0]=\emptyset$.  For $A\sub\N$, we write $A'=\set{a'}{a\in A}$ and $A''=\set{a''}{a\in A}$.  For $m,n\in\N$, we denote by $\Pmn$ the set of all set partitions of $[m]\cup[n]'$, and we write
\[
\P = \bigcup_{m,n\in\N}\Pmn
\]
for the set of all such set partitions.  For $m,n\in\N$ and $\al\in\Pmn$, we write $\bd(\al)=m$ and $\br(\al)=n$.  Then as in \cite{Martin2008},
\[
\P \equiv (\P,\N,\bd,\br,\cdot,*)
\]
is a regular $*$-category, under operations $\cdot$ and $*$ to be defined shortly.  

First we recall some terminology.  We say a nonempty subset $X\sub\N\cup\N'$ is
\bit
\item a \emph{transversal} if $X\cap \N$ and $X\cap \N'$ are both nonempty,
\item an \emph{upper nontransversal} if $X\sub \N$,
\item a \emph{lower nontransversal} if $X\sub \N'$.
\eit
If $\al\in\P$, we will write
\begin{equation}\label{eq:partn}
\al = \partI{A_1}{A_r}{C_1}{C_s}{B_1}{B_r}{D_1}{D_t}
\end{equation}
to indicate that $\al$ has transversals $A_i\cup B_i'$ ($1\leq i\leq r$), upper nontransversals $C_i$ ($1\leq i\leq s$) and lower nontransversals $D_i'$ ($1\leq i\leq t$).  This notation uniquely determines $\bd(\al)$ and $\br(\al)$, since 
\[
[\bd(\al)] = \bigcup_{i=1}^rA_i \cup \bigcup_{i=1}^sC_i \AND [\br(\al)] = \bigcup_{i=1}^rB_i \cup \bigcup_{i=1}^tD_i.
\]
Note that any of $q,r,s$ could be zero in \eqref{eq:partn}; they could even all be zero, in which case $\al=\emptyset$ is the unique element of $\P_{0,0}$.

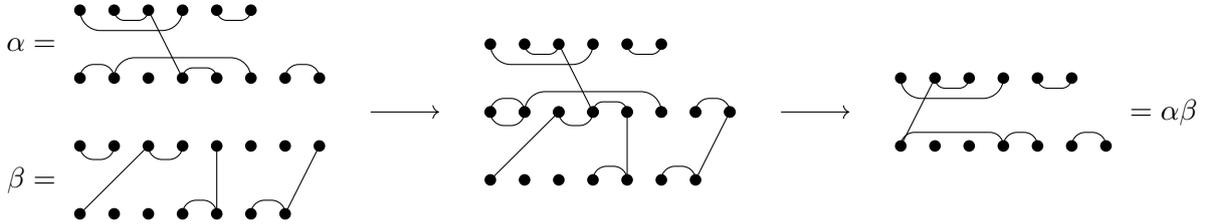
\begin{figure}[ht]
\begin{center}
\begin{tikzpicture}[scale=.45]
\begin{scope}[shift={(0,0)}]	
\uvs{1,...,6}
\lvs{1,...,8}
\uarcx14{.6}
\uarcx23{.3}
\uarcx56{.3}
\darc12
\darcx26{.6}
\darcx45{.3}
\darc78
\stline34
\draw(0.6,1)node[left]{$\al=$};
\draw[->](9.5,-1)--(11.5,-1);
\end{scope}
\begin{scope}[shift={(0,-4)}]	
\uvs{1,...,8}
\lvs{1,...,7}
\uarc12
\uarc34
\darc45
\darc67
\stline31
\stline55
\stline87
\draw(0.6,1)node[left]{$\be=$};
\end{scope}
\begin{scope}[shift={(12,-1)}]	
\uvs{1,...,6}
\lvs{1,...,8}
\uarcx14{.6}
\uarcx23{.3}
\uarcx56{.3}
\darc12
\darcx26{.6}
\darcx45{.3}
\darc78
\stline34
\draw[->](9.5,0)--(11.5,0);
\end{scope}
\begin{scope}[shift={(12,-3)}]	
\uvs{1,...,8}
\lvs{1,...,7}
\uarc12
\uarc34
\darc45
\darc67
\stline31
\stline55
\stline87
\end{scope}
\begin{scope}[shift={(24,-2)}]	
\uvs{1,...,6}
\lvs{1,...,7}
\uarcx14{.6}
\uarcx23{.3}
\uarcx56{.3}
\darc14
\darc45
\darc67
\stline21
\draw(7.4,1)node[right]{$=\al\be$};
\end{scope}
\end{tikzpicture}
\caption{Left to right: partitions $\al\in\P_{6,8}$ and $\be\in\P_{8,7}$, the product graph $\Pi(\al,\be)$, and the product $\al\be\in\P_{6,7}$.}
\label{fig:P}
\end{center}
\end{figure}

As usual, we identify a partition $\al\in\Pmn$ with any graph on vertex set $[m]\cup[n]'$ whose connected components are the blocks of $\al$.  When depicting such a graph in the plane $\RR^2$, we always draw
\bit
\item vertex $x$ at $(x,1)$ for each $1\leq x\leq m$, 
\item vertex $x'$ at $(x,0)$ for each $1\leq x\leq n$,
\item all edges in the rectangle $\bigset{(x,y)\in\mathbb R^2}{1\leq x\leq \max(m,n),\ 0\leq y\leq1}$.
\eit
For example, graphs corresponding to the partitions
\begin{align}
\label{al} \al &= \big\{ \{1,4\},\{2,3,4',5'\},\{5,6\},\{1',2',6'\},\{3'\},\{7',8'\}\big\} \in \P_{6,8} ,\\
\label{be} \be &= \big\{ \{1,2\}, \{3,4,1'\}, \{5,4',5'\}, \{6\}, \{7\}, \{8,6',7'\}, \{2'\}, \{3'\} \big\} \in \P_{8,7}
\end{align}
are pictured in Figure \ref{fig:P}.

Now let $\al\in\Pmn$ and $\be\in\Pnk$.  We define $\al_\downarrow$ to be the graph on vertex set $[m]\cup[n]''$ obtained by changing each lower vertex $x'$ of (a graph representing) $\al$ to $x''$; similarly, we define $\be^\uparrow$ to be the graph on vertex set $[n]''\cup[k]'$ obtained by changing each upper vertex $x$ of $\be$ to $x''$.  The \emph{product graph} $\Pi(\al,\be)$ is the graph on vertex set $[m]\cup[n]''\cup[k]'$ whose edge set is the union of the edge sets of $\al_\downarrow$ and $\be^\uparrow$.  The product $\al\be\in\Pmk$ is the partition such that $x,y\in[m]\cup[k]'$ belong to the same block of $\al\be$ if and only if $x,y$ are in the same connected component of $\Pi(\al,\be)$.  As an example, Figure \ref{fig:P} shows how to calculate the product
\begin{align*}
\al\be &= \big\{ \{1,4\},\{2,3,1',4',5'\},\{5,6\},\{2'\},\{3'\},\{6',7'\}\big\} \in \P_{6,7} ,
\end{align*}
where $\al\in \P_{6,8}$ and $\be\in \P_{8,7}$ are as in \eqref{al} and \eqref{be}.

The involution $*:\P\to\P$ is easier to describe.  For $\al\in\P$ as in \eqref{eq:partn}, we define
\[
\al^*=\partI{B_1}{B_r}{D_1}{D_t}{A_1}{A_r}{C_1}{C_s}.
\]
Diagrammatically, $\al^*$ is the result of reflecting (a graph representing) $\al$ in a horizontal axis.

The endomorphism monoids in $\P$ are the partition monoids $\P_n=\P_{nn}$ ($n\in\N$).  The identity of~$\P_n$ is the partition $\id_n = \bigset{\{x,x'\}}{x\in[n]}$, with respect to which the invertible elements of $\P_n$ are those of the form $\bigset{\{x,(x\pi)'\}}{x\in[n]}$, for some permutation $\pi\in\S_n$.  Such a partition will be identified with $\pi$ itself, so that the automorphism groups in $\P$ are (identified with) the symmetric groups $\S_n$.

\newpage

\subsection{Subcategories}\label{ss:K}

A partition $\al\in\Pmn$ is called
\bit
\item a \emph{Brauer partition} if each block of $\al$ has size $2$,
\item a \emph{partial Brauer partition} if each block of $\al$ has size $1$ or $2$.
\eit
We denote by $\B$ and $\PB$ the sets of all Brauer and partial Brauer partitions, respectively.  These are both subcategories of $\P$, and are both closed under the $*$ map; in particular, they are both regular $*$-categories, and they are called the \emph{Brauer} and \emph{partial Brauer categories} \cite{LZ2015,MarMaz2014}.

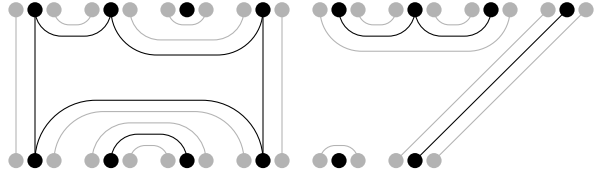
\begin{figure}[ht]
\begin{center}
\begin{tikzpicture}[x={(1,0)},y={(0,-1)}]
\uarcxx{1.25}{3.75}{.65}{gray!60}
\uarcxx{1.75}{3.25}{.5}{gray!60}
\uarcxx{2.25}{2.75}{.2}{gray!60}
\uarcxx{4.75}{5.25}{.2}{gray!60}
\uarcx14{.8}
\uarcx23{.35}
\darcxx{2.75}{3.25}{.2}{gray!60}
\darcxx{1.25}{1.75}{.2}{gray!60}
\darcxx{5.25}{5.75}{.2}{gray!60}
\darcxx{6.25}{6.75}{.2}{gray!60}
\darcxx{2.25}{3.75}{.4}{gray!60}
\darcxx{4.75}{7.25}{.55}{gray!60}
\stlinex{.75}{.75}{gray!60}
\stlinex{4.25}{4.25}{gray!60}
\stlinex{5.75}{7.75}{gray!60}
\stlinex{6.25}{8.25}{gray!60}
\darcx12{.35}
\darcx24{.6}
\darcx56{.35}
\darcx67{.35}
\stlines{1/1,4/4,6/8}
\foreach \x in {1,...,8} {
\fill (\x,0)circle(.1);
\fill[gray!60] (\x+.25,0)circle(.1);
\fill[gray!60] (\x-.25,0)circle(.1);
}
\foreach \x in {1,...,6} {
\fill (\x,2)circle(.1);
\fill[gray!60] (\x+.25,2)circle(.1);
\fill[gray!60] (\x-.25,2)circle(.1);
}
\end{tikzpicture}
\caption{A planar partition $\al$ from $\PP_{8,6}$ (black), with its corresponding Temperley-Lieb partition~$\alt$ from $\TL_{16,12}$ (grey).}
\label{f:TLiso}
\end{center}
\end{figure}

As in \cite{Jones1994_2}, a partition $\al\in\Pmn$ is \emph{planar} if (some graph representing) $\al$ can be drawn in the plane~$\RR^2$, with vertices in the locations specified above, edges in the specified rectangle, and with no edge crossings in the interior of this rectangle.
For example, with $\al,\be$ as in \eqref{al} and \eqref{be}, $\be$ is planar but $\al$ is not; cf.~Figure \ref{fig:P}.  
To each planar partition $\al$, we may associate a \emph{canonical (planar) graph}, as shown by example in Figure \ref{f:TLiso} (see the black edges, and ignore the grey edges for now); for a more  explicit definition, see \cite[Section 7]{EMRT2018}.

\begin{figure}[ht]
\begin{center}
\begin{tikzpicture}[scale=1.05]
\node[rounded corners,rectangle,draw,fill=blue!25] (P) at (0,8) {$\P$};
\node[rounded corners,rectangle,draw,fill=green!25] (PB) at (0,6) {$\PB$};
\node[rounded corners,rectangle,draw,fill=red!25] (B) at (0,4) {$\B$};
\node[rounded corners,rectangle,draw,fill=blue!15] (PP) at (3,6) {$\PP$};
\node[rounded corners,rectangle,draw,fill=green!15] (M) at (3,4) {$\M$};
\node[rounded corners,rectangle,draw,fill=red!15] (TL) at (3,2) {$\TL$};
\draw (M)--(PP)--(P)--(PB)--(B)--(TL)--(M)--(PB);
\begin{scope}[shift={(8,0)}]
\node[rounded corners,rectangle,draw,fill=blue!25] (P) at (0,8) {$\custpartn{1,2,3,4}{1,2,3,4,5,6}{\uarcx14{.6}\uarcx23{.3}\darc12\darcx26{.6}\darcx45{.3}\stline34}$};
\node[rounded corners,rectangle,draw,fill=green!25] (PB) at (0,6) {$\custpartn{1,2,3,4}{1,2,3,4,5,6}{\uarcx13{.5}\darc45\stline23}$};
\node[rounded corners,rectangle,draw,fill=red!25] (B) at (0,4) {$\custpartn{1,2,3,4}{1,2,3,4,5,6}{\stline13\stline42\uarc23\darcx16{.8}\darc45}$};
\node[rounded corners,rectangle,draw,fill=blue!15] (PP) at (3,6) {$\custpartn{1,2,3,4}{1,2,3,4,5,6}{\uarcx14{.6}\uarcx23{.2}\darc45\darc12\darc23\stline43\stline11}$};
\node[rounded corners,rectangle,draw,fill=green!15] (M) at (3,4) {$\custpartn{1,2,3,4}{1,2,3,4,5,6}{\uarcx13{.5}\darc45\stline43}$};
\node[rounded corners,rectangle,draw,fill=red!15] (TL) at (3,2) {$\custpartn{1,2,3,4}{1,2,3,4,5,6}{\uarcx12{.4}\darc34\darcx25{.8}\stline31\stline46}$};
\draw (M)--(PP)--(P)--(PB)--(B)--(TL)--(M)--(PB);
\end{scope}
\end{tikzpicture}
\end{center}
\vspace{-5mm}
\caption{Subcategories of $\P$ (left) and representative elements from each (right).}
\label{fig:subcats}
\end{figure}
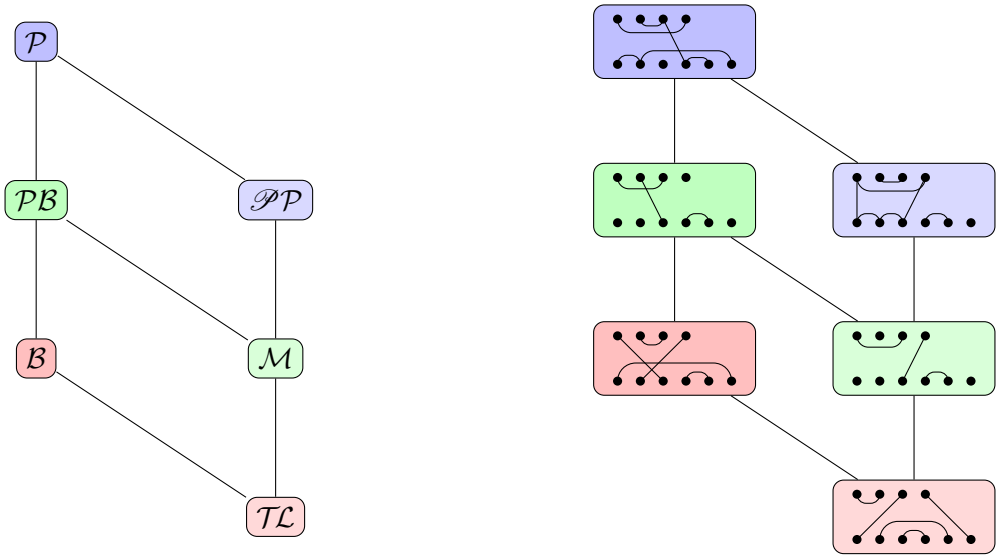

If $\al$ and $\be$ are planar, and satisfy $\br(\al)=\bd(\be)$, then the product $\al\be$ is planar as well (consider the product graph $\Pi(\al,\be)$).  It follows that the set
\[
\PP = \set{\al\in\P}{\al \text{ is planar}}
\]
is a subcategory of $\P$.  We also have corresponding planar subcategories of $\B$ and $\PB$.  These are the \emph{Temperley-Lieb} and \emph{Motzkin} categories:
\[
\TL = \B\cap\PP \AND \M = \PB\cap\PP.
\]
Again, the categories $\PP$, $\TL$ and $\M$ are all regular $*$-categories.
Figure \ref{fig:subcats} shows the relative inclusions of the above subcategories of $\P$, as well as representative elements of each.
Endomorphism monoids in these categories are Brauer monoids $\B_n$, partial Brauer monoids $\PB_n$, planar partition monoids~$\PP_n$, Motzkin monoids $\M_n$, and Temperley-Lieb monoids $\TL_n$ (also known as Jones monoids $\mathcal J_n$).  Automorphism groups in $\B$ and $\PB$ are (identified with) symmetric groups $\S_n$; automorphism groups in~$\PP$,~$\M$ and~$\TL$ are trivial.

Before moving on, we note that the planar partition category $\PP$ is isomorphic to the ``even subcategory'' $\TL^\even$ of $\TL$ defined by
\[
\TL^\even = \bigcup_{m,n\in\N}\TL_{2m,2n}.
\]
Indeed, an isomorphism $\PP\to\TL^\even:\al\mt\alt$ may be constructed as follows.  Given $\al\in\PP_{mn}$, we represent $\al$ by its canonical graph (as in \cite[Section 7]{EMRT2018}), and then construct $\alt\in\TL_{2m,2n}$ by ``tracing around'' the edges of $\al$, as shown in Figure \ref{f:TLiso}.  This observation is inspired by similar facts about planar partition and Temperley-Lieb \emph{monoids}; see for example \cite[Section 1]{HR2005}.

Next we wish to give formulae for the sizes of hom-sets in our diagram categories.
To do so, and for later use, we fix some notation for certain well-known number sequences.  More information about these numbers (including the combinatorial descriptions given below) can be found at the quoted entries on the Online Encyclopedia of Integer Sequences \cite{OEIS}.
\bit

\item The Stirling numbers of the second kind, $S(n,k)$, are defined by $S(n,1)=S(n,n)=1$ and $S(n,k)=S(n-1,k-1)+kS(n-1,k)$ for $1<k<n$; cf.~\cite[A008277]{OEIS}.  Note that $S(n,k)$ is the number of partitions of an $n$-element set into $k$ blocks.  We also define $S(0,0)=1$.

\item The $n$th Bell number is $B(n)=\sum_{k=1}^nS(n,k)$; cf.~\cite[A000110]{OEIS}.  Note that $B(n)$ is the number of partitions of an $n$-element set.  We also define $B(0)=1$.

\item The double factorial, $n!!$, is defined to be $0$ if $n$ is even, or $n\cdot(n-2)\cdots3\cdot1$ if $n$ is odd; cf.~\cite[A123023]{OEIS}.  By convention, we define $(-1)!!=1$.  Note that number of partitions of an $n$-element set into blocks of size $2$ is $(n-1)!!$.

\item The numbers $a(n)$ are defined by $a(0)=a(1)=1$ and $a(n)=a(n-1)+(n-1)a(n-2)$ for $n\geq2$; cf.~\cite[A000085]{OEIS}.  Note that $a(n)$ is the number of partitions of an $n$-element set into blocks of size $\leq2$ (and also the number of involutions of an $n$-element set).

\item The $n$th Catalan number is $C(n)=\frac1{n+1}\binom{2n}n$; cf.~\cite[A000108]{OEIS}.  Note that $C(n)$ is the number of non-crossing partitions of an $n$-element set, and also the number of non-crossing partitions of a $2n$-element set into blocks of size $2$.  By convention, we define $C(x)=0$ if $x$ is not a natural number.

\item The Motzkin triangle numbers, $\mu(n,k)$, are defined by $\mu(0,0)=1$, $\mu(n,k)=0$ if $n<k$ or $k<0$, and $\mu(n,k)=\mu(n-1,k-1)+\mu(n-1,k)+\mu(n-1,k+1)$ if $n\geq1$ and $0\leq k\leq n$; cf.~\cite[A064189]{OEIS}.  (It is not as easy to give a succinct combinatorial description of this sequence, but the above OEIS entry gives a few.)

\item The $n$th Motzkin number is $\mu(n)=\mu(n,0)$; cf.~\cite[A001006]{OEIS}.  Note that $\mu(n)$ is the number of non-crossing partitions of an $n$-element set into blocks of size $\leq2$.

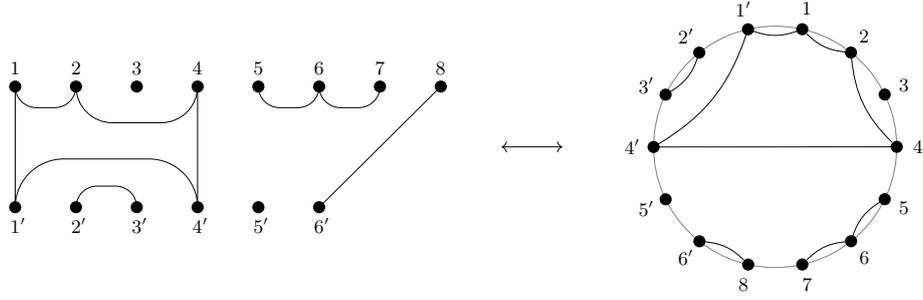
\begin{figure}[ht]
\begin{center}
\scalebox{0.8}{
\begin{tikzpicture}
\begin{scope}[x={(1,0)},y={(0,-1)}]
\uarcx14{.8}
\uarcx23{.35}
\darcx12{.35}
\darcx24{.6}
\darcx56{.35}
\darcx67{.35}
\stlines{1/1,4/4,6/8}
\foreach \x in {1,...,8} {\fill (\x,0)circle(.1); \node () at (\x,-0.3) {\footnotesize $\x$};}
\foreach \x in {1,...,6} {\fill (\x,2)circle(.1); \node () at (\x,2.3) {\footnotesize $\phantom{'}\x'$};}
\draw[<->] (9,1)--(10,1);
\end{scope}
\begin{scope}[shift={(13.5,-1)},inner sep=0.7mm]
\draw[gray] (0,0) circle (2cm);
\foreach \x in {1,...,14} {\fill ({90+12.86-25.71*\x}:2)circle(.1);}
\foreach \x in {1,...,8} {\node () at ({90+12.86-25.71*\x}:2.35) {\footnotesize $\x$};}
\foreach \x in {1,...,6} {\node () at ({90-12.86+25.71*\x}:2.35) {\footnotesize $\x'$};}
\foreach \x in {1,...,14} {\node (\x) at ({90+12.86-25.71*\x}:2) {};}
\draw (1) to [bend right = 20] (2);
\draw (2) to [bend right = 20] (4);
\draw (4) -- (11);
\draw (11) to [bend right = 20] (14);
\draw (14) to [bend right = 20] (1);
\draw (12) to [bend right = 20] (13);
\draw (5) to [bend right = 20] (6);
\draw (6) to [bend right = 20] (7);
\draw (8) to [bend right = 20] (9);
\end{scope}
\end{tikzpicture}
}
\caption{A planar partition $\al$ from $\PP_{8,6}$, with its corresponding non-crossing partition of $[8]\cup[6]'$.}
\label{f:NC}
\end{center}
\end{figure}

\eit
The next result follows quickly from the combinatorial descriptions of the above number sequences.  Note that elements of the planar categories correspond to non-crossing partitions in the way shown in Figure \ref{f:NC}.

\begin{prop}\label{p:sizes}
If $m,n\in\N$, then
\begin{align*}
|\Pmn| &= B(m+n) , & 
|\PB_{mn}| &= a(m+n) ,  &  
|\Bmn| &= (m+n-1)!! , \\
|\PP_{mn}| &= C(m+n) , & 
|\M_{mn}| &= \mu(m+n) , & 
 |\TL_{mn}| &= C(\tfrac{m+n}2) .  \epfreseq
\end{align*}
\end{prop}

\subsection{Green's relations and stability}\label{ss:greenK}

Before we even describe Green's relations in our diagram categories, we may quickly establish stability (see \eqref{e:stab} for the definition) using results of \cite{Sandwich1}.  
First, we recall that an element $u\in S$ is stable if the implications in~\eqref{e:stab} hold for all $x\in S$.  
Next, we recall that a semigroup $T$ is \emph{periodic} if for all $x\in T$, some power of $x$ is an idempotent; for example, it is well known that every finite semigroup is periodic; see \cite[Proposition~1.2.3]{Howie}.

\begin{prop}\label{p:stab}
The categories $\P$, $\PB$, $\B$, $\PP$, $\M$ and $\TL$ are all stable.  In particular, ${\J}={\D}$ in each of these categories.
\end{prop}

\pf
Let $\KK$ denote any of the stated categories.  Let $\al\in\KK$, say with $\al\in\KK_{mn}$.  Since $\al\KK_{nm}$ and~$\KK_{nm}\al$ are finite semigroups (they are subsemigroups of the finite monoids $\KK_m$ and $\KK_n$, respectively), they are periodic.  It then follows from \cite[Lemma 2.3]{Sandwich1} that $\al$ is stable.  Since $\al\in\KK$ was arbitrary, it follows that $\KK$ is itself stable.  The second assertion follows from \cite[Lemma~2.6]{DE2018}, which states that ${\J}={\D}$ in any stable partial semigroup.
\epf

Green's relations on our diagram categories may be described in terms of a number of parameters, defined as follows.  The \emph{(co)domain} and \emph{(co)kernel} of a partition $\al\in\Pmn$ are
\begin{align*}
\dom(\al) &= \set{x\in[m]}{x\text{ belongs to a transversal of }\al},\\
\codom(\al) &= \set{x\in[n]}{x'\text{ belongs to a transversal of }\al},\\
\ker(\al) &= \set{(x,y)\in[m]\times[m]}{x,y \text{ belong to the same block of }\al},\\
\coker(\al) &= \set{(x,y)\in[n]\times[n]}{x',y' \text{ belong to the same block of }\al}.
\end{align*}
The \emph{rank} of $\al$, denoted $\rank(\al)$, is the number of transversals of $\al$.  Note that $\dom(\al)$ and $\codom(\al)$ are subsets of $[m]$ and $[n]$, while $\ker(\al)$ and $\coker(\al)$ are equivalences on $[m]$ and $[n]$, and $\rank(\al)$ is an integer between $0$ and $\min(m,n)$.  
Following \cite{FL2011}, we write $N_U(\al)$ for the set of upper nontransversals of~$\al$, and $N_L(\al)$ for the set of lower nontransversals.
For $\al = \partI{A_1}{A_r}{C_1}{C_s}{B_1}{B_r}{D_1}{D_t}$, we have ${\rank(\al)=r}$, and
\begin{align*}
\dom(\al) &= \bigcup_{i=1}^rA_i ,&   N_U(\al) &= \set{C_i}{1\leq i\leq s} ,&  [m]/\ker(\al) &= \set{A_i}{1\leq i\leq r}\cup\set{C_i}{1\leq i\leq s},\\
\codom(\al) &= \bigcup_{i=1}^rB_i ,&   N_L(\al) &= \set{D_i}{1\leq i\leq t} ,&  [n]/\coker(\al) &= \set{B_i}{1\leq i\leq r}\cup\set{D_i}{1\leq i\leq t}.
\end{align*}

We have a number of obvious identities such as the following (and their duals), which are valid whenever the stated products are defined:
\begin{align}
\label{id1} \dom(\al\be)&\sub\dom(\al), & \dom(\al^*)&=\codom(\al), & N_U(\al\be) \supseteq N_U(\al),\\
\label{id2} \ker(\al\be)&\supseteq\ker(\al), & \ker(\al^*)&=\coker(\al), & \rank(\al\be\ga)\leq\rank(\be).
\end{align}

Before we give the characterisation of Green's relations and preorders on our diagram categories (Theorem \ref{t:G}), we first note a property of planarity from \cite{EMRT2018}.
Let $A=\{a_1,\ldots,a_k\}$ and ${B=\{b_1,\ldots,b_l\}}$ be nonempty finite subsets of $\N$ with $a_1<\cdots<a_k$ and $b_1<\cdots<b_l$.  We say that $A$ and $B$ are \emph{separated} if $a_k<b_1$ or $b_l<a_1$; in these cases, we write $A<B$ or $B<A$, respectively.  We say that $A$ is \emph{nested} by $B$ if there exists some $1\leq i<l$ such that $b_i<a_1$ and $a_k<b_{i+1}$; we say that $A$ and $B$ are \emph{nested} if $A$ is nested by $B$ or vice versa.  

\begin{lemma}\label{l:PP}
Let $\al=\partI{A_1}{A_r}{C_1}{C_s}{B_1}{B_r}{D_1}{D_t}\in\P$, with $\min(A_1)<\cdots<\min(A_r)$.  Then $\al$ is planar if and only if the following hold:
\bit
\item $A_1<\cdots<A_r$ and $B_1<\cdots<B_r$,
\item for all $1\leq i<j\leq s$, $C_i$ and $C_j$ are either nested or separated,
\item for all $1\leq i<j\leq t$, $D_i$ and $D_j$ are either nested or separated,
\item for all $1\leq i\leq r$ and $1\leq j\leq s$, either $A_i$ and $C_j$ are separated or else $C_j$ is nested by $A_i$, 
\item for all $1\leq i\leq r$ and $1\leq j\leq t$, either $B_i$ and $D_j$ are separated or else $D_j$ is nested by $B_i$.  
\eit
\end{lemma}

\pf
The forwards implication is essentially \cite[Lemma 7.1]{EMRT2018}; the converse is clear.
\epf

We also note an obvious parity issue involving elements of $\B$ (and $\TL$), which follows from the fact that all blocks of Brauer partitions have size $2$.  If $m,n\in\N$, then
\[
\B_{mn} \not=\emptyset \iff \TL_{mn} \not= \emptyset \iff m\equiv n\Mod2,
\]
in which case any $\al\in\B_{mn}$ satisfies $\rank(\al)\equiv m\equiv n\Mod2$.

Here is the result describing Green's relations and preorders on our diagram categories.  For the corresponding statements on diagram monoids, see \cite{Wilcox2007,DEG2017,FL2011,LF2006}.

\begin{thm}\label{t:G}
Let $\KK$ denote any of the categories $\P$, $\PB$, $\B$, $\PP$, $\M$ or $\TL$.  If $\al,\be\in\KK$, then in the category $\KK$ we have
\ben
\item \label{G1} $\al\leqR\be \iff \bd(\al)=\bd(\be)$, $\ker(\al)\supseteq\ker(\be)$ and $N_U(\al)\supseteq N_U(\be)$,
\item \label{G2} $\al\leqL\be \iff \br(\al)=\br(\be)$, $\coker(\al)\supseteq\coker(\be)$ and $N_L(\al)\supseteq N_L(\be)$,
\item \label{G3} $\al\leqJ\be \iff \begin{cases} \rank(\al)\leq\rank(\be) &\text{if $\KK$ is not $\B$ or $\TL$}\\ \rank(\al)\leq\rank(\be)\text{ and }\rank(\al)\equiv\rank(\be) \Mod2 &\text{if $\KK$ is one of $\B$ or $\TL$,}\end{cases}$
\item \label{G4} $\al\R\be \iff \ker(\al)=\ker(\be)$ and $N_U(\al)= N_U(\be)$
\item[] $\phantom{\al\R\be} \iff \dom(\al)=\dom(\be)$ and $\ker(\al)=\ker(\be)$,
\item \label{G5} $\al\L\be \iff \coker(\al)=\coker(\be)$ and $N_L(\al)= N_L(\be)$
\item[] $\phantom{\al\L\be} \iff \codom(\al)=\codom(\be)$ and $\coker(\al)=\coker(\be)$,
\item \label{G6} $\al\J\be \iff \al\D\be \iff \rank(\al)=\rank(\be)$.
\een
\end{thm}

\pf
Throughout the proof we write $\al=\partI{A_1}{A_r}{C_1}{C_s}{B_1}{B_r}{D_1}{D_t}$ and $\be=\partI{E_1}{E_q}{G_1}{G_u}{F_1}{F_q}{H_1}{H_v}$, where $\min(A_1)<\cdots<\min(A_r)$ and $\min(E_1)<\cdots<\min(E_q)$, and we also assume that $\al\in\Pmn$ and $\be\in\Pkl$.

\pfitem{\ref{G1}}  Suppose first that $\al\leqR\be$, so that $\al=\be\ga$ for some $\ga\in\P$.  Then $\bd(\al)=\bd(\be\ga)=\bd(\be)$.  By~\eqref{id1} and \eqref{id2}, we also have $\ker(\al)=\ker(\be\ga)\supseteq\ker(\be)$ and $N_U(\al)=N_U(\be\ga)\supseteq N_U(\be)$.

Conversely, suppose $\bd(\al)=\bd(\be)$, $\ker(\al)\supseteq\ker(\be)$ and $N_U(\al)\supseteq N_U(\be)$.  From the latter we have $u\leq s$, and we may assume that $G_i=C_i$ for $1\leq i\leq u$.  Together with $\bd(\al)=\bd(\be)$ and $\ker(\al)\supseteq\ker(\be)$, it follows that each $A_i$ ($1\leq i\leq r$) and each $C_j$ ($u<j\leq s$) is a union of~$E_k$'s.  Since $\be\be^*=\partI{E_1}{E_q}{C_1}{C_u}{E_1}{E_q}{C_1}{C_u}$, it quickly follows that $\al=\be\be^*\al\leqR\be$.

\pfitem{\ref{G4}}  This follows quickly from \ref{G1}, or alternatively from Lemma \ref{l:green_S}\ref{GS2}, as ${\al\al^*=\partI{A_1}{A_r}{C_1}{C_s}{A_1}{A_r}{C_1}{C_s}}$ and $\be\be^*=\partI{E_1}{E_q}{G_1}{G_u}{E_1}{E_q}{G_1}{G_u}$.

\pfitem{\ref{G2} and \ref{G5}}  These are dual to \ref{G1} and \ref{G4}.

\pfitem{\ref{G3}}  If $\al\leqJ\be$, then $\al=\ga\be\de$ for some $\ga\in\KK_{mk}$ and $\de\in\KK_{ln}$; \eqref{id2} then gives $\rank(\al)\leq\rank(\be)$.  If~$\KK$ is one of $\B$ or $\TL$, then also $\rank(\al)\equiv m\equiv k\equiv\rank(\be) \Mod2$, using the facts noted before the statement.

Conversely, suppose $r\leq q$, and also that $r\equiv q\Mod2$ if $\KK$ is either $\B$ or $\TL$.  Then $\al=\ga\be\de$ with
\[
\ga=
\Big( 
{ \scriptsize \renewcommand*{\arraystretch}{1}
\begin{array} {\c|\c|\c|\c|\c|\c|\c|\c|\cend}
A_1 \:&\: \cdots \:&\: A_r \:&\: C_1 \:& \multicolumn{4}{c|}{\cdots\cdots\cdots\cdots} &\: C_s \\ \cline{4-9}
E_1 \:&\: \cdots \:&\: E_r \:&\: E_{r+1} \:&\: \cdots \:&\: E_q \:&\: G_1 \:&\: \cdots \:&\: G_u
\rule[0mm]{0mm}{2.7mm}
\end{array} 
}
\hspace{-1.5 truemm} \Big)
\AND \de=
\Big( 
{ \scriptsize \renewcommand*{\arraystretch}{1}
\begin{array} {\c|\c|\c|\c|\c|\c|\c|\c|\cend}
F_1 \:&\: \cdots \:&\: F_r \:&\: F_{r+1} \:&\: \cdots \:&\: F_q \:&\: H_1 \:&\: \cdots \:&\: H_v \\ \cline{4-9}
B_1 \:&\: \cdots \:&\: B_r \:&\: D_1 \:& \multicolumn{4}{c|}{\cdots\cdots\cdots\cdots} &\: D_t
\rule[0mm]{0mm}{2.7mm}
\end{array} 
}
\hspace{-1.5 truemm} \Big).
\]
This clearly completes the proof (of \ref{G3}) if $\KK=\P$.  If $\al,\be$ belong to $\PB$, then $\ga,\de$ do as well, so the proof is also complete for $\KK=\PB$.  If $\al,\be$ are planar, then the above assumption about the ordering of the transversals, combined with Lemma \ref{l:PP}, ensures that $\ga,\de$ are planar, so the proof is complete for the categories $\PP$ and $\M$ as well.  
Finally, suppose $\KK$ is one of $\B$ or~$\TL$.  Note that here we have $|E_i|=|F_i|=1$ and $|G_j|=|H_k|=2$ for appropriate $i,j,k$; we write $E_i=\{e_i\}$ and $F_i=\{f_i\}$ for all $i$.  Since $r\equiv q\Mod2$, the sets $\{e_{r+1},\ldots,e_q\}$ and $\{f_{r+1},\ldots,f_q\}$ have even size, so we may replace $\ga,\de$ as above with
\[
\ga'=
\Big( 
{ \scriptsize \renewcommand*{\arraystretch}{1}
\begin{array} {\c|\c|\c|\c|\c|\c|\c|\c|\cend}
A_1 \:&\: \cdots \:&\: A_r \:&\: C_1 \:& \multicolumn{4}{c|}{\cdots\cdots\cdots\cdots\cdots\cdots} &\: C_s \\ \cline{4-9}
E_1 \:&\: \cdots \:&\: E_r \:&\: e_{r+1},e_{r+2} \:&\: \cdots \:&\: e_{q-1},e_q \:&\: G_1 \:&\: \cdots \:&\: G_u
\rule[0mm]{0mm}{2.7mm}
\end{array} 
}
\hspace{-1.5 truemm} \Big)
\anD
 \de'=
\Big( 
{ \scriptsize \renewcommand*{\arraystretch}{1}
\begin{array} {\c|\c|\c|\c|\c|\c|\c|\c|\cend}
F_1 \:&\: \cdots \:&\: F_r \:&\: f_{r+1},f_{r+2} \:&\: \cdots \:&\: f_{q-1},f_q \:&\: H_1 \:&\: \cdots \:&\: H_v \\ \cline{4-9}
B_1 \:&\: \cdots \:&\: B_r \:&\: D_1 \:& \multicolumn{4}{c|}{\cdots\cdots\cdots\cdots\cdots\cdots} &\: D_t
\rule[0mm]{0mm}{2.7mm}
\end{array} 
}
\hspace{-1.5 truemm} \Big).
\]
Then $\ga',\de'\in\B$, and we still have $\al=\ga'\be\de'$.  Moreover, if $\al,\be\in\TL$, then by Lemma \ref{l:PP}, the $e_i$ (respectively, $f_i$) are unnested by the upper (respectively, lower) nontransversals of $\be$, so it quickly follows that $\ga',\de'\in\TL$.

\pfitem{\ref{G6}}  The equivalence $\al\J\be\iff\rank(\al)=\rank(\be)$ clearly follows from \ref{G3}, and the equivalence $\al\J\be \iff \al\D\be$ from Proposition~\ref{p:stab}.
\epf

We immediately obtain the following description of the ${\J}={\D}$-classes in each diagram category, along with their ordering.  

\begin{cor}\label{c:D}
Let $\KK$ denote any of the categories $\P$, $\PB$, $\B$, $\PP$, $\M$ or $\TL$.  Let $m,n\in\N$, and suppose $m\equiv n\Mod2$ if $\KK$ is $\B$ or $\TL$.  Then the ${\J}={\D}$-classes of the hom-set $\Kmn$ are the sets
\begin{align*}
D_r = D_r(\Kmn) &= \set{\al\in\Kmn}{\rank(\al)=r}  &&\text{for each $0\leq r\leq\min(m,n)$,}\\
& &&\text{and where $r\equiv m\equiv n\Mod2$ if $\KK=\B$ or $\TL$.}
\end{align*}
These form a chain under the usual ordering of $\J$-classes: $D_q\leq D_r \iff q\leq r$.  \epfres
\end{cor}

The next result gives formulae for the number of $\R$-classes contained in the $\D$-classes of our diagram categories, as well as the sizes of $\H$-classes.  The statement uses the number sequences defined in Section \ref{ss:K}.

\newpage

\begin{prop}\label{p:D}
Let $\KK$ denote any of the categories $\P$, $\PB$, $\B$, $\PP$, $\M$ or $\TL$.  Let $m,n\in\N$, fix some $0\leq r\leq\min(m,n)$, and suppose $r\equiv m\equiv n\Mod2$ if $\KK$ is $\B$ or $\TL$.
\ben
\item \label{D1} The number of $\R$-classes contained in $D_r=D_r(\Kmn)$ is given by
\[
|D_r/{\R}| = \begin{cases}
\sum_{q=r}^m\binom qrS(m,q) &\text{if $\KK=\P$}\\[2truemm]
\binom mra(m-r) &\text{if $\KK=\PB$}\\[2truemm]
\binom mr(m-r-1)!! &\text{if $\KK=\B$}\\[2truemm]
\frac{2r+1}{2m+1}\binom{2m+1}{m-r} &\text{if $\KK=\PP$}\\[2truemm]
\mu(m,r) &\text{if $\KK=\M$}\\[2truemm]
\frac{r+1}{m+1}\binom{m+1}{(m-r)/2} &\text{if $\KK=\TL$.}
\end{cases}
\]
\item \label{D2} The size of any $\H$-class $H$ in $D_r$ is given by
\[
|H| = \begin{cases}
r! &\text{if $\KK$ is one of $\P$, $\PB$ or $\B$}\\
1 &\text{if $\KK$ is one of $\PP$, $\M$ or $\TL$.}
\end{cases}
\]
\een
\end{prop}

\pf
Part \ref{D2} is clear, and part \ref{D1} follows from arguments given elsewhere, specifically:
\bit
\item Theorems 7.5, 8.4 and 9.5 from \cite{EG2017} for $\P$, $\B$ and $\TL$, respectively, and
\item Propositions 2.7 and 2.8 from \cite{DEG2017} for $\PB$ and $\M$, respectively.
\eit
The result for $\PP$ follows from that for $\TL$, and the isomorphism $\PP\cong\TL^\even$, keeping in mind that $\rank(\alt)=2\rank(\al)$ for $\al\in\PP$; cf.~Figure \ref{f:TLiso}.
\epf

\begin{rem}\label{r:DrB}
By symmetry, we may obtain formulae for $|D_r/{\L}|$ in each of our diagram categories; we simply replace all $m$'s in the above formulae for $|D_r/{\R}|$ by $n$'s.  Then also $|D_r/{\H}|=|D_r/{\R}|\cdot|D_r/{\L}|$; multiplying this by the size of an arbitrary $\H$-class in $D_r$ gives a formula for $|D_r|$; summing over appropriate $r$ gives a formula for $|\Kmn|$, although simpler such formulae are given in Proposition \ref{p:sizes}.  

For example, in the Brauer category $\B$, we have
\begin{align*}
|D_r/{\R}| &= \tbinom mr(m-r-1)!! , & |D_r/{\H}| &= \tbinom mr\tbinom nr(m-r-1)!! (n-r-1)!!, \\
|D_r/{\L}| &= \tbinom nr(n-r-1)!! , & |D_r| &= \tbinom mr\tbinom nr(m-r-1)!! (n-r-1)!!r!.
\end{align*}
Summing over $r$, and using $|\B_{mn}| = (m+n-1)!!$, we obtain the identity
\[
(m+n-1)!! = \sum_{0\leq r\leq\min(m,n) \atop r\equiv n\Mod2} \tbinom mr\tbinom nr(m-r-1)!! (n-r-1)!!r!.
\]
\end{rem}

\section{Sandwich semigroups in diagram categories}\label{s:sandwichK}

Throughout this section $\KK$ denotes any of the categories $\P$, $\PB$, $\B$, $\PP$, $\M$ or $\TL$.  We fix some $\si\in\KK$, say with $\si\in\KK_{nm}$, with the aim of studying the sandwich semigroup $\Kmns=(\Kmn,\star_\si)$; note that if~$\KK$ is either $\B$ or $\TL$, then we have $m\equiv n\Mod2$.  Since $\Kmns\cong\KK_{nm}^{\si^*}$, we may assume by symmetry that $m\geq n$.  Note then that $r\leq n\leq m$.  We will assume that
\begin{equation}\label{e:si}
\si = \partI{X_1}{X_r}{U_1}{U_s}{Y_1}{Y_r}{V_1}{V_t}.
\end{equation}
In the case that $\KK$ is one of $\PP$, $\M$ or $\TL$, we will additionally assume that $\min(X_1)<\cdots<\min(X_r)$.  (When $\KK$ is one of $\P$, $\PB$ or $\B$, we will make no assumption on the ordering of any of the blocks of~$\si$.)
It will also be convenient to define the partition
\begin{equation}\label{e:tau}
\tau = \Big( 
{ \scriptsize \renewcommand*{\arraystretch}{1}
\begin{array} {\c|\c|\c|\c|\c|\cend}
 1 \:&\: \cdots \:&\: r \:& \multicolumn{3}{c}{} \\ \cline{4-6}
 X_1 \:&\: \cdots \:&\: X_r \:&\: U_1 \:&\: \cdots \:&\: U_s 
\rule[0mm]{0mm}{2.7mm}
\end{array} 
}
\hspace{-1.5 truemm} \Big)\in\KK_{rn}.
\end{equation}
The assumption on the ordering of the $X_i$ ensures that $\tau$ is planar if $\si$ is; cf.~Lemma \ref{l:PP}.

\subsection{Commutative diagrams and connections to nonsandwich diagram monoids}\label{ss:nonsandwichK}

We begin with some brief comments on the commutative diagrams in \eqref{e:CD2} when applied to diagram categories.
First note that $\si^*\in V(\si)$, and that
\[
\si\si^* = \partI{X_1}{X_r}{U_1}{U_s}{X_1}{X_r}{U_1}{U_s}\in\KK_n \AND \si^*\si = \partI{Y_1}{Y_r}{V_1}{V_t}{Y_1}{Y_r}{V_1}{V_t} \in \KK_m.
\]
Any partition $\al\in\si\si^*\KK_n\si\si^*$ contains the nontransversals $U_i$ and $U_i'$ for each $1\leq i\leq s$; any other block of $\al$ (whether a transversal or nontransversal) is of the form $B=\bigcup_{i\in I}X_i\cup\bigcup_{j\in J}X_j'$ for some subsets $I,J\sub[r]$, with at least one of $I,J$ nonempty.  We write $\al^\natural$ for the partition of $[r]\cup[r]'$ with each such block $B\sub[n]\cup[n]'$ replaced by $I\cup J'\sub[r]\cup[r]'$.
Note that in fact $\al^\natural = \tau\al\tau^*$ for any $\al\in\si\si^*\KK_n\si\si^*$, where $\tau$ is defined in \eqref{e:tau}; in particular, since $\tau\in\KK$, this shows that $\al^\natural\in\KK_r$ for any such $\al$.
Since $\tau^*\tau=\si\si^*$, it quickly follows that the map 
\begin{equation}\label{e:Kr}
\si\si^*\KK_n\si\si^*\to\KK_r:\al\mt\al^\natural
\end{equation}
is an isomorphism.  Combining \eqref{e:CD2} and \eqref{e:Kr}, we obtain the following commutative diagrams of semigroup surmorphisms:
\begin{equation}\label{e:CD3}
\begin{tikzcd} 
~ & \Kmns \arrow[swap]{dl}{\al\mt \al\si} \arrow{dr}{\al\mt \si\al} & \\
\KK_m\si^*\si \arrow[swap]{dr}{\be\mt \si\be\si^*} & & \si\si^*\KK_n \arrow{dl}{\be\mt \be\si\si^*}\\
& \si\si^*\KK_n\si\si^*
\arrow{dd}{\ga\mt\ga^\natural} & \\
\\
& \KK_r & 
\end{tikzcd}
\qquad\quad
\begin{tikzcd} [column sep=small]
~ & \Reg(\Kmns) \arrow[swap]{dl}{\al\mt \al\si} \arrow{dr}{\al\mt \si\al} & \\
\Reg(\KK_m\si^*\si) \arrow[swap]{dr}{\be\mt \si\be\si^*} & & \Reg(\si\si^*\KK_n) \arrow{dl}{\be\mt \be\si\si^*}\\
& \si\si^*\KK_n\si\si^*
\arrow{dd}{\ga\mt\ga^\natural} & \\
\\
& \KK_r & 
\end{tikzcd}
\end{equation}
We will use the symbol $\Psi$ to denote the surmorphism
\begin{equation}\label{e:Psi}
\Psi : \Reg(\Kmns)\to\KK_r:\al\mt(\si\al\si\si^*)^\natural
\end{equation}
induced from the second diagram in \eqref{e:CD3}; in fact, since $\si\si^*=\tau^*\tau$, we have $\al\Psi=\tau\si\al\si\si^*\tau^*=\tau\si\al\tau^*$ for all $\al\in\Reg(\Kmns)$.  It is useful to note here that
\begin{equation}\label{e:rank}
\rank(\al)=\rank(\al\Psi) \qquad\text{for all $\al\in P^\si$.}
\end{equation}
Indeed, from
\[
\al 
\J \si\al\si 
= \si\si^*\si \al \si\si^*\si 
= \tau^*\tau\si \al \tau^*\tau\si 
\leqJ \tau\si \al \tau^*
= \tau\si \al \tau^*\tau\tau^*
= \tau\si \al \si\si^*\tau^*
\leqJ \si\al\si
\J \al,
\]
it follows that $\al\J\tau\si\al\tau^*=\al\Psi$, and so $\rank(\al)=\rank(\al\Psi)$ by Theorem \ref{t:G}\ref{G6}.

By Theorem \ref{t:RG}\ref{RG1}, the restriction
\[
\Psi|_{H_\al^\si}:H_{\al}^{\si}\to H_{\al\Psi}
\]
of $\Psi$ to any $\H^\si$-class $H_\al^\si$ of $\Reg(\Kmns)$ is a bijection, and such a restriction maps (non)group $\H^\si$-classes of $\Reg(\Kmns)$ to (non)group $\H$-classes of $\KK_r$, respectively; in the case that $H_{\al}^{\si}$ is a group, $\Psi|_{H_\al^\si}$ is an isomorphism.
In particular, if $\al\in\Kmn$ is a $\star_\si$-idempotent of rank $q$, then the group $H_\al^\si$ is trivial if $\KK$ is one of $\PP$, $\M$ or $\TL$; otherwise, $H_\al^\si$ is isomorphic to the symmetric group $\S_q$.

\subsection{Green's relations and regularity in $\Kmns$}\label{ss:GreenKmns}

As discussed in Section \ref{ss:GreenSija}, Green's relations and regularity in the sandwich semigroups $\Kmns$ are governed by the sets $P_1^\si$, $P_2^\si$, $P_3^\si$ and $P^\si=P_1^\si\cap P_2^\si$.  Because of stability (cf.~Proposition \ref{p:stab}), it follows from \cite[Proposition~2.4]{Sandwich1} that $P_3^\si=P^\si$.  Since $\KK$ is regular, Proposition \ref{p:P} gives $P^\si=\Reg(\Kmns)$.  Combining all of this with Theorem~\ref{t:G}\ref{G6}, and again using stability, we obtain the following:

\begin{prop}\label{p:PK}
We have
\begin{align*}
P_1^\si &= \set{\al\in\Kmn}{\rank(\al\si)=\rank(\al)} , \\
P_2^\si &= \set{\al\in\Kmn}{\rank(\si\al)=\rank(\al)} , \\
\Reg(\Kmns) = P^\si = P_3^\si &= \set{\al\in\Kmn}{\rank(\al\si)=\rank(\si\al)=\rank(\al)} \\
\epfreseq &= \set{\al\in\Kmn}{\rank(\si\al\si)=\rank(\al)} .
\end{align*}
\end{prop}

It is possible to give combinatorial criteria for membership of these sets, but it does not give a great deal of additional insight (e.g., it does not allow us to readily compute the size of $P^\si=\Reg(\Kmns)$ in general).  For example, if $\rank(\al)=q$, then the product $\al\si$ will have rank $q$ if and only if $q$ of the $\coker(\al)\vee\ker(\si)$-classes contain elements of $\codom(\al)$, and each of these classes also contain elements of $\dom(\si)$.  In the case that $\KK=\B$ is the Brauer category, however, cleaner descriptions can be made (cf.~Proposition \ref{p:PB}), and these do lead to useful consequences (see for example Theorems~\ref{t:RegB} and~\ref{t:eBmns} and Corollary \ref{c:PB}).

Together with the above description of the sets $P_1^\si$, $P_2^\si$ and $P^\si$, characterisations of Green's relations and preorders in $\Kmns$ follow from Theorem \ref{t:Green_Sij} and Proposition \ref{p:<}.  Note that ${\J^\si}={\D^\si}$ since~$\Kmns$ is finite.  It will be convenient to describe the regular ${\J^\si}={\D^\si}$-classes.

\begin{prop}\label{p:DP}
\ben
\item \label{DP1} The regular ${\J^\si}={\D^\si}$-classes of~$\Kmns$ are precisely the sets
\begin{align*}
D_q^\si = D_q^\si(\Kmns) =D_q\cap P^\si &= \set{\al\in P^\si}{\rank(\al)=q}   &&\text{for each $0\leq q\leq r$,}\\
& &&\text{and where $q\equiv r\Mod2$ if $\KK=\B$ or $\TL$.}
\end{align*}
These form a chain under the usual ordering of $\J^\si$-classes: $D_p^\si\leq D_q^\si \iff p\leq q$.  
\item \label{DP2} The group $\H^\si$-classes in $D_q^\si$ are
\bit
\item isomorphic to the symmetric group $\S_q$ if $\KK$ is one of $\P$, $\PB$ or $\B$,
\item trivial if $\KK$ is one of $\PP$, $\M$ or $\TL$.
\eit
\een
\end{prop}

\pf
\firstpfitem{\ref{DP1}} This follows immediately from Theorem \ref{t:Green_Sij}\ref{GSij4}, Proposition \ref{p:<}\ref{<6} and Theorem \ref{t:G}\ref{G3} and \ref{G6}.  

\pfitem{\ref{DP2}}  We observed this at the end of Section \ref{ss:nonsandwichK}.
\epf

In general, $\Kmns$ has many more ${\J^\si}={\D^\si}$-classes than just these regular ones; cf.~Figures~\ref{f:emax7} and~\ref{f:emax8}.

Next we describe the maximal $\J^\si$-classes in a sandwich semigroup $\Kmns$.  Note that if $r=n$, then
$
\si = \Big( 
{ \scriptsize \renewcommand*{\arraystretch}{1}
\begin{array} {\c|\c|\c|\c|\c|\cend}
 1 \:&\: \cdots \:&\: n \:& \multicolumn{3}{c}{} \\ \cline{4-6}
 Y_1 \:&\: \cdots \:&\: Y_n \:&\: V_1 \:&\: \cdots \:&\: V_t 
\rule[0mm]{0mm}{2.7mm}
\end{array} 
}
\hspace{-1.5 truemm} \Big)
$
is right-invertible; indeed, in this case $\si\si^*=\id_n$.

\begin{prop}\label{p:maxPPBB}
Suppose $\KK$ is one of $\P$, $\PB$ or $\B$, and that $m\geq n$.
\ben
\item \label{maxPPBB1} If $r<n$, then the maximal $\J^\si$-classes of $\Kmns$ are the singleton sets $\{\al\}$ for $\al\in\Kmn$ with $\rank(\al)>r$.
\item \label{maxPPBB2} If $r=n$, then the set $D_r^\si=\set{\al\in P^\si}{\rank(\al)=r}$ is the maximum $\J^\si$-class of $\Kmns$; it is a left-group over $\S_r$.
\een
\end{prop}

\pf
We begin with \ref{maxPPBB2}.  As noted above, if $r=n$ then $\si$ is right-invertible.  It then follows from Proposition \ref{p:max3}\ref{max31} that $\Kmns$ has a maximum $\J^\si$-class, and that this is $J_{\si^*}^\si$.  Since $\si^*\in P^\si$ and $\rank(\si^*)=r$, it follows from Proposition \ref{p:DP}\ref{DP1} that $J_{\si^*}^\si=D_r^\si$.  It follows from Proposition \ref{p:max3}\ref{max32} that~$D_r^\si$ is a left-group over $H_{\si^*}^\si$.  Since $H_{\si^*}^\si\cong\S_r$ (cf.~Proposition \ref{p:DP}\ref{DP2}), the proof is complete.

\pfitem{\ref{maxPPBB1}}  
For any $\al$ with $\rank(\al)>r$, Theorem \ref{t:G}\ref{G3} gives $\al\not\leqJ\si$; it follows from Lemma \ref{l:max1} that $J_\al^\si=\{\al\}$ is a (trivial) maximal $\J^\si$-class.  It remains to show that there is no nontrivial maximal $\J^\si$-class.  By Proposition \ref{p:Jba}, it suffices to show that there exists $\al\in\Kmn$ such that in $\KK$, $(\si,\si\al\si)\in{\J}$ but $(\al,\si)\not\in{\J}$: i.e., that $\rank(\si\al\si)=r<\rank(\al)$.  To construct such an $\al$ we must consider a number of cases.

\pfcase{1}  First suppose $s$ and $t$ are both nonzero (and note that this must be the case if $\KK$ is one of~$\PB$ or~$\B$, as $r<n\leq m$).  Put
\[
\al = 
\Big( 
{ \scriptsize \renewcommand*{\arraystretch}{1}
\begin{array} {\c|\c|\c|\c|\c|\c|\cend}
 Y_1 \:&\: \cdots \:&\: Y_r \:&\: V_t \:&\: V_1 \:&\: \cdots \:&\: V_{t-1} \\ \cline{5-7}
 X_1 \:&\: \cdots \:&\: X_r \:&\: U_s \:&\: U_1 \:&\: \cdots \:&\: U_{s-1} 
\rule[0mm]{0mm}{2.7mm}
\end{array} 
}
\hspace{-1.5 truemm} \Big).
\]
Then $\si\al\si=\si$, so certainly $\rank(\si\al\si)=r$.  We also have $\rank(\al)=r+1>r$.  This completes the proof (in this case) if $\KK=\P$.  If $\KK$ is one of $\PB$ or $\B$, then $\al$ as above may not belong  to $\KK$, so we may have to modify the definition of $\al$ to ensure that it does.
\bit
\item If $|U_s|=|V_t|=2$ (the only possibility for $\KK=\B$), then we write $U_s=\{u_1,u_2\}$ and $V_t=\{v_1,v_2\}$, and we replace the transversal $V_t\cup U_s'$ of $\al$ by the pair of transversals $\{v_1,u_1'\}$ and $\{v_2,u_2'\}$.  Again we have $\si=\si\al\si$, but this time $\rank(\al)=r+2$.  
\item If $|U_s|=2$ and $|V_t|=1$, say with $U_s=\{u_1,u_2\}$ and $V_t=\{v\}$, then we replace the transversal $V_t\cup U_s'$ of $\al$ by the pair of blocks $\{v,u_1'\}$ and $\{u_2'\}$.  Here $\si=\si\al\si$ and $\rank(\al)=r+1$.
\item The case in which $|U_s|=1$ and $|V_t|=2$ is dual.
\item The case in which $|U_s|=|V_t|=1$ poses no problems, as then $\al$ (as above) belongs to $D_{r+1}(\PB)$ provided that $\si$ does.
\eit

\pfcase{2}  Next suppose $s=0$ but $t>0$, so that $\si=\Big( 
{ \scriptsize \renewcommand*{\arraystretch}{1}
\begin{array} {\c|\c|\c|\c|\c|\cend}
 X_1 \:&\: \cdots \:&\: X_r \:& \multicolumn{3}{c}{} \\ \cline{4-6}
 Y_1 \:&\: \cdots \:&\: Y_r \:&\: V_1 \:&\: \cdots \:&\: V_t 
\rule[0mm]{0mm}{2.7mm}
\end{array} 
}
\hspace{-1.5 truemm} \Big)$.
As noted above, this case (and all remaining ones) only applies to~$\KK=\P$.  Since $r<n$, we may assume without loss of generality that $|X_r|\geq2$.  Fix some $x\in X_r$, write $Z=X_r\sm\{x\}$, and define
\[
\al = \Big( 
{ \scriptsize \renewcommand*{\arraystretch}{1}
\begin{array} {\c|\c|\c|\c|\c|\c|\c|\cend}
 Y_1 \:&\: \cdots \:&\: Y_{r-1} \:&\: Y_r \:&\: V_t \:&\: V_1 \:&\: \cdots \:&\: V_{t-1} \\ \cline{6-8}
 X_1 \:&\: \cdots \:&\: X_{r-1} \:&\: x \:&\: Z \:& \multicolumn{3}{c}{}
\rule[0mm]{0mm}{2.7mm}
\end{array} 
}
\hspace{-1.5 truemm} \Big).
\]
Then we have $\si=\si\al\si$ and $\rank(\al)=r+1$.

\pfcase{3}  The case in which $t=0$ and $s>0$ is dual.

\pfcase{4}  Finally suppose $s=t=0$, so that $\si=\partpermIII{X_1}\cdots{X_r}{Y_1}\cdots{Y_r}$.  Since $r<n\leq m$ we have $|X_i|\geq2$ and $|Y_j|\geq2$ for some $i,j\in[r]$.  Renaming if necessary, we may assume that $i=r$, and that $j=r$ or $r-1$.
Fix some $x\in X_r$ and $y\in Y_j$, and write $U=X_r\sm\{x\}$ and $V=Y_j\sm\{y\}$.
\bit
\item If $j=r$, then with $\al = \partpermV{Y_1}{Y_{r-1}}yV{X_1}{X_{r-1}}xU$ we have $\si=\si\al\si$ and $\rank(\al)=r+1$.
\item If $j=r-1$, then with $\al = 
\Big(
{ \scriptsize \renewcommand*{\arraystretch}{1}
\begin{array} {\c|\c|\c|\c|\c|\c}
 Y_1 \:&\: \cdots \:&\: Y_{r-2} \:&\: Y_r \:&\: y \:&\: V \\ 
 X_1 \:&\: \cdots \:&\: X_{r-2} \:&\: X_{r-1} \:&\: x \:&\: U
\rule[0mm]{0mm}{2.7mm}
\end{array} 
}
\Big)$
we have $\si\al\si=\partpermV{X_1}{X_{r-2}}{X_{r-1}}{X_r}{Y_1}{Y_{r-2}}{Y_r}{Y_{r-1}}$; in particular, $\rank(\si\al\si)=r<r+1=\rank(\al)$.  \qedhere
\eit
\epf

\begin{rem}\label{r:max2}
In the proof of Proposition \ref{p:maxPPBB}\ref{maxPPBB1}, every case except for the $j=r-1$ subcase of Case~4 involved showing that~$\si$ has a pre-inverse not $\J$-related to $\si$ (in $\KK$); cf.~Corollary \ref{c:max}.  In the exceptional (sub)case, the~$\al$ constructed did not satisfy $\si=\si\al\si$, but one may check that $\si=\si\al\si\al\si$; in general, there may or may not be any $\al$ satisfying $\si=\si\al\si$ and $\rank(\al)>r$; this depends on the sizes of the other blocks (if any) of $\si$; cf.~Example \ref{e:max}\ref{emax8}.
\end{rem}

We now give the characterisation for the planar diagram categories.  For the statement, recall that $\Pre(\si)=\set{\al\in\Kmn}{\si=\si\al\si}$, and that $D_r(\Kmn)=\set{\al\in\Kmn}{\rank(\al)=r}$.

\begin{prop}\label{p:maxPPMTL}
Suppose $\KK$ is one of $\PP$, $\M$ or $\TL$, and that $m\geq n$.  
\ben
\item \label{maxPPMTL1} If $r<n$, then the trivial maximal $\J^\si$-classes of $\Kmns$ are the singleton sets $\{\al\}$ for $\al\in\Kmn$ with ${\rank(\al)>r}$.  Moreover, the following are equivalent:
\bit
\item $\Kmns$ has a nontrivial maximal $\J^\si$-class,
\item $\Pre(\si)\sub D_r(\Kmn)$,
\item $\Pre(\si)=V(\si)$,
\eit
in which case the nontrivial maximal $\J^\si$-class is the set ${D_r^\si=\set{\al\in P^\si}{\rank(\al)=r}}$.
\item \label{maxPPMTL2} If $r=n$, then the set $D_r^\si=\set{\al\in P^\si}{\rank(\al)=r}$ is the maximum $\J^\si$-class of $\Kmns$; it is a left-zero semigroup.  
\een
\end{prop}

\pf
\firstpfitem{\ref{maxPPMTL1}}  The statement concerning trivial maximal $\J^\si$-classes again follows from Lemma \ref{l:max1} and Theorem \ref{t:G}\ref{G3}.  The assertion regarding the existence of a nontrivial maximal $\J^\si$-class (which would have to be~$J_{\si^*}^\si=D_r^\si$, by Lemma \ref{l:max2}\ref{max22}) follows from Proposition \ref{p:max2}, since $\KK$ is regular, stable and $\H^\si$-trivial.

\pfitem{\ref{maxPPMTL2}}  This is proved in the same way as Proposition \ref{p:maxPPBB}\ref{maxPPBB2}, but keeping in mind that group $\H$-classes of $\KK_r$ are all trivial if $\KK$ is any of $\PP$, $\M$ or $\TL$ (cf.~Proposition \ref{p:DP}\ref{DP2}).
\epf

\begin{rem}
Sometimes the containment $\Pre(\si)\sub D_r(\Kmn)$ holds, and sometimes it does not; it depends on the structure of $\si$.  Example \ref{e:max}\ref{emax7} gives examples in both cases when $\KK=\TL$; it is easy to construct examples for the categories $\PP$ and $\M$ (by hand or with GAP \cite{GAP}).
\end{rem}

It will also be convenient to describe the minimum $\J^\si$-class of a sandwich semigroup $\Kmns$.  Note that every finite semigroup has a minimum $\J$-class, and this coincides with the minimal ideal of the semigroup.

\begin{prop}\label{p:min}
If $z$ is the smallest possible rank of partitions from $\Kmn$, then the minimal ideal of $\Kmns$ is $D_z=D_z^\si$.  Further, we have $D_z/{\R^\si}=D_z/{\R}$ and $D_z/{\L^\si}=D_z/{\L}$.
\end{prop}

\pf
For any $\al\in D_z$, we have $\rank(\si\al\si)\leq\rank(\al)=z$, so that $\rank(\si\al\si)=z$ (as $z$ is the smallest possible rank).  It follows that $D_z\sub P^\si$ (cf.~Proposition \ref{p:PK}), and so $D_z^\si=D_z\cap P^\si=D_z$.  By Proposition~\ref{p:<}\ref{<3} and Theorem \ref{t:G}\ref{G3}, $D_z$ is the minimum $\J^\si$-class.  

For the last assertion, suppose $\al\in D_z$.  Since $\al\in P^\si$, Theorem \ref{t:Green_Sij}\ref{GS1} gives $R_\al^\si=R_\al\cap P_2^\si$.  Since $R_\al\sub D_z\sub P^\si\sub P_2^\si$, it follows that $R_\al^\si=R_\al\cap P_2^\si=R_\al$.  The statement for $\L^\si$-classes is dual.
\epf

\subsection{The idempotent-generated subsemigroup}\label{ss:IGKmns}

We write $E(\Kmns) = \set{\al\in\Kmn}{\al=\al\star_\si\al}$ for the set of all $\star_\si$-idempotents from $\Kmns$, and ${\E(\Kmns) = \la E(\Kmns)\ra}$ for the idempotent-generated subsemigroup of $\Kmns$.  Using \eqref{e:ESija} and the isomorphism~\eqref{e:Kr}, it follows that
\begin{equation}\label{e:EK}
E(\Kmns) = E(\KK_r) \Psi^{-1} \AND \E(\Kmns) = \E(\KK_r) \Psi^{-1},
\end{equation}
where $\Psi:P^\si=\Reg(\Kmns)\to\KK_r$ is the surmorphism from \eqref{e:Psi}.  
Since the idempotent-generated subsemigroups $\E(\KK_r)$ are well understood for each $\KK$ \cite{East2011,MM2007,DEG2017,HR2005,BDP2002},~\eqref{e:EK} is in principal enough to completely describe $\E(\Kmns)$.  However, in most cases we can give a more precise statement.  
Note that ${\E(\Kmns)=\E(\KK_r)\Psi^{-1}\sub P^\si}$.  The next result shows that we have equality when $\KK$ is $\TL$ or $\PP$.

\begin{thm}\label{t:ETLPP}
If $\KK$ is one of $\TL$ or $\PP$, then $\E(\Kmns)=P^\si=\Reg(\Kmns)$.
\end{thm}

\pf
Since $\KK_r$ is idempotent-generated \cite{BDP2002,HR2005}, we have $\E(\KK_r)=\KK_r$.  The result then follows from~\eqref{e:EK}, and the fact that $\Psi:\Reg(\Kmns)\to\KK_r$ is surjective.
\epf

For the other categories, we first prove a simple lemma:

\begin{lemma}\label{l:Dq}
Suppose $0\leq q\leq r$, and that $q\equiv r\Mod2$ if $\KK$ is $\B$ or $\TL$.  Then
\[
E(D_q(\KK_r))\Psi^{-1} = E(D_q^\si).
\]
\end{lemma}

\pf
This follows from $E(\Kmns) = E(\KK_r) \Psi^{-1}$, and the fact that $\rank(\al\Psi)=\rank(\al)$ for $\al\in P^\si$; cf.~\eqref{e:EK} and \eqref{e:rank}.
\epf

Next we consider the categories $\P$ and $\B$.  If $\KK$ denotes either of these, then by \cite[Proposition 16]{East2011} and \cite[Proposition 2]{MM2007} we have $\E(\KK_r) = \{\id_r\} \cup (\KK_r\sm\S_r)$.  (The papers \cite{East2011} and \cite{MM2007} were concerned with presentations for the singular ideal $\KK_r\sm\S_r$.)  It follows that
\[
\E(\Kmns) = \E(\KK_r)\Psi^{-1} = \id_r\Psi^{-1} \cup \bigcup_{q<r} D_q^\si.
\]
Note that the union is over $\set{q}{0\leq q<r}$ for $\KK=\P$ or over $\set{q}{0\leq q<r,\ q\equiv r\Mod2}$ for $\KK=\B$.  In both cases it follows that
\[
\E(\Kmns) = \id_r\Psi^{-1} \cup (P^\si\sm D_r^\si).
\]

\begin{thm}\label{t:EPB}
If $\KK$ is one of $\P$ or $\B$, then $\E(\Kmns)=V(\si)\cup(P^\si\sm D_r^\si)$.
\end{thm}

\pf
It remains to observe, using Lemma \ref{l:Dq} and Proposition \ref{p:Jba2}\ref{Jba22}, that
\[
\id_r\Psi^{-1} = E(D_r(\KK_r))\Psi^{-1} = E(D_r^\si) = E(J_{\si^*}^\si) = V(\si).  \qedhere
\]
\epf

It was shown in \cite[Theorem 3.18]{DEG2017} that the idempotent-generated subsemigroup of the partial Brauer monoid $\PB_r$ is given by
\[
\E(\PB_r) = E\big(D_r(\PB_r)\cup D_{r-1}(\PB_r)\big) \cup \bigcup_{q=0}^{r-2} D_q(\PB_r).
\]
As above, and using Lemma \ref{l:Dq}, we quickly obtain the following:

\begin{thm}\label{t:EPB2}
We have $\displaystyle\E(\PBmns)=E(D_r^\si\cup D_{r-1}^\si) \cup \bigcup_{q=0}^{r-2} D_q^\si$.  \epfres
\end{thm}

The situation for the Motzkin category $\M$ is more complicated, due to the far more complex structure of the idempotent-generated subsemigroup of $\M_r$; see \cite[Theorem~4.17]{DEG2017}.  Thus, we leave our description of $\E(\Mmns)$ at \eqref{e:EK}.

\subsection{Ideals of $\Reg(\Kmns)$}\label{ss:idealsK}

In this section we study the ideals of the regular subsemigroup $P^\si=\Reg(\Kmns)$.  The ideals of the monoids $\KK_r$ are often (but not always) idempotent-generated, and it turns out that there is a close parallel with idempotent-generation in ideals of $P^\si$.

Throughout Section \ref{ss:idealsK}, it will be convenient to define the indexing set
\[
Q = \set{\rank(\al)}{\al\in P^\si} = \begin{cases}
\set{q}{0\leq q\leq r} &\text{if $\KK$ is one of $\P$, $\PB$, $\PP$ or $\M$}\\
\set{q}{0\leq q\leq r,\ q\equiv r\Mod2} &\text{if $\KK$ is one of $\B$ or $\TL$.}
\end{cases}
\]
By Proposition \ref{p:DP}, the ${\J^\si}={\D^\si}$-classes $D_q^\si$ ($q\in Q$) of $P^\si$ form a chain: $D_p^\si\leq D_q^\si \iff p\leq q$.  It follows quickly from this that the ideals of $P^\si$ are precisely the sets
\begin{align*}
I_q = I_q(P^\si) &= \bigcup_{p\leq q} D_p^\si &&\hspace{-3cm}\text{for each $q\in Q$.}
\intertext{This mirrors the situation with the ideals of the monoid $\KK_r$, which are all of the form}
I_q(\KK_r) &= \bigcup_{p\leq q} D_p(\KK_r) &&\hspace{-3cm}\text{for each $q\in Q$.}
\end{align*}
Keeping in mind that $\PP_n\cong\TL_{2n}$ (cf.~\cite[Section 1]{HR2005}), it was shown in \cite{DEG2017,EG2017} that
\begin{equation}\label{e:mu}
\text{$I_q(\KK_r)$ is idempotent-generated $\iff$ $q\leq\mu$, where } \mu = \mu(\KK) = \begin{cases}
r &\text{if $\KK=\PP$ or $\TL$}\\
r-1 &\text{if $\KK=\P$}\\
r-2 &\text{if $\KK=\PB$ or $\B$}\\
\lfloor\frac r2\rfloor-1 &\text{if $\KK=\M$.}
\end{cases}
\end{equation}
In fact, it was shown in \cite[Theorems 7.5, 8.4 and 9.5]{EG2017} that if $\KK$ is any of $\P$, $\B$, $\PP$ or $\TL$, then a proper ideal $I_q(\KK_r)$, $q<r$, is generated by the idempotents in its top $\D$-class: i.e., $I_q(\KK_r) = \big\la E(D_q(\KK_r)\big\ra$.  This is true of approximately half the ideals of $\PB_r$, but only of the minimal ideal of $\M_r$ \cite{DEG2017}.  

In the following proof, it will be convenient to write $\ol\al=\al\Psi$ for $\al\in P^\si$.  If $\K$ is any of Green's relations, we will also write $\Kh_\al^\si$ for the $\gKh^\si$-class of $\al\in P^\si$.  (As in Section \ref{ss:nonsandwich}, and using the isomorphism $(\si\Kmn\si,\star_{\si^*})\cong\KK_r$, the $\gKh^\si$ relation is defined on $P^\si$ by $\al\gKh^\si\be\iff\ol\al\K\ol\be$ in~$\KK_r$.)  Throughout the proof, we make frequent use of \eqref{e:rank}, which says that $\rank(\al)=\rank(\ol\al)$ for all~$\al\in P^\si$.

\begin{thm}\label{t:IqK}
Let $\KK$ be one of $\P$, $\PB$, $\B$, $\PP$, $\M$ or $\TL$, and let $\mu=\mu(\KK)$ be as in \eqref{e:mu}.  Then for any $q\in Q$, the ideal $I_q=I_q(P^\si)$ is idempotent-generated if and only if $q\leq\mu$.
\end{thm}

\pf
Since $\Psi:P^\si\to\KK_r$ preserves ranks and is surjective, it follows that $I_q(\KK_r)=I_q(P^\si)\Psi$ for all $q\in Q$.  Thus, if $I_q=I_q(P^\si)$ is idempotent-generated, then so too is $I_q(\KK_r)$; it then follows from \eqref{e:mu} that $q\leq\mu$.

Conversely, suppose $q\leq\mu$, and let $\al\in I_q$.  To complete the proof, we must show that $\al$ is a product of idempotents from $I_q$.  

Write $\rank(\al)=p$, and note that $\ol\al\in D_p(\KK_r)\sub I_p(\KK_r)$.  Since $p\leq q\leq\mu$, it follows from~\eqref{e:mu} that $I_p(\KK_r)$ is idempotent-generated.  Thus, $\ol\al = \be_1\cdots\be_k$ for some idempotents $\be_1,\ldots,\be_k\in E(I_p(\KK_r))$.  Since ${p\geq\rank(\be_i)\geq\rank(\be_1\cdots\be_k)=p}$ for all~$1\leq i\leq k$, it follows in fact that $\be_1,\ldots,\be_k\in E(D_p(\KK_r))$.
For each $i$, fix some $\ga_i\in P^\si$ such that $\ol\ga_i=\be_i$, noting that $\ga_i\in E(D_p^\si)\sub E(I_q)$ by Lemma \ref{l:Dq}.

Since $\rank(\ol\al)=p=\rank(\be_1)$, we have $\be_1\J \ol\al = \be_1(\be_2\cdots\be_k)$, so it follows from stability that $\ol\ga_1=\be_1\R\be_1(\be_2\cdots\be_k)=\ol\al$: i.e., $\ga_1\gRh^\si\al$, which gives $\Rh_\al^\si=\Rh_{\ga_1}^\si$.  Since $\rank(\al)=\rank(\ga_1)$, we have $\al\D^\si\ga_1$, so we may fix some $\de\in R_\al^\si\cap L_{\ga_1}^\si$.  Note then that
\[
\de\in R_\al^\si\cap L_{\ga_1}^\si \sub \Rh_\al^\si\cap \Lh_{\ga_1}^\si = \Rh_{\ga_1}^\si\cap \Lh_{\ga_1}^\si = \Hh_{\ga_1}^\si.
\]
Since $\ol\ga_1$ is an idempotent of $\KK_r$, it follows from Theorem \ref{t:RG}\ref{RG3} that $\Hh_{\ga_1}^\si$ is a rectangular group, and hence $H_{\de}^\si$ is a group.  Thus, without loss of generality, we may assume that $\de$ is the identity element of this group.  Since $\de\in\Hh_{\ga_1}^\si$, we have $\ol\de\H\ol\ga_1$; but $\ol\de$ and $\ol\ga_1$ are idempotents of $\KK_r$, so it follows that~$\ol\de=\ol\ga_1$.  

The key conclusion of the previous paragraph is that there exists an idempotent $\de$ such that $\de\R^\si\al$ and $\ol\de=\ol\ga_1$.  By a symmetrical argument, there exists an idempotent $\ve$ such that $\ve\L^\si\al$ and $\ol\ve=\ol\ga_k$.  

Now consider the element $\rho = \de\star_\si(\ga_1\star_\si\cdots\star_\si\ga_k)\star_\si\ve$.  Then since $\ol\de=\ol\ga_1$ and $\ol\ve=\ol\ga_k$ are idempotents,
\[
\ol\rho = \ol\de (\ol\ga_1\cdots\ol\ga_k)\ol\ve = \ol\ga_1\cdots\ol\ga_k = \ol\al.
\]
Consequently, $\rank(\rho)=\rank(\ol\rho)=\rank(\ol\al)=p=\rank(\de)$, so that $\rho\J^\si\de$.  But then 
\[
\de \J^\si \rho =  \de\star_\si(\ga_1\star_\si\cdots\star_\si\ga_k)\star_\si\ve,
\]
so that $\de \R^\si  \de\star_\si(\ga_1\star_\si\cdots\star_\si\ga_k)\star_\si\ve$, by stability: i.e., $\de\R^\si\rho$.  Since $\de\R^\si\al$, it follows that $\rho\R^\si\al$.  A symmetrical argument shows that $\rho\L^\si\ve\L^\si\al$, so it follows that $\rho\H^\si\al$.  Thus, since $\ol\rho=\ol\al$ (shown above), and since $\Psi$ is injective when restricted to any $\H^\si$-class (cf.~Theorem \ref{t:RG}\ref{RG1}), it follows that $\al=\rho$.  But $\rho$ is a product of idempotents, so the proof is complete.
\epf

\begin{rem}\label{r:IqK}
It follows from the above proof that if $q\leq\mu$, then any element of $D_q^\si$ is a product of idempotents from $E(D_q^\si)$.  We will see in Section \ref{ss:ideals} that a stronger statement holds when $\KK=\B$, namely that any element of $I_q^\si$ is a product of idempotents from $E(D_q^\si)$.  This mirrors the situation with the monoids $\P_r$, $\B_r$, $\PP_r$ and $\TL_r$, as noted just after \eqref{e:mu}.

However, this is not the case for the other diagram categories.  For example, consider the partition $\si = \custpartn{1,2,3,4}{1,2,3,4}{\stline11\stline22\stline33\stline44\uarc34\darc34}$ from $\P_4=\P_{4,4}$.  The regular subsemigroup $P^\si=\Reg(\P_4^\si)$ has four ${\J^\si}={\D^\si}$-classes: $D_0^\si<D_1^\si<D_2^\si<D_3^\si$.  The proper ideal $I_2=D_0^\si\cup D_1^\si\cup D_2^\si$ is idempotent-generated (by Theorem~\ref{t:IqK}), but it is not generated by $E(D_2^\si)$.  Indeed, GAP \cite{GAP} shows that $I_2$ has size $2476$, while the subsemigroup generated by $E(D_2^\si)$ has size $2332$.
\end{rem}

\section{The Brauer category}\label{s:B}

It turns out that the Brauer category $\B$ has many properties not shared by any of the other diagram categories studied in this paper, many of which make $\B$ more amenable to further analysis.  For example:
\bit
\item Sandwich elements $\si,\tau\in\B_{nm}$ give rise to isomorphic variants $\Bmns\cong\Bmn^\tau$ if $\rank(\si)=\rank(\tau)$; this leads to a neat classification of the isomorphisms between sandwich semigroups in $\B$; see Section \ref{ss:IsoB}.
\item Every sandwich semigroup $\Bmns$ with $\rank(\si)<\min(m,n)$ is generated by its trivial maximal $\J^\si$-classes, a fact that leads to neat formulae for the rank of $\Bmns$; see Section \ref{ss:RankB}.
\item The sets $P_1^\si$ and $P_2^\si$ may be readily described in a way that yields convenient description of the regular subsemigroup of a sandwich semigroup; see Section \ref{ss:RegB}.  Among other things, this means that Green's classes and idempotents in the regular subsemigroup $P^\si=\Reg(\Bmns)$ may be enumerated by analysing joins of certain equivalence relations, as described in Section \ref{ss:ka}.
\item The regular subsemigroup $P^\si=\Reg(\Bmns)$ has the so-called \emph{MI-domination} property, as defined in \cite[Section 4]{Sandwich1}.  We show this in Section \ref{ss:MIB}, and then apply this to calculate the (idempotent) ranks of $P^\si$ and the idempotent-generated subsemigroup $\E(\Bmns)$.
\item The proper ideals of $P^\si=\Reg(\Bmns)$ are all generated by the idempotents in their top $\D^\si$-class, and simple formulae exist for their rank and idempotent rank (which are equal); see Section \ref{ss:ideals}.
\eit
We believe it would be interesting (but challenging) to investigate the corresponding problems for the other diagram categories.

\subsection{Isomorphism of sandwich semigroups in $\B$}\label{ss:IsoB}

In this section, we give necessary and sufficient conditions for two sandwich semigroups $\Bmns$ and $\Bklt$ to be isomorphic; see Theorem \ref{t:iso}.  We begin with a simple lemma; as well as being a key ingredient in the proof of Theorem \ref{t:iso}, it will also be used to simplify some calculations in Sections \ref{ss:RankB}--\ref{ss:ideals}.

\begin{lemma}\label{l:iso1}
If $m,n\in\N$, and if $\si,\tau\in\B_{nm}$ are such that $\rank(\si)=\rank(\tau)$, then $\Bmns\cong\B_{mn}^\tau$.
\end{lemma}

\pf
Since $\rank(\si)=\rank(\tau)$, we have $\si=\pi_1\tau\pi_2$ for some permutations $\pi_1\in\S_n$ and $\pi_2\in\S_m$.  It is then easy to check that the map $\al\mt\pi_2\al\pi_1$ determines an isomorphism $\Bmns\to\B_{mn}^\tau$.
\epf

\begin{figure}[ht]
\begin{center}
\begin{tikzpicture}
\foreach \x in {1,...,9} {\node[above] () at (1.9*\x,0) {\includegraphics[scale=0.3]{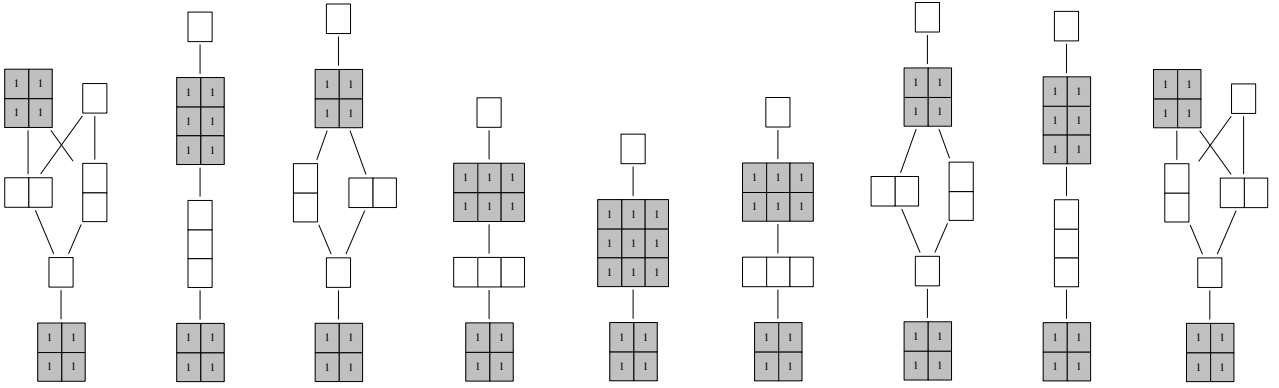}};}
\end{tikzpicture}
\caption{The variants $\TL_4^\si$ for each $\si\in D_2(\TL_4)$.}
\label{f:TL4}
\end{center}
\end{figure}

\begin{rem}\label{r:TL4}
Lemma \ref{l:iso1} does not hold in any of the categories $\P$, $\PB$, $\PP$, $\M$ or $\TL$.  Example~\ref{e:max}\ref{emax7} gives two elements $\si,\tau$ from $\TL_4=\TL_{4,4}$ with $\rank(\si)=\rank(\tau)=2$ for which $\TL_4^\si\not\cong\TL_4^\tau$; cf.~Figure \ref{f:emax7}.  There are nine partitions from $\TL_4$ of rank $2$, and Figure \ref{f:TL4} gives eggbox diagrams of the corrsponding nine variants of $\TL_4$, all produced using GAP \cite{GAP}.  Up to isomorphism, there are five such variants; up to isomorphism and anti-isomorphism, there are four.  It would be interesting to attempt to classify the isomorphism classes of sandwich semigroups in $\TL$, and indeed in all the other diagram~categories.
\end{rem}

The next result gives another situation in which distinct sandwich semigroups can be isomorphic.  There is an obvious dual, but we will not state it.  Here $\Z^+=\{1,2,3,\ldots\}$ is the set of positive integers.

\begin{lemma}\label{l:iso2}
If $q\in\Z^+$, then
\ben
\item \label{iso21} $\B_{2q,0}^\si$ is a left-zero semigroup of size $(2q-1)!!$ for any $\si\in\B_{0,2q}$,
\item \label{iso22} $\B_{2q-1,1}^\si$ is a left-zero semigroup of size $(2q-1)!!$ for any $\si\in\B_{1,2q-1}$.
\een
\end{lemma}

\pf
We just prove \ref{iso21}, as \ref{iso22} is virtually identical.  If $\al,\be\in\B_{2q,0}$, then ${\al\star_\si\be=\al\si\be=\al\id_0=\al}$, so $\B_{2q,0}^\si$ is a left-zero semigroup.  The size follows from Proposition \ref{p:sizes}.
\epf

We are now ready to prove the classification result.

\begin{thm}\label{t:iso}
Let $m,n,k,l\in\N$ with $m\equiv n\Mod2$ and $k\equiv l\Mod2$, and let $\si\in\B_{nm}$ and $\tau\in\B_{lk}$ with $r=\rank(\si)$ and $s=\rank(\tau)$.  Then $\Bmns \cong \Bklt$ if and only if one of the following holds:
\ben
\item \label{iso1} $(m,n,r)=(k,l,s)$, 
\item \label{iso2} $m+n\leq2$ and $k+l\leq2$,
\item \label{iso3} renaming if necessary, $(m,n,r)=(2q,0,0)$ and $(k,l,s)=(2q-1,1,1)$ for some $q\in\Z^+$,
\item \label{iso4} renaming if necessary, $(m,n,r)=(0,2q,0)$ and $(k,l,s)=(1,2q-1,1)$ for some $q\in\Z^+$.
\een
\end{thm}

\pf
For the forwards implication,
\bit
\item \ref{iso1}$\implies$$\Bmns \cong \Bklt$ by Lemma \ref{l:iso1},
\item \ref{iso2}$\implies$$\Bmns \cong \Bklt$, because both semigroups have size $1$,
\item (\ref{iso3} or \ref{iso4})$\implies$$\Bmns \cong \Bklt$ by Lemma \ref{l:iso2} and its dual.
\eit
Conversely, suppose $\Bmns \cong \Bklt$, and suppose further that \ref{iso2} does not hold.  By Proposition \ref{p:sizes}, we have ${(m+n-1)!!=(k+l-1)!!}$.  Since $m+n$ and $k+l$ are not both $\leq2$, and since $x!!$ is strictly increasing for odd $x\geq1$, it follows that $m+n=k+l$.  We consider two cases.

\pfcase{1}
Suppose first that $(m,n)=(k,l)$.  Then $r\equiv s\Mod2$.  The largest group $\H$-classes in $\Bmns$ and $\Bklt$ are isomorphic to $\S_r$ and $\S_s$ (cf.~Proposition \ref{p:DP}\ref{DP2}), and these have sizes $r!$ and $s!$, respectively.  Thus, $r!=s!$.  Since $x!$ is strictly increasing for even $x\geq0$, and also for odd $x\geq1$, it follows that $r=s$, so that \ref{iso1} holds in this case.

\pfcase{2}
Now suppose $(m,n)\not=(k,l)$, say with $m>k$ (and $n<l$).  Let $z\in\{0,1\}$ be such that $z\equiv m\Mod2$.  By Proposition \ref{p:min}, $D_z(\Bmn)$ is the minimal ideal of $\Bmns$.  By Propositions \ref{p:min} and~\ref{p:D}\ref{D1} (and the dual of the latter), 
\[
|D_z(\Bmn)/{\R^\si}| =|D_z(\Bmn)/{\R}| = f(m) \AND |D_z(\Bmn)/{\L^\si}| = |D_z(\Bmn)/{\L}| = f(n),
\]
where $f:\N\to\N$ is defined by
\[
f(x) = \begin{cases}
x!! &\text{if $x$ is odd}\\
(x-1)!! &\text{if $x$ is even.}
\end{cases}
\]
Similarly, if $w\in\{0,1\}$ is such that $w\equiv k\Mod2$, then the minimal ideal of $\Bklt$ is $D_w(\B_{kl})$, and we have
\[
|D_z(\B_{kl})/{\R^\tau}| = f(k) \AND |D_z(\B_{kl})/{\L^\tau}| = f(l).
\]
Since $\Bmns\cong\Bklt$, it follows that $f(m)=f(k)$ and $f(n)=f(l)$.  For $x,y\in\N$ with $x<y$, we have $f(x)=f(y)$ if and only if $(x,y)=(0,1)$ or $(2q-1,2q)$ for some $q\in\Z^+$.  Since $m>k$ and $n<l$, it follows that
\bit
\item $(m,k)=(1,0)$ or $(2q,2q-1)$ for some $q\in\Z^+$, and
\item $(n,l)=(0,1)$ or $(2p-1,2p)$ for some $p\in\Z^+$.
\eit
Since $m\equiv n\Mod2$ and $k\equiv l\Mod2$, it follows that
\bit
\item $(m,k)=(1,0)$ and $(n,l)=(2p-1,2p)$ for some $p\in\Z^+$, or
\item $(m,k)=(2q,2q-1)$ and $(n,l)=(0,1)$ for some $q\in\Z^+$.
\eit
Thus, one of \ref{iso3} or \ref{iso4} holds in this case.
\epf

\begin{rem}
One could readily deduce a classification up to isomorphism and anti-isomorphism using the fact that $\Bmns$ and $\B_{nm}^{\si^*}$ are anti-isomorphic.
\end{rem}

\subsection{The rank of $\Bmns$}\label{ss:RankB}

Recall that the \emph{rank} of a semigroup $S$, denoted $\rank(S)$, is the smallest cardinality of a generating set for $S$.  In this section we calculate the rank of an arbitrary sandwich semigroup in the Brauer category $\B$.  To this end, we fix $m,n\in\N$ with $m\equiv n\Mod2$, and some $\si\in\B_{nm}$.  We may assume by symmetry that $m\geq n$, and we write $r=\rank(\si)$.  Further, by Lemma \ref{l:iso1}, we may assume without loss of generality that 
\begin{equation}\label{e:siB}
\si = \partI{1}{r}{r+1,r+2}{n-1,n}{1}{r}{r+1,r+2}{m-1,m}.
\end{equation}
The value of $\rank(\Bmns)$ depends on whether $r<n$ or $r=n$; see Theorems~\ref{t:rankB1} and~\ref{t:rankB2}.  Note that if $r=n=m$, then $\si=\id_n$, so that $\Bmns=\B_n$ is the Brauer monoid of degree $n$.  Since the rank of $\B_n$ is known to be $3$ for $n\geq3$ (cf.~\cite[Lemma 2.1]{Mazorchuk2002}), we will exclude the case of $r=n=m$.

We begin with a lemma that will be of use in both cases.  Throughout Section \ref{ss:RankB}, we write $D_q=D_q(\Bmn)$ for each $0\leq q\leq n$ with $q\equiv n\Mod2$.

\begin{lemma}\label{l:q+2}
If $\al\in D_q$, where $q\leq r$ and $q<n$, then $\al=\be\star_\si\ga$ for some $\be,\ga\in D_{q+2}$.
\end{lemma}

\pf
Write $\al = \partI{a_1}{a_q}{c_1,d_1}{c_s,d_s}{b_1}{b_q}{e_1,f_1}{e_t,f_t}$.  Since $q<n$, we have $s,t\geq1$.  Then with
\[
\be = \Big( 
{ \scriptsize \renewcommand*{\arraystretch}{1}
\begin{array} {\c|\c|\c|\c|\c|\c|\c|\cend}
a_1 \:&\: \cdots \:&\: a_q \:&\: c_s\:&\: d_s \:&\: c_1,d_1 \:& \: \cdots \:&\: c_{s-1},d_{s-1} \\ \cline{6-8}
1 \:&\: \cdots \:&\: q \:&\: n-1 \:&\: n \:&\: q+1,q+2 \:&\: \cdots \:&\:n-3,n-2
\rule[0mm]{0mm}{2.7mm}
\end{array} 
}
\hspace{-1.5 truemm} \Big)
\ANd
\ga = \Big( 
{ \scriptsize \renewcommand*{\arraystretch}{1}
\begin{array} {\c|\c|\c|\c|\c|\c|\c|\cend}
1 \:&\: \cdots \:&\: q \:&\: m-1 \:&\: m \:&\: q+1,q+2 \:&\: \cdots \:&\:m-3,m-2 \\ \cline{6-8}
b_1 \:&\: \cdots \:&\: b_q \:&\: e_t\:&\: f_t \:&\: e_1,f_1 \:& \: \cdots \:&\: e_{t-1},f_{t-1}
\rule[0mm]{0mm}{2.7mm}
\end{array} 
}
\hspace{-1.5 truemm} \Big),
\]
one may easily check (with separate diagrams for the $r=n<m$ and $r<n\leq m$ cases) that $\al=\be\si\ga$.
\epf

\begin{thm}\label{t:rankB1}
If $r<n\leq m$, then $\Bmns = \la\Om\ra$, where $\Om=\set{\al\in\Bmn}{\rank(\al)>r}$.  Furthermore, every generating set for $\Bmns$ contains $\Om$, and so
\[
\rank(\Bmns) = |\Om| = \sum_{r<q\leq n \atop r\equiv n\Mod2} \tbinom mq\tbinom nq(m-q-1)!! (n-q-1)!!q!.
\]
\end{thm}

\pf
Note that $\Om = D_{r+2}\cup D_{r+4}\cup\cdots\cup D_n$, so the formula for $|\Om|$ follows from the formulae for~$|D_q|$ given in Remark \ref{r:DrB}.  By Proposition \ref{p:maxPPBB}\ref{maxPPBB1}, $\Om$ is the union of all the maximal $\J^\si$-classes of $\Bmns$, and all of these are singletons; it follows immediately that every generating set for $\Bmns$ contains $\Om$.

It remains to show that $\Bmns=\la\Om\ra$.  To do this, we must show that $\la\Om\ra$ contains $D_q$ for each $0\leq q\leq n$ with $q\equiv n\Mod2$.  We do this by descending induction on $q$, noting that $D_q\sub\Om\sub\la\Om\ra$ for $q\geq r+2$.  Now suppose $q\leq r$, and fix some $\al\in D_q$.  By Lemma \ref{l:q+2}, we have $\al=\be\star_\si\ga$ for some $\be,\ga\in D_{q+2}$.  By induction we have $\be,\ga\in\la\Om\ra$, so it follows that $\al=\be\star_\si\ga\in\la\Om\ra$, as required.
\epf

We now treat the case in which $r=n<m$.  Here we have 
$
\si= \Big( 
{ \scriptsize \renewcommand*{\arraystretch}{1}
\begin{array} {\c|\c|\c|\c|\c|\cend}
 1 \:&\: \cdots \:&\: n \:& \multicolumn{3}{c}{} \\ \cline{4-6}
 1 \:&\: \cdots \:&\: n \:&\: n+1,n+2 \:&\: \cdots \:&\: m-1,m 
\rule[0mm]{0mm}{2.7mm}
\end{array} 
}
\hspace{-1.5 truemm} \Big)$, and we note that $\si$ is right-invertible; indeed, $\si\si^*=\id_n$.

\begin{thm}\label{t:rankB2}
If $r=n<m$, then $\rank(\Bmns) = \binom mn(m-n-1)!!$.
\end{thm}

\pf
We first claim that $\Bmns=\la D_n\ra$.  Indeed, it may be proved by descending induction that $D_q\sub\la D_n\ra$ for all $0\leq q\leq n$ with $q\equiv n\Mod2$.  The $q=n$ case is obvious, and the inductive step again follows from Lemma \ref{l:q+2}.

Next we note that $D_n=J_{\si^*}$ is the maximum $\J$-class in the hom-set $\Bmn$ (cf.~Lemma \ref{l:max4} and Corollary \ref{c:D}).  Since~$\si^*$ is a right-inverse of $\si$, $D_n^\si=J_{\si^*}^\si$ is the maximum $\J^\si$-class in $\Bmns$, and this is a left-group over $H_{\si^*}^\si\cong\S_n$ (cf.~Propositions \ref{p:max3} and \ref{p:maxPPBB}\ref{maxPPBB2}).  Now,
\[
\be = \Big( 
{ \scriptsize \renewcommand*{\arraystretch}{1}
\begin{array} {\c|\c|\c|\c|\c|\c|\cend}
 1 \:&\: \cdots \:&\: n-1 \:&\: m \:&\: n,n+1 \:&\: \cdots \:&\: m-2,m-1  \\ \cline{5-7}
 1 \:&\: \cdots \:&\: n-1 \:&\: n  \:& \multicolumn{3}{c}{}
\rule[0mm]{0mm}{2.7mm}
\end{array} 
}
\hspace{-1.5 truemm} \Big)
\]
is a right-inverse of $\si$, and hence an idempotent of $D_n^\si$.  Since $\be\not=\si^*$ (as $n<m$), it follows that ${|D_n^\si/{\H^\si}|\geq2\geq\rank(\S_n)}$.  It follows from Proposition \ref{p:RI}\ref{pRI2} and Proposition \ref{p:D}\ref{D1} that
\[
\rank(\Bmns)=\rank(\la D_n\ra)=|D_n/{\H}|=|D_n/{\R}|=\tbinom mn(m-n-1)!!. \qedhere
\]
\epf

\begin{rem}
Although there is not generally a unique generating set for $\Bmns$ of minimum cardinality for $r=n<m$, it is possible to classify the generating sets of minimum cardinality (modulo a classification of the generating sets of the symmetric group $\S_n$).  Indeed, following the proof of Proposition~\ref{p:RI}, these are all of the form $X_1\cup X_2$, where:
\bit
\item $X_1$ is a cross-section of the $\H$-classes of $\Bmn$ contained in $D_n^\si$ such that $D_n^\si=\la X_1\ra$, and
\item $X_2$ is an arbitrary cross-section of the $\H$-classes of $\Bmn$ contained in $D_n\sm D_n^\si$.
\eit
Writing $k=|D_n^\si/{\H}|$, sets $X_1$ of the above form are in one-one correspondence with the ordered lists $(\pi_1,\ldots,\pi_k)$ where $\S_n$ is generated by $\{\pi_1,\ldots,\pi_k\}$; see the proof of \cite[Proposition 4.11]{Sandwich1}.
\end{rem}

\newpage

\begin{rem}
Theorems \ref{t:rankB1} and \ref{t:rankB2} do not have direct analogues for the other diagram categories.  In general:
\bit
\item If $r<n\leq m$, then $\Kmns$ is not generated by its maximal $\J^\si$-classes (regardless of whether there is a nontrivial maximal $\J^\si$-class or not).
\item If $r=n<m$, then $\Kmns$ is not generated by $D_n$.
\eit
All of this may be verified with GAP \cite{GAP}.  Calculating the ranks of sandwich semigroups in other diagram categories appears to be an interesting (and challenging) problem.
\end{rem}

\subsection{A combinatorial digression}\label{ss:ka}

In this section we introduce certain numbers that will play an important role in all that follows.

Let $\ve$ be an equivalence on a set $X$.  We say $\ve$ is a \emph{2-equivalence} if every $\ve$-class has size $2$.  If $|X|=m$, then there are $(m-1)!!$ 2-equivalences on $X$; recall that we define $(-1)!!=1$, and $k!!=0$ if~$k$ is even.  We say $\ve$ is a \emph{1-2-equivalence} if each $\ve$-class has size 1 or 2.  The \emph{rank} of a 1-2-equivalence~$\ve$, denoted $\rank(\ve)$, is defined to be the number of singleton $\ve$-classes; note that $\rank(\ve)\equiv |X|\Mod2$ if $|X|$ is finite.  (Of course the 1-2-equivalences on $[m]$ are simply the kernels of elements of $\Bmn$; for $\al\in\Bmn$, the rank of the 1-2-equivalence $\ker(\al)$ is just $\rank(\al)$.)

For $m,q\in\N$ with $q\leq m$ and $q\equiv m\Mod2$, let $\ka(m,q)$ be the number of 1-2-equivalences $\ve$ of an $m$-element set with $\rank(\ve)=q$.  Then
\begin{equation}\label{e:kamq}
\ka(m,q)=\tbinom mq(m-q-1)!!.
\end{equation}
Indeed, to specify such an $\ve$, we first choose the~$q$ singletons in $\binom mq$ ways, and then pair the remaining $m-q$ elements of $X$ in $(m-q-1)!!$ ways (cf.~Proposition \ref{p:D}\ref{D1}.)

Suppose $m,r,q\in\N$ are such $q\leq r\leq m$ and $q\equiv r\equiv m\Mod2$.  Fix a set $X$ with $|X|=m$, and a 1-2-equivalence $\eta$ on $X$ with $\rank(\eta)=r$.  We define $\ka(m,r,q)$ to be the number of 1-2-equivalences~$\ve$ on $X$ such that:
\bit
\item $\rank(\ve)=q$, and
\item the join $\ve\vee\eta$ has precisely $q$ classes of odd size.
\eit
Clearly the definition of $\ka(m,r,q)$ depends only on the cardinality of the set $X$ and the rank of $\eta$.  (Note that if $\ve$ is an arbitrary 1-2-equivalence on $X$, then every $\ve\vee\eta$-class is a union of $\ve$-classes, at most two of which can be singletons.)

When working with 1-2-equivalences below, one can visualise a join $\ve\vee\eta$ as follows.  We draw the vertices from $X$ on a horizontal row, and indicate the nontrivial $\ve$-classes (or $\eta$-classes) by drawing edges below (or above) the vertices, respectively; the $\ve\vee\eta$-classes are then the connected components of this graph.

\begin{lemma}\label{l:ka1}
If $m,r\in\N$ are such that $r\leq m$ and $r\equiv m\Mod2$, then $\ka(m,r,r)=\frac{(m+r-1)!!}{(2r-1)!!}$.
\end{lemma}

\pf
Without loss of generality we may assume that $X=\{1,\ldots,m\}$, and that the $\eta$-classes are $\{1\},\ldots,\{r\}$ and $\{r+1,r+2\},\ldots,\{m-1,m\}$.  
Define the numbers $\lam(m,r)$, for $m,r\in\N$ with $r\leq m$ and $r\equiv m\Mod2$, as follows:
\ben
\item \label{la1} $\lam(m,r)=(m-1)!!$ if $r=0$,
\item \label{la2} $\lam(m,r)=1$ if $r=m$,
\item \label{la3} $\lam(m,r)=\lam(m-1,r-1)+(m-r)\lam(m-2,r)$ if $0<r<m$.
\een
We may prove the lemma by showing that $\ka(m,r,r)$ and $\frac{(m+r-1)!!}{(2r-1)!!}$ both satisfy the recurrence for~$\lam(m,r)$.  The latter is easily checked.  

\begin{figure}[ht]
\begin{center}
\scalebox{0.85}{
\begin{tikzpicture}[scale=1.3]
\tikzstyle{vertex}=[circle,draw=black, fill=white, inner sep = 0.03cm,minimum size=25pt]
\begin{scope}[shift={(0,-2)}]
\uarcx19{.7}
\end{scope}
\darcx45{.5}
\darcx78{.5}
\darcx9{10}{.5}
\darcx{11}{12}{.5}
\darcx{14}{15}{.5}
\foreach \x/\y in {1/1,3/r,9/a,15/m} {\node[vertex] (\x) at (\x,0) {\large $\y$};}
\foreach \x/\y in {
4/r\!+\!1,
5/r\!+\!2,
7/a\!-\!2,
8/a\!-\!1,
10/a\!+\!1,
11/a\!+\!2,
12/a\!+\!3,
14/m\!-\!1
} {\node[vertex] (\x) at (\x,0) {\small $\y$};}
\draw[dotted] (1.5,0)--+(1,0);
\draw[dotted] (5.5,0)--+(1,0);
\draw[dotted] (12.5,0)--+(1,0);
\draw[-{latex}] (8,-1)--(8,-2.2);
\begin{scope}[shift={(0,-3)}]
\darcx45{.5}
\darcx78{.5}
\darcx{11}{12}{.5}
\darcx{14}{15}{.5}
\foreach \x/\y in {1/1,3/r,15/m} {\node[vertex] (\x) at (\x,0) {\large $\y$};}
\foreach \x/\y in {
4/r\!+\!1,
5/r\!+\!2,
7/a\!-\!2,
8/a\!-\!1,
11/a\!+\!2,
12/a\!+\!3,
14/m\!-\!1
} {\node[vertex] (\x) at (\x,0) {\small $\y$};}
\draw[dotted] (1.5,0)--+(1,0);
\draw[dotted] (5.5,0)--+(1,0);
\draw[dotted] (12.5,0)--+(1,0);
\end{scope}
\end{tikzpicture}
}
\caption{Graphs constructed during the proof of Lemma \ref{l:ka1}.}
\label{f:etave}
\end{center}
\end{figure}
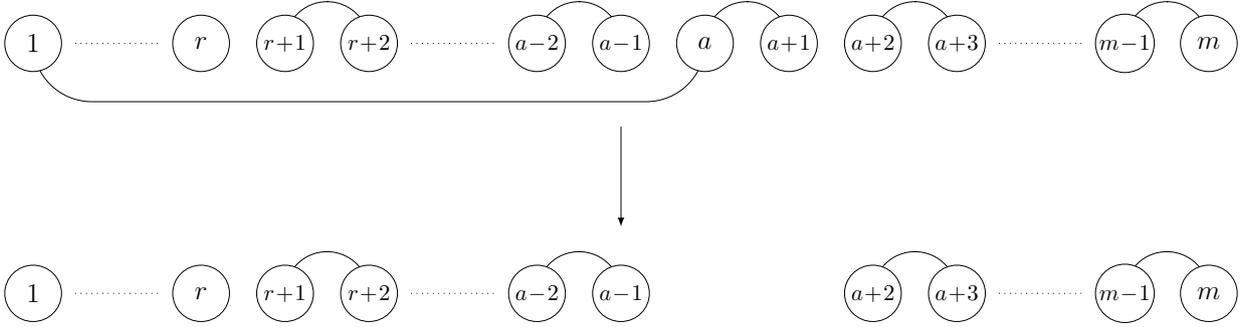

For the former, first note that \ref{la1} holds because $\ka(m,0,0)$ is simply the number of 2-equivalences on $X$.  Item \ref{la2} is clear.  For \ref{la3}, suppose $0<r<m$.  Consider a 1-2-equivalence $\ve$ of the desired form, and let $A$ be the $\ve$-class of~$1$.  Then one of the following holds:
\[
A=\{1\} \OR A=\{1,a\} \text{ for some $r<a\leq m$.}
\]
There are $\ka(m-1,r-1,r-1)$ and $(m-r)\ka(m-2,r,r)$ possibilities for $\ve$ in these two cases, respectively.  Adding these gives \ref{la3}.  
To see why there are $(m-r)\ka(m-2,r,r)$ possibilities for $\ve$ in the second case, we argue as follows.  We first draw the edges from $\eta$ above the vertices $1,\ldots,m$ as indicated in the top diagram from Figure \ref{f:etave}.  Once $a\in\{r+1,\ldots,m\}$ is chosen, we draw the edge $\{1,a\}$ below the vertices, as also indicated in Figure \ref{f:etave}.  The number of ways to complete this picture to obtain a graph with $r$ odd-sized components is the same as the number of ways to complete the bottom diagram from Figure \ref{f:etave} to obtain a graph with $r$ odd-sized components.  By definition, there are $\ka(m-2,r,r)$ ways to do the latter.
\epf

\begin{lemma}\label{l:ka}
If $m,r,q\in\N$ are such $q\leq r\leq m$ and $q\equiv r\equiv m\Mod2$, then
\[
\ka(m,r,q) = \binom rq \frac{(r-q-1)!!(m+q-1)!!}{(r+q-1)!!}.
\]
\end{lemma}

\pf
Again we will assume that $X=\{1,\ldots,m\}$, and that the $\eta$-classes are $\{1\},\ldots,\{r\}$ and ${\{r+1,r+2\},\ldots,\{m-1,m\}}$.  Define the numbers $\lam(m,r,q)$, for $m,r,q\in\N$ with $q\leq r\leq m$ and $q\equiv r\equiv m\Mod2$, as follows:
\ben
\item \label{ka1} $\lam(m,r,q)=(m-1)!!$ if $q=0$,
\item \label{ka3} $\lam(m,r,q)=\binom mq(m-q-1)!!$ if $m=r$,
\item \label{ka4} $\lam(m,r,q)=\frac{(m+r-1)!!}{(2r-1)!!}$ if $r=q$,
\item \label{ka5} $\lam(m,r,q)=\lam(m-1,r-1,q-1)+(r-1)\lam(m-2,r-2,q)+(m-r)\lam(m-2,r,q)$ for $0<q<r<m$.
\een
We will show that $\ka(m,r,q)$ and $\binom rq \frac{(r-q-1)!!(m+q-1)!!}{(r+q-1)!!}$ both satisfy the recurrence for $\lam(m,r,q)$.  Again the latter is easily checked, though the calculation is somewhat involved for \ref{ka5}.  The rest of the proof is devoted to the former.

\pfitem{\ref{ka1} and \ref{ka3}}  If $r=m$ or if $q=0$, then $\ve\vee\eta$ has $q$ odd-sized blocks for \emph{every} 1-2-equivalence~$\ve$ with $\rank(\ve)=q$; thus, $\ka(m,r,q)=\ka(m,q)$, and we then apply \eqref{e:kamq}.

\pfitem{\ref{ka4}}  This was proved in Lemma \ref{l:ka1}.

\pfitem{\ref{ka5}}  Suppose $0<q<r<m$, consider a 1-2-equivalence $\ve$ of the desired form, and let $A$ be the $\ve$-class of $1$.  Then one of the following holds:
\[
A=\{1\} \OR A=\{1,a\} \text{ for some $1<a\leq r$} \OR A=\{1,a\} \text{ for some $r<a\leq m$.}
\]
As in the proof of Lemma \ref{l:ka1}, there are $\ka(m-1,r-1,q-1)$, $(r-1)\ka(m-2,r-2,q)$ and ${(m-r)\ka(m-2,r,q)}$ possibilities for $\ve$ in these three cases, respectively.  Adding these gives \ref{ka5}.
\epf

\subsection{The structure of the regular subsemigroup of $\Bmns$}\label{ss:RegB}

Throughout Section \ref{ss:RegB} we fix $m,n\in\N$ with $m\equiv n\Mod2$, and some $\si\in\B_{nm}$ of rank $r$.  Our goal is to give a structural description of the regular subsemigroup $P^\si=\Reg(\Bmns)$ of the sandwich semigroup $\Bmns$.  It is of no advantage to assume that $m\geq n$ holds, or that $m=n=r$ does not hold, so we make no such assumption during this section.

Let $\al\in\Bmn$, and consider the product graph $\Pi(\al,\si)$; this graph is defined in Section \ref{ss:P}, and it is used in the construction of the product $\al\si$.  Because of the form of Brauer partitions, every component in the graph $\Pi(\al,\si)$ is a path or a loop, and these can be of the following types:
\begin{enumerate}[label=\textup{(C\arabic*)},leftmargin=9mm]
\item \label{c1} $x\lar\al z$ for some $x,z\in[m]$,
\item \label{c2} $x'\lar\si z'$ for some $x,z\in[m]$,
\item \label{c3} $x\lar\al y_1''\lar\si y_2''\lar\al \cdots\lar\si y_{2k}''\lar\al z$ for some $x,z\in[m]$, $k\in\Z^+$ and $y_1,y_2,\ldots,y_{2k}\in[n]$,
\item \label{c4} $x'\lar\si y_1''\lar\al y_2''\lar\si \cdots\lar\al y_{2k}''\lar\si z'$ for some $x,z\in[m]$, $k\in\Z^+$ and $y_1,y_2,\ldots,y_{2k}\in[n]$,
\item \label{c5} $y_1''\lar\al y_2''\lar\si\cdots\lar\al y_{2k}''\lar\si y_1''$ for some $k\in\Z^+$ and $y_1,y_2,\ldots,y_{2k}\in[n]$,
\item \label{c6} $x\lar\al y_1''\lar\si y_2''\lar\al \cdots\lar\al y_{2k-1}''\lar\si z'$ for some $x,z\in[m]$, $k\in\Z^+$ and ${y_1,y_2,\ldots,y_{2k-1}\in[n]}$.
\een
Components of types \ref{c1} and \ref{c3} lead to upper nontransversals $\{x,z\}$ in the product $\al\si$; components of types \ref{c2} and \ref{c4} lead to lower nontransversals $\{x',z'\}$ of $\al\si$; components of type \ref{c5} are closed loops lying in the middle row of $\Pi(\al,\si)$, and these are discarded when forming $\al\si$; finally, components of type \ref{c6} lead to transversals $\{x,z'\}$ of $\al\si$.  Also note that every equivalence class of the join $\coker(\al)\vee\ker(\si)$ is of the form
\bit
\item $\{y_1,\ldots,y_{2k}\}$ for some component of $\Pi(\al,\si)$ of type \ref{c3}, \ref{c4} or \ref{c5}, or
\item $\{y_1,\ldots,y_{2k-1}\}$ for some component of $\Pi(\al,\si)$ of type \ref{c6}.
\eit
In particular, keeping Proposition \ref{p:PK} in mind, and writing $q=\rank(\al)$,
\begin{align*}
\al\in P_1^\si \iff \rank(\al\si)=q &\iff \text{$\Pi(\al,\si)$ has no components of type \ref{c3}}\\
&\iff \text{$\Pi(\al,\si)$ has $q$ components of type \ref{c6}}\\
&\iff \text{$\coker(\al)\vee\ker(\si)$ has $q$ classes of odd size}\\
&\iff \text{$\coker(\al)\vee\ker(\si)$ separates $\codom(\al)$.}
\end{align*}
(Recall that an equivalence $\ve$ on a set $X$ \emph{separates} a subset $A$ of $X$ if each $\ve$-class contains at most one element of $A$.)  Combining this with the dual statement concerning products of the form $\si\al$, we have proved the following:

\begin{prop}\label{p:PB}
We have
\begin{align*}
P_1^\si &= \bigset{\al\in\Bmn}{\text{$\coker(\al)\vee\ker(\si)$ separates $\codom(\al)$}} , \\
P_2^\si &= \bigset{\al\in\Bmn}{\text{$\ker(\al)\vee\coker(\si)$ separates $\dom(\al)$}} , \\
\Reg(\Bmns) = P^\si = P_3^\si &= \big\{\al\in\Bmn:\text{$\coker(\al)\vee\ker(\si)$ separates $\codom(\al)$}\\
\epfreseq &\hspace{3.6cm}\text{and $\ker(\al)\vee\coker(\si)$ separates $\dom(\al)$}\big\}.  
\end{align*}
\end{prop}

\begin{rem}
Of course Proposition \ref{p:PB} holds for the Temperley-Lieb category $\TL$ as well.  However, it does not hold for any of the other categories.  For example, with $\al = \custpartn{1,2}{1,2}{\stline11}$ and $\si = \custpartn{1,2}{1,2}{\uarc12\darc12}$, both from~$\M_2$ (and hence both from $\PB_2$, $\PP_2$ and $\P_2$), we have $\rank(\al\si)=\rank(\si\al)=0<\rank(\al)$, even though ${\coker(\al)\vee\ker(\si)=\ker(\al)\vee\coker(\si)}$ separates $\dom(\al)=\codom(\al)$.
\end{rem}

Recall from Section \ref{ss:nonsandwichK} that we have an surmorphism
\[
\Psi:P^\si=\Reg(\Bmns)\to\B_r:\al\mt(\si\al\si\si^*)^\natural = \tau\si\al\tau^*.
\]
Here $\tau\in\B_{rn}$ is as defined in \eqref{e:tau}.  In what follows, it will be convenient to write $\ol\al=\al\Psi$ for $\al\in P^\si$.  

As in Section \ref{ss:nonsandwich}, for each of Green's relations $\K$ we have a relation $\gKh^\si$ on $P^\si$ defined so that $\al \gKh^\si \be$ in $\Bmns$ if and only if $\ol\al\K\ol\be$ in $\KK_r$, and we recall that ${\gDh^\si}={\D^\si}$.  In particular, the $\D^\si$-classes $D_q^\si=D_q(P^\si)$ of $P^\si$ map onto the $\D$-classes $D_q(\B_r)$ of $\B_r$ for each $q$.

The next result uses the $\gKh^\si$ relations to give a finer description of the internal structure of the $\D^\si$-classes of $P^\si$.  In the statement and proof, we denote the $\gKh^\si$-class of $\al\in P^\si$ by $\Kh_\al^\si$.  By Lemma~\ref{l:iso1}, it suffices to assume that $\si$ has the form in \eqref{e:siB}, in which case
$
\tau = \Big( 
{ \scriptsize \renewcommand*{\arraystretch}{1}
\begin{array} {\c|\c|\c|\c|\c|\cend}
 1 \:&\: \cdots \:&\: r \:& \multicolumn{3}{c}{} \\ \cline{4-6}
 1 \:&\: \cdots \:&\: r \:&\: r,r+1 \:&\: \cdots \:&\: n-1,n 
\rule[0mm]{0mm}{2.7mm}
\end{array} 
}
\hspace{-1.5 truemm} \Big)
$.

\begin{thm}\label{t:RegB}
Let $0\leq q\leq r$ with $q\equiv r\Mod2$.
\ben
\item \label{RB1} $D_q^\si$ contains $\binom rq(r-q-1)!!$ $\gRh^\si$-classes, each of which contains $\frac{(m+q-1)!!}{(r+q-1)!!}$ $\R^\si$-classes.
\item \label{RB2} $D_q^\si$ contains $\binom rq(r-q-1)!!$ $\gLh^\si$-classes, each of which contains $\frac{(n+q-1)!!}{(r+q-1)!!}$ $\L^\si$-classes.
\item \label{RB3} $D_q^\si$ contains ${\binom rq}^{\!2}(r-q-1)!!^2$ $\gHh^\si$-classes, each of which contains $\frac{(m+q-1)!!(n+q-1)!!}{(r+q-1)!!^2}$ $\H^\si$-classes.
\item \label{RB4} Each $\H^\si$-class in $D_q^\si$ has size $q!$, and group $\H^\si$-classes in $D_q^\si$ are isomorphic to the symmetric group $\S_q$.
\item \label{RB5} An $\H^\si$-class $H_\al^\si\sub D_q^\si$ is a group if and only if $H_{\ol\al}\sub D_q(\B_r)$ is a group $\H$-class of $\B_r$, in which case $\Hh_\al^\si$ is a $\frac{(m+q-1)!!}{(r+q-1)!!}\times\frac{(n+q-1)!!}{(r+q-1)!!}$ rectangular group over $\S_q$.
\een
\end{thm}

\pf
\firstpfitem{\ref{RB1}}  By definition, the $\gRh^\si$-classes in $D_q^\si$ are in one-one correspondence with the $\R$-classes in $D_q(\B_r)$; by Proposition \ref{p:D}\ref{D1}, there are $\binom rq(r-q-1)!!$ of these.  

An $\R^\si$-class in $D_q^\si$ is uniquely determined by the common kernel of all its elements.  Such a kernel is a 1-2-equivalence on $[m]$ of rank $q$ (as defined in Section \ref{ss:ka}) whose join with $\coker(\si)$ has $q$ classes of odd size (cf.~Proposition \ref{p:PB} and the discussion before it).  By definition, and using Lemma \ref{l:ka}, there are $\ka(m,r,q)=\binom rq \frac{(r-q-1)!!(m+q-1)!!}{(r+q-1)!!}$ such equivalences.

Dividing $\ka(m,r,q)$ by the number of $\gRh^\si$-classes gives $\frac{(m+q-1)!!}{(r+q-1)!!}$, so the proof of \ref{RB1} will be complete if we can show that each $\gRh^\si$-class contains the same number of $\R^\si$-classes.  To do so, consider two $\gRh^\si$-classes $\Rh_\al^\si$ and $\Rh_\be^\si$, where $\al,\be\in D_q^\si$.  
Since $\Rh_\al^\si$ and $\Rh_\be^\si$ are unions of $\R^\si$-classes, it is enough to show that these $\gRh^\si$-classes have the same size.

Since $\al\D^\si\be$, we have $\al\R^\si\ga\L^\si\be$ for some $\ga\in D_q^\si$.  Since $\ga\in R_\al^\si\sub\Rh_\al^\si$, we have $\Rh_\al^\si=\Rh_\ga^\si$, so we may in fact assume without loss of generality that $\al\L^\si\be$.  It follows that $\ol\al\L\ol\be$ in $\B_r$, so that $\coker(\ol\al)=\coker(\ol\be)$, whence $\ol\be=\pi\ol\al$ for some permutation $\pi\in\S_r$.  Define the permutations
\[
\rho = \partpermVI1r{r+1}m{1\pi}{r\pi}{r+1}m\in\S_m \AND \varrho = \partpermVI1r{r+1}n{1\pi}{r\pi}{r+1}n\in\S_n.
\]
It is easily checked that $\tau\si\cdot\rho=\pi\cdot\tau\si$ and $\si\rho=\varrho\si$.  

We claim that there is a well-defined map
\[
\th:\Rh_\al^\si\to\Rh_\be^\si:\ga\mt\rho\ga .
\]
To prove this, we must show that $\th$ does indeed map $\Rh_\al^\si$ into $\Rh_\be^\si$.  To do so, let $\ga\in\Rh_\al^\si$; we must show that 
\bena\bmc2
\item \label{xi1} $\ga\th\in P^\si$, and 
\item \label{xi2} $\ga\th\gRh^\si\be$: i.e., $\ol{\ga\th} \R \ol\be$.
\emc\een
Since $\ga\in\Rh_\al^\si$, it follows in particular that $\ga\in P^\si$, so that $\rank(\si\ga\si)=\rank(\ga)$.  Since $\rho\in\S_m$ and $\varrho\in\S_n$, and using $\si\rho=\varrho\si$, we have
\[
\rank(\si(\ga\th)\si) = \rank(\si\rho\ga\si) = \rank(\varrho\si\ga\si) = \rank(\si\ga\si) = \rank(\ga) = \rank(\rho\ga) = \rank(\ga\th),
\]
which verifies \ref{xi1}.  For \ref{xi2}, using the facts that $\tau\si\cdot\rho=\pi\cdot\tau\si$, that $\ol\ga\R\ol\al$ (as $\ga\in\Rh_\al^\si$), and that $\R$ is a left-congruence, we see that
\[
\ol{\ga\th} = \ol{\rho\ga} = \tau\si(\rho\ga)\tau^* = \pi(\tau\si\ga\tau^*)=\pi\ol\ga \R \pi\ol\al = \ol\be.
\]
This completes the proof that $\th$ is well defined.  Since $\rho$ is a permutation, $\th$ is injective, so it follows that $|\Rh_\al^\si|\leq|\Rh_\be^\si|$.  The reverse inequality follows by symmetry, and so $|\Rh_\al^\si|=|\Rh_\be^\si|$.  As noted above, this completes the proof of \ref{RB1}.

\pfitem{\ref{RB2}}  This is dual to \ref{RB1}.

\pfitem{\ref{RB3}}  This follows immediately from \ref{RB1} and \ref{RB2}.

\pfitem{\ref{RB4}}  This was proved in Proposition \ref{p:DP}\ref{DP2}.

\pfitem{\ref{RB5}}  This follows from Theorem \ref{t:RG}\ref{RG2} and \ref{RG3}, and parts \ref{RB1}, \ref{RB2} and \ref{RB4}.
\epf

\begin{rem}\label{r:B}
Figure \ref{f:B} gives eggbox diagrams for the regular semigroups $\Reg(\B_{66}^{\si_1})$ and $\Reg(\B_{64}^{\si_2})$, where $\si_1\in\B_{66}$ and $\si_2\in\B_{46}$ both have rank $4$; also pictured is the Brauer monoid $\B_4$; all were produced using GAP \cite{GAP}.  These illustrate the relationships between $\gKh^\si$-classes of $\Reg(\Bmns)$ and $\K$-classes of $\B_r$, as described in Theorem \ref{t:RegB}.  (See also Figures \ref{f:emax7} and \ref{f:emax8} for similar comparisons between $\Reg(\Kmns)$ and $\KK_r$ in the categories $\TL$ and~$\P$.  Note that in diagram categories other than $\B$ it is possible for $\gRh^\si$-classes in the same $\D^\si$-class to contain different numbers of $\R^\si$-classes.)

One may also compare Figure \ref{f:B} with the combinatorial formulae given in Theorem~\ref{t:RegB}.  For example, the $\D^{\si_1}$-classes $D_4^{\si_1}$, $D_2^{\si_1}$ and $D_0^{\si_1}$ of $\Reg(\B_{66}^{\si_1})$ contain $\ka(6,4,4)=9$, $\ka(6,4,2)=42$ and $\ka(6,4,0)=15$ $\R^{\si_1}$-classes, respectively.  
\end{rem}

\begin{figure}[ht]
\begin{center}
\includegraphics[height=16.6cm]{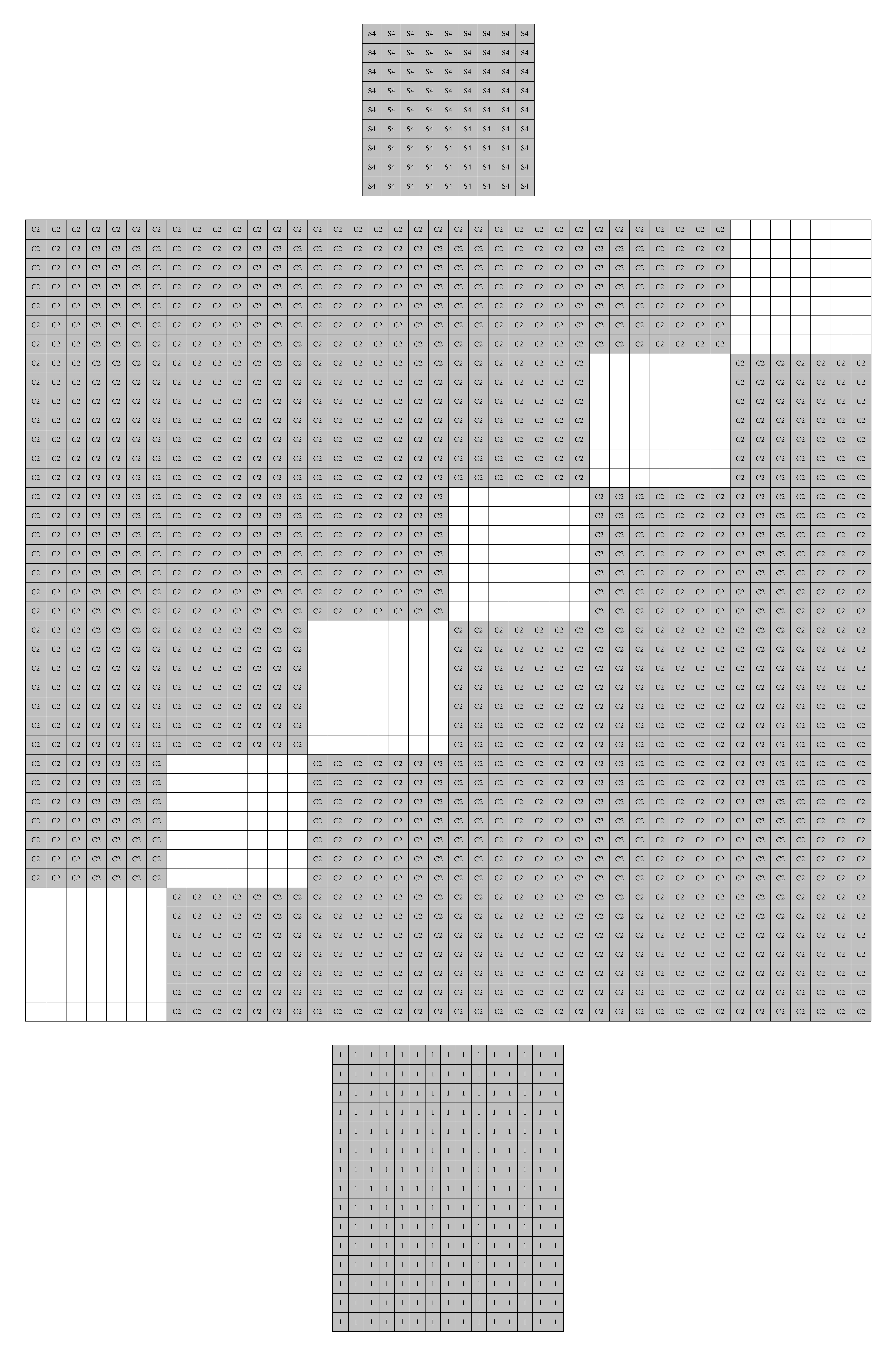} 
\qquad
\includegraphics[height=16.6cm]{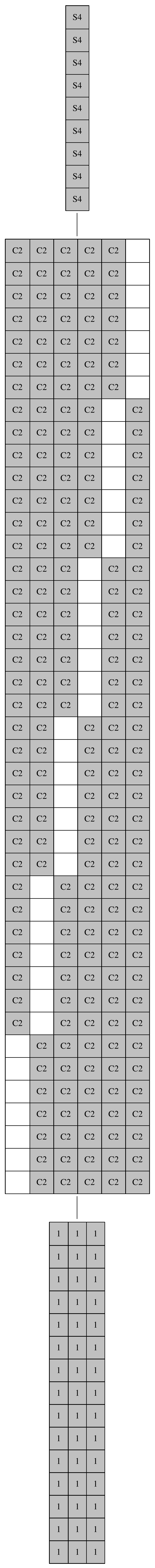} 
\qquad
\includegraphics[height=4.5cm]{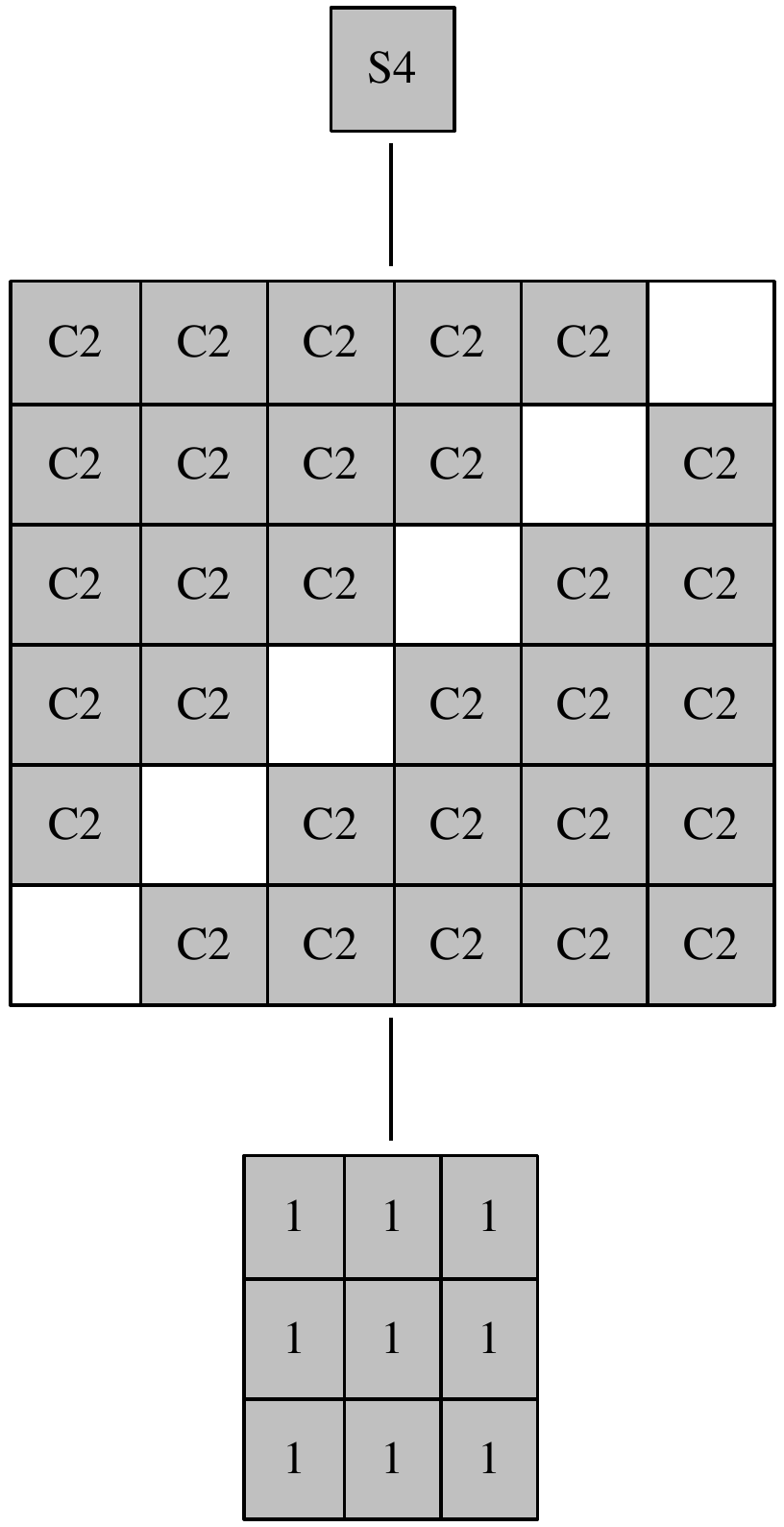} 
\caption[blah]{Left to right: eggbox diagrams of $\Reg(\B_{66}^{\si_1})$, $\Reg(\B_{64}^{\si_2})$ and $\B_4$; cf.~Remark \ref{r:B}.}
\label{f:B}
\end{center}
\end{figure}

Multiplying the two numbers from Theorem \ref{t:RegB}\ref{RB3} gives the number of $\H^\si$-classes in $D_q^\si$; multiplying this by $q!$ gives the size of $D_q^\si$; summing over appropriate $q$ gives the size of $P^\si$:

\begin{cor}\label{c:PB}
The size of the regular subsemigroup $P^\si=\Reg(\Bmns)$ is given by
\[
|P^\si| = \sum_{0\leq q\leq r \atop q\equiv r\Mod2} \ka(m,r,q)\ka(n,r,q)q! = \sum_{0\leq q\leq r \atop q\equiv r\Mod2} {\binom rq}^{\!\!2}\frac{(r-q-1)!!^2(m+q-1)!!(n+q-1)!!}{(r+q-1)!!^2}q!.  \epfreseq
\]
\end{cor}

As another application of Theorem \ref{t:RegB}, we may calculate the number of idempotents in the sandwich semigroup $\Bmns$.  As usual, for a subset $A$ of a semigroup $T$, we write $E(A)$ for the set of all idempotents of $T$ contained in $A$; we also write $e(A)=|E(A)|$ for the number of such idempotents.  Several formulae were given in \cite{DEEFHHL1} for the number $e(D_q(\B_r))$ of rank-$q$ idempotents from the Brauer monoid~$\B_r$.  Specifically, \cite[Theorem 30]{DEEFHHL1} says that $e(D_q(\B_r))=\binom rq(r-q-1)!!a_{rq}$, where $a_{rq}$ is defined by the recurrence
\[
a_{rr}=1\text{ for all $r$} \COMMA a_{r0}=(r-1)!!\text{ for even $r$}  \COMMA a_{rq}=a_{r-1,q-1}+(r-q)a_{r-2,q}\text{ for $1<q<r$.} 
\]
In fact, by comparing this with the proof of Lemma \ref{l:ka1}, we see that $a_{rq}=\ka(r,q,q)=\frac{(r+q-1)!!}{(2q-1)!!}$, so that in fact 
\begin{equation}\label{e:erq}
e(D_q(\B_r))=\binom rq\frac{(r-q-1)!!(r+q-1)!!}{(2q-1)!!}.
\end{equation}
Now, each idempotent of $D_q(\B_r)$ corresponds to a group $\H$-class in $D_q(\B_r)$; by Theorem \ref{t:RegB}\ref{RB5}, the preimage of this group $\H$-class is a $\frac{(m+q-1)!!}{(r+q-1)!!}\times\frac{(n+q-1)!!}{(r+q-1)!!}$ rectangular group contained in $D_q^\si$.  This rectangular group of course has $\frac{(m+q-1)!!(n+q-1)!!}{(r+q-1)!!^2}$ idempotents.  Mutiplying this by $e(D_q(\B_r))$ (cf.~\eqref{e:erq}) gives the size of $E(D_q^\si)$, and summing over appropriate $q$ gives the size of $E(P^\si)=E(\Bmns)$.

\begin{thm}\label{t:eBmns}
The number of idempotents of the sandwich semigroup $\Bmns$ is given by
\[
e(\Bmns) = \sum_{0\leq q\leq r \atop q\equiv r\Mod2} \binom rq \frac{(r-q-1)!!(m+q-1)!!(n+q-1)!!}{(r+q-1)!!(2q-1)!!}.  
\]
The $q$th term in the above sum gives the number of idempotents from $D_q^\si$.  \epfres
\end{thm}

\begin{rem}
Since the idempotents of $\Bmns$ are precisely the post-inverses of $\si$, Theorem \ref{t:eBmns} also gives the size of $\Post(\si)$.
\end{rem}

\begin{rem}
Using \eqref{e:erq}, we may give a simpler expression for the number of idempotents in a Brauer monoid $\B_r$ than that given in \cite[Theorem 30]{DEEFHHL1}:
\[
e(\B_r) = \sum_{0\leq q\leq r \atop q\equiv r\Mod2} \binom rq\frac{(r-q-1)!!(r+q-1)!!}{(2q-1)!!}.
\]
This formula is also the $m=n=r$ case of Theorem \ref{t:eBmns}.  An alternative formula for $e(\B_r)$ may be found in \cite[Proposition 4.10]{Larsson}.
\end{rem}

\subsection{MI-domination and the ranks of $\Reg(\Bmns)$ and $\E(\Bmns)$}\label{ss:MIB}

The concept of \emph{MI-domination} (see below for the definition) was introduced in \cite{Sandwich1}.  Among other things, it was shown in \cite[Theorems 4.16 and 4.17]{Sandwich1} that if the regular subsemigroup ${P^a=\Reg(S_{ij}^a)}$ of a sandwich semigroup $S_{ij}^a$ in a regular category $S$ is MI-dominated, then convenient formulae exist for the rank of $P^\si$, and also for the rank and idempotent rank for the idempotent-generated subsemigroup~$\E(P^a)=\E(S_{ij}^a)$.  
In this section we show that sandwich semigroups in the Brauer category $\B$ are MI-dominated, and then use this to deduce the values of the various (idempotent) ranks mentioned above.

Recall that the set $E(T)$ of idempotents of a semigroup $T$ is partially ordered under the relation~$\pre$ defined by $e\pre f \iff e=ef=fe \iff e=fef$.  We have already implicitly referred to this partial order in Section \ref{ss:RBG} when discussing primitive idempotents.  (This order extends in a natural way to the whole semigroup $T$, but we do not need to give the definition \cite{Mitsch1986}.  We also note that $\pre$ is simply the restriction of the ${\leqH}={\leqR}\cap{\leqL}$ preorder to $E(S)$.)
It is easy to check that every mid-identity $u\in\MI(T)$ is maximal in the $\pre$ order: i.e., $u\pre e\implies e=u$ for all $e\in E(T)$.  
As in \cite{Sandwich1}, we say the semigroup $T$ is \emph{MI-dominated} if it is regular and if every idempotent of $T$ is $\pre$-below a mid-identity of $T$.

If $S$ is a partial semigroup, and if $a\in S_{ji}$ is regular, then any preinverse $b\in\Pre(a)$ is clearly a mid-identity of the regular semigroup $P^a=\Reg(S_{ij}^a)$, and indeed of $S_{ij}^a$ itself.  We have already noted (cf.~Proposition \ref{p:Jba2}) that if $S$ is stable and regular, then $\MI(P^a)=V(a)$; it also follows from \cite[Proposition 4.5]{Sandwich1} that $\MI(P^a) =\Hh_b^a$ for any $b\in V(a)$, where the latter denotes the $\gHh^a$-class of $b$.

We now work towards showing that the regular subsemigroup $P^\si=\Reg(\Bmns)$ of a sandwich semigroup $\Bmns$ is MI-dominated.  For the rest of Section \ref{ss:MIB}, we fix $m,n\in\N$ with $m\equiv n\Mod2$, and some $\si\in\B_{nm}$ with $r=\rank(\si)$.  As usual, by Lemma \ref{l:iso1}, we may assume that $\si$ has the form given in \eqref{e:siB}.

Here is the main technical result we need; there is an obvious dual, but we will not state it.

\begin{lemma}\label{l:MI}
If $\al\in P_2^\si$, then $\al=\lam\star_\si\al$ for some $\lam\in\MI(P^\si)$.
\end{lemma}

\pf
Write $\al=\partI{a_1}{a_q}{C_1}{C_s}{b_1}{b_q}{D_1}{D_t}$.  Since $\al\L\si\al$, we have $\codom(\si\al)=\codom(\al)=\{b_1,\ldots,b_q\}$.  Since also $\dom(\si\al)\sub\dom(\si)=[r]$, it follows that the transversals of $\si\al$ are $\{x_1,b_1'\},\ldots,\{x_q,b_q'\}$, for some $x_1,\ldots,x_q\in[r]$.  There are three kinds of nontransversals of $\si\al$:
\bit
\item The lower nontransversals $D_1',\ldots,D_t'$ of $\al$ are still nontransversals of $\si\al$; in fact, these are all the lower nontransversals of $\si\al$.
\item The upper nontransversals $\{r+1,r+2\},\ldots,\{n-1,n\}$ of $\si$ are still nontransversals of $\si\al$.
\item Any other upper nontransversal of $\si\al$ is contained in $[r]$, and there are $k=\frac{r-q}2$ of these.  Suppose these are $\{y_1,z_1\},\ldots,\{y_k,z_k\}$.
\eit
Thus, we have
\begin{equation}\label{e:sial}
\si\al = \Big( 
{ \scriptsize \renewcommand*{\arraystretch}{1}
\begin{array} {\c|\c|\c|\c|\c|\c|\c|\c|\cend}
x_1 \:&\: \cdots \:&\: x_q \:&\: y_1,z_1 \:&\: \cdots \:&\: y_k,z_k \:&\: r+1,r+2 \:&\: \cdots \:&\: n-1,n \\ \cline{4-9}
b_1 \:&\: \cdots \:&\: b_q \:&\: D_1 \:& \multicolumn{4}{c|}{\cdots\cdots\cdots\cdots\cdots\cdots\cdots} &\: D_t
\rule[0mm]{0mm}{2.7mm}
\end{array} 
}
\hspace{-1.5 truemm} \Big).
\end{equation}
We will construct the element $\lam$ in stages.  First:
\ben
\item \label{lam1} $\{a_1,x_1'\},\ldots,\{a_q,x_q'\}$ will all be transversals of $\lam$, and
\item \label{lam2} $\{r+1,r+2\}',\ldots,\{n-1,n\}'$ will all be lower nontransversals of $\lam$.  In fact, these will be all the lower nontransversals.
\een
We still have to construct $r-q=2k$ further transversals, and all the upper nontransversals of $\lam$.  To do the former, consider some $1\leq i\leq k$.  Since $\{y_i,z_i\}$ is a nontransversal of $\si\al$, and since $\{y_i,y_i''\}$ and $\{z_i,z_i''\}$ are both edges in the product graph $\Pi(\si,\al)$, there is a path in $\Pi(\si,\al)$ from $y_i''$ to $z_i''$.  The first edge in this path is of the form $y_i''\to w_i''$ for some upper nontransversal $\{y_i,w_i\}$ of $\al$; renaming if necessary, we may assume that $C_i=\{y_i,w_i\}$.  Now we stipulate that
\ben \setcounter{enumi}{2}
\item \label{lam3} $\{y_1,y_1'\},\ldots,\{y_k,y_k'\}$ and $\{w_1,z_1'\},\ldots,\{w_k,z_k'\}$ will all be transversals of $\lam$.
\een
The elements from $[m]\cup[n]'$ involved in the blocks listed in \ref{lam1}--\ref{lam3} are precisely the elements of $\{a_1,\ldots,a_q\}\cup\{y_1,\ldots,y_k\}\cup\{w_1,\ldots,w_k\}\cup[n]'$.  The remaining elements of $[m]$ are precisely the elements of $C_{k+1},\ldots,C_s$, so we may complete the definition of $\lam$ by specifying that
\ben \setcounter{enumi}{3}
\item \label{lam4} $C_{k+1},\ldots,C_s$ will all be upper nontransversals of $\lam$.
\een
To summarise, we have defined 
\begin{equation}\label{e:lam}
\lam = \Big( 
{ \scriptsize \renewcommand*{\arraystretch}{1}
\begin{array} {\c|\c|\c|\c|\c|\c|\c|\c|\c|\c|\c|\cend}
a_1 \:&\: \cdots \:&\: a_q \:&\: y_1 \:&\: \cdots \:&\: y_k \:&\: w_1 \:&\: \cdots \:&\: w_k \:&\: C_{k+1} \:& \: \cdots \:&\: C_s \\ \cline{10-12}
x_1 \:&\: \cdots \:&\: x_q \:&\: y_1 \:&\: \cdots \:&\: y_k \:&\: z_1 \:&\: \cdots \:&\: z_k \:&\: r+1,r+2 \:& \: \cdots \:&\: n-1,n
\rule[0mm]{0mm}{2.7mm}
\end{array} 
}
\hspace{-1.5 truemm} \Big).
\end{equation}
To complete the proof, we must show that $\al=\lam\si\al$, and that $\lam\in\MI(P^\si)$.  The first of these is clear, in light of \eqref{e:sial} and \eqref{e:lam}.  Since $\MI(P^\si)=V(\si)$, as noted before the statement of the lemma, it remains to show that $\lam=\lam\si\lam$ and $\si=\si\lam\si$; in fact, since $\si\J\lam$ (as $\rank(\lam)=r$) it suffices to do the latter, by Lemma \ref{l:axa}\ref{axa2}.  We begin by identifying the transversals of $\si\lam$.

First note that since $\{y_1,y_1'\},\ldots,\{y_k,y_k'\}$ are all transversals of both $\si$ and $\lam$, these are all transversals of $\si\lam$ as well.

Next consider some $1\leq i\leq q$.  Since $\{x_i,b_i'\}$ is a transversal of $\si\al$, there is a path in the product graph $\Pi(\si,\al)$ of the form
\begin{equation}\label{e:xibi}
x_i \lar\si x_i'' \lar\al u_1'' \lar\si u_2''\lar\al \cdots \lar\si u_{2l}''=a_i'' \lar\al b_i' \quad\text{for some $l\in\N$ and some $u_1,\ldots,u_{2l}\in[m]$.}
\end{equation}
Note that we could have $l=0$, in which case the above path is simply $x_i \lar\si x_i''=a_i''\lar\al b_i'$.  Now, all the edges in \eqref{e:xibi} coming from $\si$ are also of course in the product graph $\Pi(\si,\lam)$.  Next note that the only upper nontransversals of $\al$ that are not blocks of $\lam$ are $C_1,\ldots,C_k$, and these are all involved in components of $\Pi(\si,\al)$ of type \ref{c3}, as enumerated at the beginning of Section \ref{ss:RegB}.  Since the path~\eqref{e:xibi} is of type \ref{c6}, it follows that every edge in \eqref{e:xibi} coming from $\al$, apart from $a_i''\lar\al b_i'$, is also in the product graph $\Pi(\si,\lam)$.  Since $\lam$ contains the transversal $\{a_i,x_i'\}$, it follows that $\Pi(\si,\lam)$ contains the path
\[
x_i \lar\si x_i'' \lar\lam u_1'' \lar\si u_2''\lar\lam \cdots \lar\si u_{2l}''=a_i'' \lar\lam x_i'.
\]
In particular, $\si\lam$ contains the transversals $\{x_1,x_1'\},\ldots,\{x_q,x_q'\}$.

Now consider some $1\leq i\leq k$.  We know that $\{y_i,z_i\}$ is an upper nontransversal of $\si\al$, with ${y_i,z_i\in[r]}$, and that the component of the product graph $\Pi(\si,\al)$ containing $y_i,z_i$ begins with the edges ${y_i\lar\si y_i''\lar\al w_i''}$.  Let the full component be
\begin{equation}\label{e:yizi} 
y_i \lar\si y_i'' \lar\al w_i'' \lar\si v_1''\lar\al v_2'' \lar\si \cdots \lar\al v_{2l}''=z_i'' \lar\si z_i \quad\text{where $l\in\N$ and $v_1,\ldots,v_{2l}\in[m]$.}
\end{equation}
See Figure \ref{f:MI}.
All edges in \eqref{e:yizi} but $y_i'' \lar\al w_i''$ belong to the product graph $\Pi(\si,\lam)$ as well; since $\lam$ contains the transversal $\{w_i,z_i'\}$, it follows that $\Pi(\si,\lam)$ contains the path
\[
z_i' \lar\lam w_i'' \lar\si v_1''\lar\lam v_2'' \lar\si \cdots \lar\lam v_{2l}''=z_i'' \lar\si z_i.
\]
In particular, $\si\lam$ contains the transversals $\{z_1,z_1'\},\ldots,\{z_k,z_k'\}$.

The previous three paragraphs show that $\si\lam$ contains the transversals $\{1,1'\},\ldots,\{r,r'\}$.  Since $\si$ contains these transversals, so too therefore does $\si\lam\si$.  Since $\si\lam\si$ contains all the (upper and lower) nontransversals of $\si$, it follows that $\si=\si\lam\si$.  As noted above, this completes the proof.
\epf

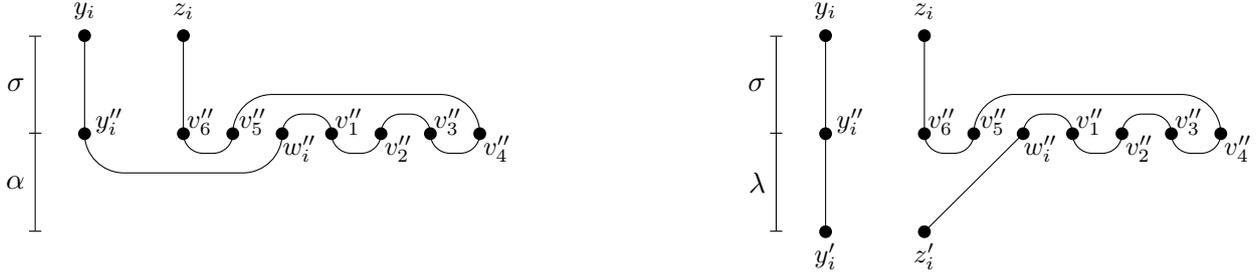
\begin{figure}[ht]
\begin{center}
\begin{tikzpicture}[scale=.65]
\begin{scope}[shift={(0,0)}]	
\uvxs{1,3}{.13}
\lvxs{1,3,4,5,6,7,8,9}{.13}
\darc56
\darc78
\darcx49{.8}
\stline11
\stline33
\draw[|-|] (0,2)--(0,0);
\draw(0.0,1)node[left]{$\si$};
\node () at (1,2.5) {\small $y_i$};
\node () at (3,2.5) {\small $z_i$};
\end{scope}
\begin{scope}[shift={(0,-2)}]	
\uarc34
\uarc67
\uarc89
\uarcx15{.8}
\draw[-|] (0,2)--(0,0);
\draw(0.0,1)node[left]{$\al$};
\node () at (5.35,1.7) {\small $w_i''$};
\node () at (6.35,2.25) {\small $v_1''$};
\node () at (7.35,1.7) {\small $v_2''$};
\node () at (8.35,2.25) {\small $v_3''$};
\node () at (9.35,1.7) {\small $v_4''$};
\node () at (4.4,2.25) {\small $v_5''$};
\node () at (3.35,2.25) {\small $v_6''$};
\node () at (1.5,2.25) {\small $y_i''$};
\end{scope}
\begin{scope}[shift={(15,0)}]	
\uvxs{1,3}{.13}
\lvxs{1,3,4,5,6,7,8,9}{.13}
\darc56
\darc78
\darcx49{.8}
\stline11
\stline33
\draw[|-|] (0,2)--(0,0);
\draw(0.0,1)node[left]{$\si$};
\node () at (1,2.5) {\small $y_i$};
\node () at (3,2.5) {\small $z_i$};
\end{scope}
\begin{scope}[shift={(15,-2)}]	
\lvxs{1,3}{.13}
\uarc34
\uarc67
\uarc89
\stline11
\stline53
\draw[-|] (0,2)--(0,0);
\draw(0.0,1)node[left]{$\lam$};
\node () at (5.35,1.7) {\small $w_i''$};
\node () at (6.35,2.25) {\small $v_1''$};
\node () at (7.35,1.7) {\small $v_2''$};
\node () at (8.35,2.25) {\small $v_3''$};
\node () at (9.35,1.7) {\small $v_4''$};
\node () at (4.4,2.25) {\small $v_5''$};
\node () at (3.35,2.25) {\small $v_6''$};
\node () at (1.5,2.25) {\small $y_i''$};
\node () at (1,-.5) {\small $y_i'$};
\node () at (3,-.5) {\small $z_i'$};
\end{scope}
\end{tikzpicture}
\caption{Left: a component of type \eqref{e:yizi} in the product graph $\Pi(\si,\al)$.  Right: the corresponding two components of $\Pi(\si,\lam)$.}
\label{f:MI}
\end{center}
\end{figure}

We may now deduce the following:

\begin{prop}\label{p:MI}
The semigroup $P^\si=\Reg(\Bmns)$ is MI-dominated.
\end{prop}

\pf
Let $\al\in E(P^\si)$.  Since $\al\in P^\si=P_1^\si\cap P_2^\si$, Lemma \ref{l:MI} and its dual tell us that $\al=\lam\star_\si\al\star_\si\rho$ for some $\lam,\rho\in \MI(P^\si)$.  Put $\ve=\lam\star_\si\rho$, noting that $\ve\in\MI(P^\si)$, since $\MI(P^\si)$ is a subsemigroup.  Since $\lam,\rho$ are mid-identities, we have $\al=\lam\star_\si\al\star_\si\rho=\lam\star_\si\rho\star_\si\al\star_\si\lam\star_\si\rho=\ve\star_\si\al\star_\si\ve$, which shows that $\al\pre\ve$.
\epf

\begin{rem}
If $\KK$ is any of the other diagram categories studied in this paper, then the regular subsemigroup $\Reg(\Kmns)$ is not MI-dominated in general, as may be verified by GAP \cite{GAP}.  For example, consider the partition $\si=\custpartn{1,2,3}{1,2,3}{\stline11\stline32\stline33\darc23}$ from $\P_3=\P_{3,3}$.  Using Proposition \ref{p:Jba2}\ref{Jba22} and \ref{Jba24}, GAP calculates:
\[
\MI(\P_3^\si) = V(\si) = \left\{
\custpartn{1,2,3}{1,2,3}{\stline11\stline23},
\custpartn{1,2,3}{1,2,3}{\stline11\stline12\stline23\darc12},
\custpartn{1,2,3}{1,2,3}{\stline11\stline22\stline23\darc23},
\custpartn{1,2,3}{1,2,3}{\stline11\stline33},
\custpartn{1,2,3}{1,2,3}{\stline11\stline12\stline33\darc12},
\custpartn{1,2,3}{1,2,3}{\stline11\stline32\stline33\darc23},
\custpartn{1,2,3}{1,2,3}{\stline11\stline22\stline33\darc23\uarc23},
\custpartn{1,2,3}{1,2,3}{\stline11\stline12\stline23\stline33\darc12\uarc23},
\custpartn{1,2,3}{1,2,3}{\stline11\stline23\stline33\uarc23}
\right\}.
\]
It also tells us that $E(\P_3^\si)$ has size $99$, but the set $\set{\ve\star_\si\al\star_\si\ve}{\ve\in\MI(\P_3^\si),\ \al\in E(\P_3^\si)}$ has size~$83$.  From this, it quickly follows that $\Reg(\P_3^\si)$ is not MI-dominated.  As an explicit example, consider~$\al=\custpartn{1,2,3}{1,2,3}{\stline33\darc12\uarc12}\in\P_3$.  Then one may check that $\al=\al\si\al$, so that $\al\in E(\P_3^\si)$.  But $\al$ is not~$\leqR$-below any of the above mid-identities (cf.~Theorem \ref{t:G}\ref{G1}), so it follows that $\al$ is not $\pre$-below any mid-identity.
\end{rem}

Armed with Proposition \ref{p:MI}, we may quickly deduce information about the ranks of the regular subsemigroup $P^\si=\Reg(\Bmns)$ and the idempotent-generated subsemigroup $\E(\Bmns)$, as well as the idempotent rank of the latter.  (Recall that the \emph{idempotent rank} of an idempotent-generated semigroup~$T$, denoted $\idrank(T)$, is the smallest size of a generating set for $T$ consisting of idempotents.)

For the next proof, we also need the concept of \emph{relative rank} \cite{HRH1998}.  If $T$ is a semigroup, and if $A\sub T$, then the \emph{relative rank of $T$ modulo $A$}, denoted $\relrank TA$, is the smallest size of a subset $B\sub T$ such that~$T$ is generated by $A\cup B$.  
It follows immediately from \cite[Lemma 2.1]{Mazorchuk2002} that
\begin{equation}\label{e:BrSr}
\relrank{\B_r}{\S_r}=\begin{cases}
1 &\text{if $r\geq2$}\\
0 &\text{if $r\leq1$.}
\end{cases}
\end{equation}
We will also need to refer to $\Hh_{\si^*}^\si$, the $\gHh^\si$-class of $\si^*$.  Note that the $q=r$ case of Theorem \ref{t:RegB}\ref{RB3} says that $D_r^\si$ is a single $\gHh^\si$-class; since $\si^*\in D_r^\si$, it follows in fact that $\Hh_{\si^*}^\si=D_r^\si$.  It follows from the same theorem that
\begin{equation}\label{e:DrRL}
|\Hh_{\si^*}^\si/{\R^\si}| = \tfrac{(m+r-1)!!}{(2r-1)!!} \AND |\Hh_{\si^*}^\si/{\L^\si}| = \tfrac{(n+r-1)!!}{(2r-1)!!}.
\end{equation}

The next result gives the rank of the regular semigroup $P^\si=\Reg(\Bmns)$.  If $r=m=n$, then $P^\si=\B_n$.  We have already noted that $\rank(\B_n)=3$ for $n\geq3$, so we will exclude the $r=m=n$ case from the statement.  Also, since $\Bmns$ is anti-isomorphic to $\B_{nm}^{\si^*}$, it suffices to assume that $m\geq n$.

\begin{thm}\label{t:rankP}
If $m\geq n$, and if $r=m=n$ does not hold, then the rank of the regular semigroup~${P^\si=\Reg(\Bmns)}$ is given by
\[
\rank(P^\si) = \frac{(m+r-1)!!}{(2r-1)!!} + \begin{cases}
1 &\text{if $r\geq2$}\\
0 &\text{if $r\leq1$.}
\end{cases}
\]
\end{thm}

\pf
By \cite[Theorem 4.16]{Sandwich1}, and since $P^\si$ is MI-dominated (cf.~Proposition \ref{p:MI}), we have
\[
\rank(P^\si) = \relrank{\B_r}{\S_r} + \max\big(|\Hh_{\si^*}^\si/{\R^\si}|,|\Hh_{\si^*}^\si/{\L^\si}|,\rank(\S_r)\big).
\]
The result now follows quickly from \eqref{e:BrSr} and \eqref{e:DrRL}, along with the fact \cite{Moo} that
\[
\rank(\S_r)=\begin{cases}
2&\text{if $r\geq3$}\\
1&\text{if $r\leq2$.}
\end{cases}
\qedhere
\]
\epf

In the next statement we again assume by symmetry that $m\geq n$, but we do not need to exclude the possibility that $r=m=n$.

\begin{thm}\label{t:rankEB}
If $m\geq n$, then the rank and idempotent rank of the idempotent-generated semigroup~$\E(\Bmns)$ are given by
\[
\rank(\E(\Bmns)) = \idrank(\E(\Bmns)) =\frac{(m+r-1)!! }{(2r-1)!!} +  \binom r2.
\]
\end{thm}

\pf
By \cite[Theorem 4.17]{Sandwich1}, and since $P^\si$ is MI-dominated (cf.~Proposition \ref{p:MI}), we have
\begin{align*}
\rank(\E(\Bmns)) &= \rank(\E(\B_r)) + \max\big(|\Hh_{\si^*}^\si/{\R^\si}|,|\Hh_{\si^*}^\si/{\L^\si}|\big) - 1, \\[2mm]
\idrank(\E(\Bmns)) &= \idrank(\E(\B_r)) + \max\big(|\Hh_{\si^*}^\si/{\R^\si}|,|\Hh_{\si^*}^\si/{\L^\si}|\big) - 1.
\end{align*}
The result now follows from \eqref{e:DrRL}, along with the fact \cite[Proposition 2]{MM2007} that
\[
\rank(\E(\B_r))=\idrank(\E(\B_r))=1+\tbinom r2.  \qedhere
\]
\epf

\subsection{Ideals of $\Reg(\Bmns)$}\label{ss:ideals}

We now turn our attention to the ideals of the regular subsemigroup $P^\si=\Reg(\Bmns)$.  These are of the form $I_q=I_q(P^\si)$, for each $0\leq q\leq r$ with $q\equiv r\Mod2$, in the notation of Section \ref{ss:idealsK}.  Theorem \ref{t:IqK} shows that any proper ideal $I_q$ ($q<r$) is idempotent-generated.  In this section we wish to improve this result, and show that $I_q$ is generated by the idempotents from $E(D_q^\si)$.  We will also calculate the rank and idempotent rank of each proper ideal $I_q$.
Before we state the result (Theorem \ref{t:IqB}), we require two results from \cite{Sandwich1}:

\begin{lemma}[{cf.~\cite[Proposition 4.3]{Sandwich1}}]\label{l:eTe}
If $T$ is a regular semigroup, then $T$ is MI-dominated if and only if $\displaystyle{T=\bigcup_{e\in\MI(T)} eTe}$. \epfres
\end{lemma}

The next result makes use of the surmorphism
\[
\Phi: P^a=\Reg(S_{ij}^a)\to abS_jab:x\mt axab,
\]
defined for a regular partial semigroup $S$, and for $a\in S_{ji}$ and $b\in V(a)$; cf.~\eqref{e:Phi}.

\begin{lemma}[{cf.~\cite[Proposition 4.8]{Sandwich1}}]\label{l:ePe}
For any $e\in V(a)$, the restriction of $\Phi$ to $e\star_a P^a\star_ae$ is an isomorphism onto $abS_jab$. \epfres
\end{lemma}

\begin{rem}\label{r:ePe}
In the case that $S=\B$, and using the isomorphism $\si\si^*\B_n\si\si^*\to\B_r:\al\mt\al^\natural$ from~\eqref{e:Kr}, it follows that for any $\ve\in V(\si)$, the map
\[
\ve\star_\si P^\si\star_\si\ve \to \B_r: \al\mt(\si\al\si\si^*)^\natural = \tau\si\al\tau^*
\]
is an isomorphism.
\end{rem}

We also require the following basic fact about generating sets of stable semigroups (as usual, we omit the dual statement):

\begin{lemma}\label{l:rankJR}
If $J$ is a maximal $\J$-class of a semigroup $T$, and if every element of $J$ is stable, then $\rank(T)\geq|J/{\R}|$.
\end{lemma}

\pf
Suppose $T=\la X\ra$.  The proof will be complete if we can show that $X$ contains at least one element from each $\R$-class of $J$.  To do so, let $z\in J$ be arbitrary.  Then $z=x_1\cdots x_k$ for some $x_1,\ldots,x_k\in X$.  Since $z = x_1\cdots x_k \leqJ x_1$, maximality of $J=J_z$ gives $z\J x_1$.  But then $x_1\J z = x_1(x_2\cdots x_k)$, so it follows from stability that $x_1\R x_1(x_2\cdots x_k)=z$.
\epf

\begin{rem}
Lemma \ref{l:rankJR} could also be proved using the Rees Theorem \cite[Theorem 3.2.3]{Howie} and results about generating sets for completely 0-simple semigroups; see for example \cite[Theorem 7]{Gray2014} or \cite[Lemma 3.1]{Ruskuc1994}.
\end{rem}

We now return our attention to the ideals of $P^\si$.  

\begin{lemma}\label{l:Dp}
For any $0\leq p\leq r-4$ with $p\equiv r\Mod2$, we have $D_p^\si\sub\la D_{p+2}^\si\ra$.
\end{lemma}

\pf
Let $\al\in D_p^\si$.  Since $\ol\al\in D_p(\B_r)\sub I_{p+2}(\B_r)=\big\la D_{p+2}(\B_r)\big\ra$, using \cite[Theorem 8.4]{EG2017} in the last step, we have $\ol\al=\be_1\cdots\be_k$ for some $\be_1,\ldots,\be_k\in D_{p+2}(\B_r)$.  
For each~$i$, we have $\be_i=\ol\ga_i$ for some $\ga_i\in D_{p+2}^\si$.

Since $P^\si$ is MI-dominated (cf.~Proposition \ref{p:MI}), it follows from Lemma \ref{l:eTe} that $\al\in\ve\star_\si P^\si\star_\si\ve$ for some $\ve\in\MI(P^\si)$.  Now put $\de=\ve\star_\si(\ga_1\star_\si\cdots\star_\si\ga_k)\star_\si\ve$.  Since $\MI(P^\si)=E(J_{\si^*}^\si)=E(D_r^\si)$, by Proposition \ref{p:Jba2}\ref{Jba24}, it follows that $\ol\ve\in E(D_r(\B_r))=\{\id_r\}$, and so $\ol\de=\id_r(\ol\ga_1\cdots\ol\ga_k)\id_r=\ol\al$.  But $\al$ and $\de$ both belong to $\ve\star_\si P^\si\star_\si\ve$, so by Lemma \ref{l:ePe} (cf.~Remark \ref{r:ePe}), it follows that $\al=\de$.  Thus,
\[
\al = (\ve\star_\si\ga_1)\star_\si\ga_2\star_\si\cdots\star_\si\ga_{k-1}\star_\si(\ga_k\star_\si\ve).
\]
Since $\ga_2,\ldots,\ga_{k-1}\in D_{p+2}^\si$, the proof will be complete if we can show that $\ve\star_\si\ga_1$ and $\ga_k\star_\si\ve$ also belong to $D_{p+2}^\si$.  But this follows from 
$
p+2 = \rank(\ol\ga_1) = \rank(\id_r\ol\ga_1) = \rank(\ol{\ve\star_\si\ga_1}) = \rank(\ve\star_\si\ga_1),
$
and the analogous calculation for $\ga_k\star_\si\ve$.
\epf

We are now ready to prove our final main result.  As usual, we may assume that $m\geq n$.

\begin{thm}\label{t:IqB}
If $m\geq n$, then for any $0\leq q<r$ with $q\equiv r\Mod2$, we have $I_q = \big\la E(D_q^\si)\big\ra$.  Moreover, 
\[
\rank(I_q) = \idrank(I_q) = \binom rq\frac{(r-q-1)!!(m+q-1)!!}{(r+q-1)!!}.
\]
\end{thm}

\pf
We first note that $I_q = \la D_q^\si\ra$.  Indeed, it follows from Lemma \ref{l:Dp} and a simple descending induction that ${D_p^\si \sub \la D_q^\si\ra}$ for all $0\leq p\leq q$ with $p\equiv q\Mod2$.

Next note that the top $\J^\si$-class of $I_q$ is $D_q^\si$, and this has $\ka(m,r,q)=\binom rq \frac{(r-q-1)!!(m+q-1)!!}{(r+q-1)!!}$ $\R^\si$-classes; cf.~Theorem \ref{t:RegB}\ref{RB1}.  Since $I_q$ is stable, it follows from Lemma \ref{l:rankJR} that
\[
\rank(I_q) \geq |D_q^\si/{\R^\si}| = \binom rq \frac{(r-q-1)!!(m+q-1)!!}{(r+q-1)!!}.
\]
Thus, since $\idrank(T)\geq\rank(T)$ for any idempotent-generated semigroup $T$, it remains to construct a generating set of the desired size, consisting of idempotents.

By \cite[Theorem 8.4]{EG2017}, we have $I_q(\B_r) = \big\la E(D_q(\B_r)) \big\ra$ and $\idrank(I_q(\B_r)) = \binom rq(r-q-1)!!$, so we may fix some $\Om\sub E(D_q(\B_r))$ such that $I_q(\B_r)=\la\Om\ra$ and $|\Om|=\binom rq(r-q-1)!!$.  Let $\Ga=\Om\Psi^{-1}$.  First we claim that $I_q=\la\Ga\ra$.

To prove this claim, since we already know that $I_q=\la D_q^\si\ra$, it suffices to show that $D_q^\si\sub\la\Ga\ra$.  With this in mind, let $\al\in D_q^\si$ be arbitrary.  Since $\ol\al\in D_q(\B_r)$, we have $\ol\al=\be_1\cdots\be_k$ for some $\be_1,\ldots,\be_k\in\Om$.  We then follow the proof of Theorem \ref{t:IqK}, and conclude that
\[
\al = \de \star_\si(\ga_1\star_\si\cdots\star_\si\ga_k)\star_\si\ve,
\]
where $\ol\ga_i=\be_i$ for each $i$, and where $\de\in\Hh_{\ga_1}^\si$ and $\ve\in\Hh_{\ga_k}^\si$ are both idempotents.  By definition, we have $\ga_i\in\Om\Psi^{-1}=\Ga$ for all $i$ (since each $\be_i\in\Om$).  Since $\de\in\Hh_{\ga_1}^\si$, we have $\ol\de\H\ol\ga_1$; but since $\ol\de$ and~$\ol\ga_1$ are both idempotents, it follows that in fact $\ol\de=\ol\ga_1=\be_1\in\Om$, which means that $\de\in\Om\Psi^{-1}=\Ga$.  Similarly, $\ve\in\Ga$.  It follows that $\al\in\la\Ga\ra$.  This completes the proof of the claim that $I_q=\la\Ga\ra$.

For each $\al\in\Om$, let $\Ga_\al=\al\Psi^{-1}$, so that $\Ga=\bigcup_{\al\in\Om}\Ga_\al$.  By Theorem \ref{t:RegB}\ref{RB5} and \eqref{e:EK}, it follows that $\Ga_\al = E(H_\al\Psi^{-1})$ is a $\frac{(m+q-1)!!}{(r+q-1)!!}\times\frac{(n+q-1)!!}{(r+q-1)!!}$ rectangular band.
It follows from results of Ru\v{s}kuc \cite{Ruskuc1994} (see also \cite[Proposition 4.11]{Sandwich1}) that the rank of a $\rho\times\lam$ rectangular band is equal to $\max(\rho,\lam)$.  Thus, for each $\al\in\Om$, we may fix a subset $\Xi_\al$ of $\Ga_\al$ such that
\[
\Ga_\al = \la\Xi_\al\ra \AND |\Xi_\al| = \max\left(\frac{(m+q-1)!!}{(r+q-1)!!},\frac{(n+q-1)!!}{(r+q-1)!!}\right) = \frac{(m+q-1)!!}{(r+q-1)!!}.
\]
Then with $\Xi=\bigcup_{\al\in\Om}\Xi_\al$, we have $\Ga = \bigcup_{\al\in\Om}\Ga_\al = \bigcup_{\al\in\Om}\la\Xi_\al\ra \sub \la\Xi\ra$, so that $I_q=\la\Ga\ra\sub\la\Xi\ra$, whence $I_q=\la\Xi\ra$.  But $\Xi$ consists of idempotents, and has the desired size, so the proof is complete.
\epf

\begin{rem}
We have already noted that while ideals of $P^\si=\Reg(\Kmns)$ in other diagram categories are often generated by their idempotents, they are generally not generated by the idempotents in their top $\D^\si$-class; cf.~Theorem \ref{t:IqK} and Remark \ref{r:IqK}.  It would therefore be a very interesting problem to try and calculate the (idempotent) ranks of these ideals in general.
\end{rem}

\section*{Acknowledgements}

The first and second authors are supported by Grant Nos.~174018 and 174019, respectively, of the Ministry of Education, Science, and Technological Development of the Republic of Serbia.
The third author thanks the University of Novi Sad for its hospitality during his visit in 2016.

\footnotesize
\def\bibspacing{-1.1pt}
\bibliography{biblio}
\bibliographystyle{abbrv}

\end{document}